\def\simplelatex{\iftrue}
\documentclass[11pt]{article}
\usepackage[francais]{babel}
\simplelatex%
\else%
\usepackage[latin1]{inputenc}
\usepackage[T1]{fontenc}
\fi%
\usepackage{amsmath}
\usepackage{amsfonts}
\usepackage{amssymb}
\usepackage{amsthm}
\simplelatex%
\let\mathscr=\mathcal
\else%
\usepackage{mathrsfs}
\fi%
\simplelatex%
\let\url=\texttt
\else%
\usepackage{url}
\fi%
\newtheorem{theo}{Th\'{e}or\`{e}me}[section]

\newtheorem{prop}[theo]{Proposition}
\newtheorem{lemm}[theo]{Lemme}
\newtheorem{coro}[theo]{Corollaire}
\newtheorem{rema}[theo]{Remarque}
\newtheorem{Defi}[theo]{D\'{e}finition}

\newtheorem{conj}[theo]{Conjecture}
\newtheorem{question}[theo]{Question}

\newcommand{\Spec}{\mathop\mathrm{Spec}\nolimits}

\newcommand{\Pic}{\mathop\mathrm{Pic}\nolimits}

\def\oi{\hskip1mm {\buildrel \simeq \over \rightarrow} \hskip1mm}
\def\loi{\hskip1mm {\buildrel \simeq \over \leftarrow} \hskip1mm}
\def\ovF{{\overline F}}
\def\ovV{{\overline V}}
\def\K{{\mathcal K}}
\def\et{\rm{\acute et}}
\def\O{  {\mathcal O}  }
\def\Ker{{\rm Ker}}
\def\Coker{{\rm Coker}}
\def\Im{{\rm Im}}
\def\Spec{{\rm Spec}}
\def\dim{{\rm dim}}
\def\deg{{\rm deg}}

\def\cqfd{\qed}
\def\H{{\mathcal H}}

\newcommand{\Q}{{\mathbb Q}}
\newcommand{\C}{{\mathbb C}}

\newcommand{\Z}{{\mathbb Z}}
\newcommand{\N}{{\mathbb N}}

\newcommand{\A}{{\mathbb A}}
 \def\et{{\acute et}}

\simplelatex%
\let\bimu=\mu
\let\mathbi=\mathbf
\else%
\DeclareFontFamily{OML}{cmmib}{\skewchar\font127 }
\DeclareFontShape{OML}{cmmib}{m}{it}%
       {<5><6><7><8><9>gen*cmmib%
        <10><10.95>cmmib10%
        <12><14.4><17.28><20.74><24.88>cmmib12%
        }{}
\DeclareSymbolFont{biletters}{OML}{cmmib}{m}{it}
\DeclareSymbolFontAlphabet{\mathbi}{biletters}
\DeclareMathSymbol{\bimu}{\mathord}{biletters}{"16}
\fi%
\simplelatex%

\fi%

\begin{document}

\title
{Cohomologie non ramifi\'{e}e et conjecture
de Hodge enti\`{e}re}
 \author{Jean-Louis Colliot-Th\'{e}l\`{e}ne et Claire Voisin}
 \date{\footnotesize    10 septembre 2011}

 \maketitle

Abstract : Building upon the Bloch--Kato conjecture in Milnor K-theory,
we relate the third unramified cohomology group with $\Q/\Z$ coefficients
with a group which measures the failure of the integral Hodge conjecture
in degree 4.
As a first  consequence, a geometric theorem of the second-named author  implies
that the third unramified cohomology group with $\Q/\Z$ coefficients
vanishes on all  uniruled threefolds.
As a second consequence, a 1989 example by Ojanguren and the first named author
implies that the integral Hodge conjecture in degree 4  fails for unirational varieties
of dimension at least 6.
For certain classes of threefolds fibered over a curve, we establish  a relation
between the integral Hodge conjecture  and the computation of
the index of the generic fibre.

\medskip

R\'{e}sum\'{e} : En nous appuyant sur la conjecture de Bloch--Kato en K-th\'{e}orie de Milnor,
nous \'{e}tablissons un lien g\'{e}n\'{e}ral entre le d\'{e}faut de la conjecture
de Hodge enti\`{e}re pour la cohomologie de degr\'{e} 4 et  le troisi\`{e}me
groupe de cohomologie non  ramifi\'{e}e \`{a} coefficients $\Q/\Z$.
Ceci  permet  de montrer que
sur un solide \footnote{en anglais, threefold} 
unir\'{e}gl\'{e}  le troisi\`{e}me
groupe de cohomologie non  ramifi\'{e}e \`{a} coefficients $\Q/\Z$
s'annule, ce que la K-th\'{e}orie alg\'{e}brique ne permet d'obtenir
que dans certains cas.
Ceci permet \`{a} l'inverse de d\'{e}duire d'exemples ayant leur
source en  K-th\'{e}orie que la conjecture de Hodge enti\`{e}re
pour la cohomologie de degr\'{e} 4 peut \^{e}tre en d\'{e}faut
pour les vari\'{e}t\'{e}s rationnellement connexes.
Pour certaines  familles \`{a} un param\`{e}tre de surfaces, on \'{e}tablit
un lien entre  la  conjecture de Hodge enti\`{e}re et l'indice de la fibre g\'{e}n\'{e}rique.

\section{Introduction}

Soit $X$ une vari\'{e}t\'{e} projective, lisse, connexe,  sur le corps $\C$ des complexes, de dimension $d$.
Soit $\Z(1)=\Z(2\pi i) \subset \C$ et, pour $r>0$, soit $\Z(r)=\Z(1)^{\otimes r}$.
Pour $i\geq 0$, soit $CH^{i}(X)$ le groupe des cycles  alg\'{e}briques de codimension $i$
sur $X$ modulo l'\'{e}quivalence rationnelle.

Pour tout entier $i \geq 1$ on sait d\'{e}finir (cf. \cite[IV.11.1.2]{voisinlivre})
des applications  classe de cycles
$$ CH^{i}(X) \to H^{2i}(X(\C),\Z(i))$$
qui g\'{e}n\'{e}ralisent l'application $\Pic(X) \to H^2(X(\C),\Z(1))$ d\'{e}duite de la suite
de l'exponentielle.
Via les isomorphismes de dualit\'{e} de Poincar\'{e}
entre cohomologie enti\`{e}re et
homologie enti\`{e}re
$$ H^{2d-2j}(X(\C),\Z(d-j)) \simeq H_{2j}(X(\C),\Z),$$
ceci d\'{e}finit aussi
des applications   classe de cycles homologiques
$$CH_{j}(X) \to H_{2j}(X(\C),\Z).$$

L'image de $CH^{i}(X)$  dans  $ H^{2i}(X(\C),\Z(i))$ sera not\'{e}e $H^{2i}_{alg}(X(\C),\Z(i))$.

Dans la cohomologie
$H^{2i}(X(\C),\Q(i))$,
on a le sous-groupe
 $Hdg^{2i}(X,\Q)$ des classes
de Hodge : ce sont celles dont l'image dans $H^{2i}(X(\C),\C (i))$ est de type $(i,i)$
pour la d\'{e}composition de Hodge. On a une inclusion
$H^{2i}_{alg}(X(\C),\Z(i))\otimes \Q \subset Hdg^{2i}(X,\Q).$
La conjecture de Hodge affirme que cette inclusion est une \'{e}galit\'{e}.
C'est connu pour $i=1$ et pour $i=d-1$, o\`u $d=\dim(X)$.

On d\'{e}finit le groupe des classes de Hodge enti\`{e}res  $Hdg^{2i}(X,\Z) \subset
H^{2i}(X(\C),\Z(i))$ comme l'image r\'{e}ciproque de $Hdg^{2i}(X,\Q)$
dans $H^{2i}(X(\C),\Z(i))$ par l'application naturelle de la cohomologie enti\`{e}re
vers la cohomologie rationnelle.

 De ce qui pr\'{e}c\`{e}de on tire l'inclusion
$$H^{2i}_{alg}(X(\C),\Z(i)) \subset Hdg^{2i}(X,\Z).$$
On note  $Z^{2i}(X): = Hdg^{2i}(X,\Z)/H^{2i}_{alg}(X(\C),\Z(i))$.
En utilisant la dualit\'{e} de Poincar\'{e}, on d\'{e}finit les groupes $Z_{2j}(X)$.

On a  l'inclusion $Z^{2i}(X)  \subset H^{2i} (X(\C),\Z(i))/H^{2i}_{alg}(X(\C),\Z(i))$,
les sous-groupes de torsion de chacun de ces deux groupes co\"{\i}ncident,
et sont conjecturalement \'{e}gaux \`{a} $Z^{2i}(X)$
(cet \'{e}nonc\'{e} est \'{e}quivalent \`{a} la conjecture de Hodge en degr\'{e} $2i$).
En particulier chaque groupe $Z^{2i}(X)$ est conjecturalement un groupe fini.
Ceci est connu pour $i=1$ et $i=d-1$.
 Pour $i=2$,
c'est connu lorsque le groupe de Chow des z\'{e}ro-cycles $CH_{0}(X)$
est support\'{e} sur un ferm\'{e} alg\'{e}brique de $X$ de dimension
 $\leq 3$ (th\'{e}or\`{e}me de Bloch--Srinivas \cite{blochsrinivas}).

Une version enti\`{e}re de la conjecture de Hodge dirait que les groupes $Z^{2i}(X)$
sont nuls.  On sait que c'est le cas pour $i=1$, mais d\'{e}j\`{a} pour $i=2$,
il existe des contre-exemples (Atiyah-Hirzebruch~\cite{atiyahhirzebruch}, Koll\'ar \cite{kollar}).

On cherche \`{a} d\'{e}terminer des classes de vari\'{e}t\'{e}s pour lesquelles,
pour certaines valeurs de $i$, ce   groupe est toujours nul.

On s'int\'{e}resse particuli\`{e}rement aux cas $i=2$ et $i=d-1$, pour la raison suivante :
Le groupe
 $Z^4(X)$
et le  groupe  fini  $Z_{2}(X) \simeq Z^{2d-2}(X)$
sont des invariants birationnels des vari\'{e}t\'{e}s projectives et lisses. La d\'{e}monstration de
ce fait dans
 \cite[Lemme 15]{voisinjapjmath}  utilise le th\'{e}or\`{e}me d'Hironaka, le comportement
 de la cohomologie et des groupes de Chow par
\'{e}clatement de sous-vari\'{e}t\'{e}s lisses, et, pour $Z^4(X)$,  le th\'{e}or\`{e}me de Lefschetz sur les classes
de type $(1,1)$.

L'\'{e}tude de la cohomologie \'{e}tale des vari\'{e}t\'{e}s \`{a} coefficients finis
(qui pour les vari\'{e}t\'{e}s sur $\C$  \'{e}quivaut \`{a} l'\'{e}tude de la cohomologie
de Betti \`{a} coefficients finis) a men\'{e} \`{a} introduire d'autres invariants birationnels
permettant en particulier de d\'{e}tecter la non rationalit\'{e} de certaines vari\'{e}t\'{e}s unirationnelles.
Les invariants en question sont les groupes de cohomologie non ramifi\'{e}e $H^{i}_{nr}(X,\Z/n)$
introduits dans \cite{colliotojanguren}.
La preuve de leur invariance birationnelle
utilise la conjecture de Gersten en cohomologie \'{e}tale, \'{e}tablie par Bloch et Ogus \cite{blochogus}.

L'observation de cette double invariance birationnelle nous a amen\'{e}s \`{a} nous demander
s'il existe un lien entre la torsion de $Z^4(X)$ et les groupes $H^3_{nr}(X,\Z/n)$.
Il est aussi notoire que   tant  les
 groupes $Z^{2i}(X)$ que les groupes de cohomologie non ramifi\'{e}e
sont difficiles \`{a} calculer, et que leur comportement  \og  en famille \fg \, est myst\'{e}rieux.

Des r\'{e}sultats  sur $Z^{4}(X)$ ont \'{e}t\'{e} obtenus
par des m\'{e}thodes de g\'{e}om\'{e}trie complexe (\cite{voisinuniruled}, \cite{voisinjapjmath}).
 Des r\'{e}sultats
sur les groupes $H^{3}_{nr}(X)$ \`{a} coefficients de torsion ont \'{e}t\'{e} obtenus par des m\'{e}thodes de $K$-th\'{e}orie alg\'{e}brique
(\cite{colliotojanguren}, \cite{kahn}).

\medskip

La th\'{e}orie de Bloch--Ogus, combin\'{e}e \`{a} la conjecture de Bloch--Kato,
  maintenant un th\'{e}or\`{e}me de   Voevodsky
\cite{voevodskymodl} et Rost,
et \`{a} un argument de Bloch et Srinivas,
permet d'\'{e}tablir un lien entre les deux types d'invariants.
Pour les vari\'{e}t\'{e}s de dimension 3, ce
 lien avait \'{e}t\'{e} remarqu\'{e} en 1992 par Barbieri-Viale \cite{barbieriviale}.
  En prenant appui sur le cas g\'{e}n\'{e}ral de la conjecture de Bloch--Kato,
nous montrons au \S 3 que le lien existe en toute dimension.
Le th\'{e}or\`{e}me g\'{e}n\'{e}ral est le th\'{e}or\`{e}me \ref{principal1}.
Citons-en ici une cons\'{e}quence  (Th\'{e}or\`{e}me \ref{blochsrinivasgriffnul}) :

\begin{theo}\label{theointro10avril}
Soit $X$ une vari\'{e}t\'{e},
projective et lisse  sur les complexes, dont le groupe de Chow des z\'{e}ro-cycles
$CH_{0}(X)$ est support\'{e} sur une surface. On a   un isomorphisme de groupes finis
$$ H^3_{nr}(X,\Q/\Z(2)) \oi  Z^4(X),$$
o\`u le premier groupe est l'union de ses sous-groupes $H^3_{nr}(X,\mu_{n}^{\otimes 2})$.
\end{theo}

Cela nous permet de traduire et comparer les r\'{e}sultats obtenus
par la g\'{e}om\'{e}trie alg\'{e}brique et par la $K$-th\'{e}orie alg\'{e}brique.

La combinaison du th\'{e}or\`{e}me  \ref{theointro10avril} et
d'un th\'{e}or\`{e}me obtenu par des m\'{e}thodes
de g\'{e}om\'{e}trie complexe \cite{voisinuniruled} donne :

\begin{theo}\label{illustr1}
  Soit $X$ un solide projectif et lisse  sur les complexes. Si $X$ est
  unir\'{e}gl\'{e}, alors $H^3_{nr}(X, \Q/\Z(2))=0.$
\end{theo}
Pour les solides fibr\'{e}s en coniques sur une surface, ce r\'{e}sultat avait \'{e}t\'{e} obtenu
via la $K$-th\'{e}orie alg\'{e}brique en 1989 \cite{parimala}.

\medskip

La combinaison du th\'{e}or\`{e}me \ref{theointro10avril}
et d'un contre-exemple \`{a} la rationalit\'{e} obtenu par
des m\'{e}thodes de $K$-th\'{e}orie alg\'{e}brique \cite{colliotojanguren} donne :

\begin{theo}\label{illustr2}
Il existe une vari\'{e}t\'{e} $X$ projective, lisse, connexe sur $\C$,
unirationnelle, de dimension 6, avec $Z^4(X)=Z^4(X)_{\rm tors} \neq 0$.
\end{theo}

Ceci r\'{e}pond n\'{e}gativement \`{a} une question soulev\'{e}e dans
\cite{voisinjapjmath} : Si $X$ est une vari\'{e}t\'{e} rationnellement connexe,
le groupe $Z^4(X)=Z^{4}(X)_{\rm tors}$ est-il nul ?

 \medskip

Donnons des indications plus d\'{e}taill\'{e}es sur la structure de l'article. Le paragraphe~2
est consacr\'{e} \`{a} des rappels sur la th\'{e}orie de Bloch et Ogus et sur la conjecture
de Bloch--Kato. Au paragraphe~3, on commence par
\'{e}tendre un argument de Bloch--Srinivas en dimension quelconque.
Ceci nous permet ensuite de montrer que le groupe fini $Z^4(X)_{\rm tors}$ est
le quotient du groupe de cohomologie non   ramifi\'{e}e $H^3_{nr}(X,\Q/\Z(2))$
par son sous-groupe divisible maximal,  r\'{e}sultat qui implique le th\'{e}or\`{e}me
\ref{theointro10avril}.
Nous \'{e}tablissons un r\'{e}sultat
analogue pour le groupe fini $Z^{2d-2}(X)$ (th\'{e}or\`{e}me \ref{principal3}).
Il y a des liens entre le sous-groupe divisible maximal de $H^3_{nr}(X,\Q/\Z(2))$,
le groupe de Griffiths en codimension 2 et le groupe de cohomologie coh\'{e}rente
$H^3(X,\O_{X})$. Ces liens, en partie conjecturaux, sont discut\'{e}s au para\-graphe~4.
Au para\-graphe 5, on passe en revue les exemples, d'origines tr\`{e}s diverses,
de vari\'{e}t\'{e}s $X$ pour lesquelles divers auteurs ont montr\'{e} que soit le
groupe $Z^4(X)_{\rm tors}$, soit le groupe $H^3_{nr}(X,\Z/n)$ est non nul.
{\it On construit un exemple de solide $X$ dont conjecturalement
le groupe de Chow $CH_{0}(X)$ est r\'{e}duit \`{a} $\Z$ mais pour lequel n\'{e}anmoins
$Z^4(X) \neq 0$}. Au paragraphe 6, on  \'{e}tablit le th\'{e}or\`{e}me \ref{illustr1}.
 La nullit\'{e} de $Z^4(X)$  a aussi \'{e}t\'{e} \'{e}tablie pour les
 hypersurfaces cubiques lisses de dimension 4 (\cite{voisinjapjmath}). Motiv\'{e}s par ce r\'{e}sultat,
 pour l'espace total d'une fibration en solides au-dessus
 d'une courbe, nous donnons des conditions suffisantes  sur la fibration
 qui impliquent
 l'annulation de $Z^4(X)$(Th\'{e}or\`{e}me \ref{theocritereabeljacobi}).

 Au paragraphe 7, on s'int\'{e}resse de fa\c con g\'{e}n\'{e}rale
 aux vari\'{e}t\'{e}s $X$ munies d'une fibration  $\pi : X \to \Gamma$ sur une courbe.
 L'indice d'une telle fibration $\pi$ est le pgcd des degr\'{e}s des multisections.
Lorsque les groupes de cohomologie coh\'{e}rente $H^{i}(X_{t},\O_{X_{t}})$
des fibres lisses sont nuls pour $i\geq 1$,
  on \'{e}tablit  par des m\'{e}thodes g\'{e}om\'{e}triques  ind\'{e}pendantes
  des paragraphes pr\'{e}c\'{e}dents
  des
 liens entre la nullit\'{e} de $Z_{2}(X)$ et le fait que l'indice  de la fibration soit \'{e}gal \`{a} 1 (Th\'{e}or\`{e}mes
 \ref{propindice1} et \ref{propindice2}).

La fibre g\'{e}n\'{e}rique $V/\C(\Gamma)$ d'une fibration comme ci-dessus
est une vari\'{e}t\'{e} projective et lisse sur un corps $F$  de caract\'{e}ristique
z\'{e}ro et de dimension cohomologique $cd(F)\leq 1$.
Au paragraphe  \ref{varsurcd1},  pour tout tel corps $F$ et toute vari\'{e}t\'{e} $V$
projective et lisse sur $F$,
en d\'{e}veloppant des techniques de $K$-th\'{e}orie alg\'{e}brique et de cohomologie galoisienne
(\cite{blochbook}, \cite{ctraskind}, \cite{kahn}),
on \'{e}tablit  des liens entre  le groupe de Chow de
codimension 2 de $V$ et la cohomologie non ramifi\'{e}e de degr\'{e}~3 de $V$
 (Th\'{e}or\`{e}me \ref{cdFleq1} et, sous l'hypoth\`{e}se $H^2(V,{\mathcal O}_{V})=0$, Th\'{e}or\`{e}me \ref{H2OXtrivial}).
Ceci  \'{e}tend certains r\'{e}sultats de  Bruno Kahn  \cite{kahn}.
Les applications aux vari\'{e}t\'{e}s fibr\'{e}es au-dessus de courbes sont donn\'{e}es
au paragraphe \ref{fibrationcourbe}. On \'{e}tablit ainsi par
cette m\'{e}thode ind\'{e}pendante
   l'annulation  (Th\'{e}or\`{e}me  \ref{famillesurfacesrationnelles}) de $H^3_{nr}(X,\Q/\Z(2))$
pour les solides $X$  fibr\'{e}s en surfaces rationnelles sur une courbe (ce
qui est un cas particulier du th\'{e}or\`{e}me \ref{illustr1}).  En  combinant $K$-th\'{e}orie alg\'{e}brique
et  m\'{e}thodes
g\'{e}om\'{e}triques, on obtient, pour les solides  $X$
fibr\'{e}s au-dessus d'une courbe $\Gamma$,
dont la fibre g\'{e}n\'{e}rique $V/\C(\Gamma)$ satisfait
$H^{i}(V,\O_{V})=0$ pour $i=1,\,2$,
 d'autres liens entre le groupe
 $H^3_{nr}(X,\Q/\Z(2))$, le groupe $Z_{2}(X)$
et l'indice de la fibration (Th\'{e}or\`{e}mes  \ref{Hitrivialbis} et   \ref{corconclus8}).

\medskip

{\bf Notations}

\'Etant donn\'{e} un groupe ab\'{e}lien $A$, un entier $n>0 $ et un nombre premier $l$,
on note $A[n] \subset A$ le sous-groupe form\'{e} des \'{e}l\'{e}ments annul\'{e}s par $n$, on note
$A\{l\} \subset A$ le sous-groupe de torsion $l$-primaire, et on note $A_{\rm tors} \subset A$
le sous-groupe de torsion.

\section{Rappels}
\subsection{La th\'{e}orie de Bloch--Ogus pour la cohomologie de Betti
\label{secrappelsbo}}

\medskip

On renvoie \`{a} \cite[\S 4]{artin}
 pour plus de d\'{e}tails sur ce qui suit.   Soit $X$ une vari\'{e}t\'{e}  alg\'{e}brique sur $\C$.
On dispose de l'espace topologique $X(\C)$.
On note $X_{cl}$ le site des isomorphismes locaux $U \to X(\C)$.
On a une fl\`{e}che de sites  $\delta : X_{cl} \to X(\C)$. Les topos associ\'{e}s sont \'{e}quivalents.

On a   le diagramme commutatif de morphismes de sites suivant
$$\begin{array}{ccccccccccccccccccccccccccccccccccccccccccccccccccc}
X_{cl}  &  \buildrel{\delta}\over \longrightarrow  & X(\C) \cr
 \downarrow{}{f}  &   &  \downarrow{h}{}  \cr
X_{\et} &  \buildrel{g}\over \longrightarrow & X_{Zar}
\end{array} $$

On note
$$\pi : X_{cl} \to X_{Zar} $$
le morphisme de sites $h \circ \delta = g \circ f$.

Soit $A$ un groupe ab\'{e}lien. On note $A(1)=A \otimes_{\Z} \Z(1)$
et $A(-1)= Hom_{\Z}(\Z(1),A)$. Ceci permet de d\'{e}finir $A(n)$
pour tout entier $n \in \Z$.

\begin{Defi} (Bloch--Ogus \cite{blochogus})
 Soit $X$ une vari\'{e}t\'{e} alg\'{e}brique sur $\C$.
 Les faisceaux $\mathcal{H}^i_{X}(A)$ sur $X$ sont les faisceaux pour la topologie de Zariski sur $X$
  d\'{e}finis par
$$\mathcal{H}^i_{X}(A):=R^i\pi_*A= R^ih_{*}(A).$$
\end{Defi}
En d'autres termes, $\mathcal{H}^i_{X}(A)$ est le faisceau pour la topologie de Zariski sur $X$
associ\'{e} au pr\'{e}faisceau pour la dite topologie d\'{e}fini par $U  \mapsto H^{i}(U(\C),A)$.
Pour une discussion de la d\'{e}finition suivante dans le contexte de la cohomologie  \'{e}tale
\`{a} coefficients finis, on consultera
 \cite{colliotojanguren} et \cite{ctsantabarbara}.
\begin{Defi}
 Soit $X$ une vari\'{e}t\'{e} alg\'{e}brique sur $\C$.
 On d\'{e}finit  le $i$-i\`{e}me groupe de    cohomologie non ramifi\'{e}e
de $X$  \`{a} coefficients dans $A$ par la formule
$$H^i_{nr}(X,A)=H^0(X,\mathcal{H}^i_{X}(A)).$$
\end{Defi}

Un r\'{e}sultat central de  l'article de Bloch  et Ogus \cite{blochogus} est la d\'{e}monstration de la conjecture
de Gersten pour diverses th\'{e}ories cohomologiques.
Ceci a comme cons\'{e}quence
 le th\'{e}or\`{e}me  suivant. Pour toute sous-vari\'{e}t\'{e} ferm\'{e}e int\`{e}gre $D$ de $X$, on note $i_D:D\rightarrow X$ l'inclusion
et  on note
$$H^i_{B}(\C(D),A) = \lim_{\overset{\rightarrow}
{\overset{U\subset D}
 {{\rm ouvert\,\, de \,\,Zariski  \,\,non  \,\,vide }}}} H^i(U(\C),A)$$ en tout point de $D$.
  Lorsque $D'\subset D$ est de codimension $1$, on a une fl\`{e}che induite par le r\'{e}sidu topologique (sur la normalisation de $D$)   (cf.  \cite[p. 417]{voisinlivre}) :
   $$Res_{D,D'}:H^i_{B}(\C(D),A)\rightarrow H^{i-1}_{B}(\C(D'),A(-1)).$$

Supposons $X$ irr\'{e}ductible. Pour $r\geq 0$ on note $X^{(r)}$ l'ensemble des sous-sch\'{e}mas
ferm\'{e}s de codimension $r$ dans $X$. Pour un groupe ab\'{e}lien $M$ et une sous-vari\'{e}t\'{e} ferm\'{e}e int\`{e}gre
$i_{D} : D \subset X$, on note $i_{D¬*}M$ l'image par $i_{D*}$ du faisceau constant $M_{D}$ sur $D$ d\'{e}fini par le groupe $M$.

\begin{theo}\label{boresolution} (\cite[Thm 4.2]{blochogus})
Pour toute vari\'{e}t\'{e} lisse connexe sur $\C$ et  tout entier $i\geq 1$,
on a une suite exacte
 de faisceaux   pour la topologie de Zariski sur $X$
 $$0\rightarrow \mathcal{H}^i_{X}(A)
\rightarrow i_{X*}H^i_{B}(\C(X),A)
\stackrel{\partial}{\rightarrow}
\bigoplus_{   D \in X^{(1)}  }
i_{D*}H^{i-1}_{B}(\C(D),A(-1))
\stackrel{\partial}{\rightarrow}\ldots
\stackrel{\partial}{\rightarrow}
\bigoplus_{D   \in X^{(i)}  }
 i_{D*}A_D(-i)\rightarrow0.$$
\end{theo}
Ici les fl\`{e}ches $\partial$ sont induites par les fl\`{e}ches $Res_{D,D'}$ lorsque $D'\subset D$ (et sont
nulles sinon). Le faisceau $A_D(-i)$ sur $D$ s'identifie bien s\^{u}r au faisceau constant
de fibre $H^{0}_{B}(\C(D),A(-i))$.

\medskip

Les cons\'{e}quences de ce th\'{e}or\`{e}me sont consid\'{e}rables.
Tout d'abord, notant $CH^k(X)/alg$ le groupe des cycles de codimension $k$ de $X$ modulo \'{e}quivalence alg\'{e}brique, on obtient la formule de Bloch--Ogus:
\begin{coro} \label{formuledebo}(\cite[Cor. 7.4]{blochogus})
Si $X$ est une vari\'{e}t\'{e} projective,
lisse et connexe sur $\C$, on a un isomorphisme canonique
 $$CH^k(X)/alg= H^k(X,\mathcal{H}^k_{X}(\Z(k))).$$
 \end{coro}

D'apr\`{e}s le
 th\'{e}or\`{e}me \ref{boresolution}, le faisceau $\mathcal{H}_{X}^i$ a une r\'{e}solution acyclique de longueur $\leq i$. On a donc le
 r\'{e}sultat
d'annulation suivant.
\begin{coro}\label{boannulation} Pour $X$   lisse,   $A$  un groupe ab\'{e}lien  et $r>i$, on
 a $H^r(X,\mathcal{H}_{X}^i(A))=0$.
\end{coro}

Bloch et Ogus d\'{e}duisent de cette annulation, par analyse de la suite spectrale de Leray de l'application continue
$\pi$, le r\'{e}sultat suivant, crucial pour la suite de cet article.
\begin{theo} \label{blochogussuitecourte}
 (\cite[Ex. 7.5]{blochogus})
Soit $X$ une vari\'{e}t\'{e} lisse sur $\C$ et soit $A$ un groupe ab\'{e}lien.
Alors la suite spectrale de Leray du morphisme de sites
$\pi :   X(\C)  \to X $ fournit une suite exacte
$$H^3(X(\C),A)\rightarrow H^3_{nr}(X,A)\stackrel{d_2}{\rightarrow} H^2(X,\mathcal{H}_{X}^2(A))  \to H^4(X(\C),A)$$
 Pour $A=\Z(2)$ et $X$ de plus projective,
   la suite exacte ci-dessus donne la suite
exacte
\begin{eqnarray}\label{exacte17jan}
 H^3(X(\C),\Z(2)) \to H^3_{nr}(X,\Z(2)) \to CH^2(X)/alg  \stackrel{c}{\rightarrow}    H^4(X(\C),\Z(2))
 \end{eqnarray}
o\`u l'application $c : CH^2(X)/alg \to H^4(X(\C),\Z(2))$ est    induite par l'application   classe de cycles  pour la cohomologie de Betti.
\end{theo}

\medskip

Le noyau de $  c :  CH^2(X)/alg \to   H^4(X(\C),\Z(2)) $ est par d\'{e}finition  le groupe de Griffiths
${\rm Griff}^2(X)$ des cycles de codimension $2$ de $X$ homologues \`{a} $0$ modulo \'{e}quivalence alg\'{e}brique.
La th\'{e}orie de Bloch--Ogus donne donc le

\begin{theo}(Bloch--Ogus, \cite[Cor. 7.4]{blochogus}) \label{theogriffblochogus}
Pour $X$ projective et lisse, il y a un isomorphisme
$$ H^3_{nr}(X, \Z(2))/ \Im[ H^3(X(\C),\Z(2))] \oi  {\rm Griff}^2(X).$$
\end{theo}

Voici une autre application du th\'{e}or\`{e}me \ref{boresolution}.

\begin{theo}
Soit $A$ un groupe ab\'{e}lien. Pour tout entier $i\geq 0$,

(i)  le groupe $H^{i}_{nr}(X,A)$,

(ii) l'image de $H^{i}(X(\C),A)$ dans $H^{i}_{nr}(X,A) \subset H^{i}_{B}(\C(X),A)$,

(iii) le quotient du premier groupe par
le second,

\noindent sont des
  invariants birationnels des vari\'{e}t\'{e}s projectives et lisses sur $\C$, invariants qui s'annulent
sur les vari\'{e}t\'{e}s rationnelles.
\end{theo}

{\bf D\'{e}monstration.}
L'invariance birationnelle de $H^{i}_{nr}(X,A)$ est une cons\'{e}\-quence
du th\'{e}or\`{e}me \ref{boresolution} (validit\'{e} de la conjecture de Gersten), voir
 \cite[Thm. 4.4.1]{ctsantabarbara}.

L'invariance birationnelle de l'image de $H^{i}(X(\C),A)$ dans $H^{i}_{B}(\C(X),A)$,
est \'{e}tablie par Grothendieck dans \cite[III (9.2)]{grothendieck}. Comme on est
en caract\'{e}ristique z\'{e}ro, on peut utiliser le th\'{e}or\`{e}me d'Hironaka et
r\'{e}duire la d\'{e}monstration
au cas de l'\'{e}clatement d'une sous-vari\'{e}t\'{e} ferm\'{e}e  lisse. On sait comment
se calcule la cohomologie d'un tel \'{e}clatement : tout ce qu'on ajoute
dans la cohomologie est support\'{e} dans le diviseur exceptionnel et
donc s'annule par passage au corps des fonctions. \cqfd

\subsection{$K$-th\'{e}orie alg\'{e}brique, cohomologie \'{e}tale et conjecture de Bloch--Kato }\label{secrappelblochkato}

Sur tout sch\'{e}ma $X$ et tout entier $i\geq 0$, on note ${\mathcal K}_{i,X}$  le faisceau
  pour la topologie
de Zariski  sur $X$ associ\'{e} au pr\'{e}faisceau $U \mapsto K_{i}(H^0(U,{\mathcal O}_{X}))$.
Le cas $n=0$ du th\'{e}or\`{e}me suivant (conjecture de Gersten  pour la $K$-th\'{e}orie alg\'{e}brique)
est un c\'{e}l\`{e}bre r\'{e}sultat de Quillen. On trouvera une d\'{e}monstration du cas g\'{e}n\'{e}ral
dans \cite{ctkh}.

\begin{theo}\label{quillenresolution}
Soit $X$ une   vari\'{e}t\'{e} lisse connexe sur un corps $F$.
Pour tout entier $i\geq 1$, et tout entier $n \geq 0$,
on a une suite exacte
 de faisceaux   pour la topologie de Zariski sur $X$
 $$0\rightarrow \mathcal{K}_{i,X}/n
\rightarrow i_{X*}K_{i}(F(X))/n
\stackrel{\partial}{\rightarrow}
\bigoplus_{   D \in X^{(1)}  }
i_{D*}K_{i-1}(F(D))/n
\stackrel{\partial}{\rightarrow}\ldots
\stackrel{\partial}{\rightarrow}
\bigoplus_{D   \in X^{(i)}  }
 i_{D*}\Z /n \rightarrow 0.$$
\end{theo}

\`A tout anneau commutatif unitaire $A$ et tout entier $i\geq 0$ et tout entier $i \geq 0$
on associe le groupe   $K_{i}^M(A)$ quotient du produit tensoriel
$A^{\times}\otimes_{\Z} \ldots \otimes_{\Z}A^{\times}$ ($i$ fois)  par le sous-groupe engendr\'{e} par les \'{e}l\'{e}ments
$x_{1} \otimes \ldots \otimes x_{i}$ avec $x_{r}+x_{s}=1$ pour un couple d'entiers $r<s$.
Pour $A$ un corps, c'est la d\'{e}finition de Milnor.
On a  $K_{1}^M(A)=A^*$.
Sur tout sch\'{e}ma $X$ et tout entier $i\geq 0$, on note ${\mathcal K}^M_{i,X}$  le faisceau
  pour la topologie
de Zariski  sur $X$ associ\'{e} au pr\'{e}faisceau $U \mapsto K^M_{i}(H^0(U,{\mathcal O}_{X}))$.
L'analogue  du th\'{e}or\`{e}me de Quillen pour la $K$-th\'{e}orie de Milnor a fait l'objet
de travaux de nombreux auteurs. On a maintenant le r\'{e}sultat suivant
(M. Kerz \cite[Thm. 7.1]{kerz}) :

\begin{theo}\label{kerzresolution}
Soit $X$ une   vari\'{e}t\'{e} lisse connexe sur un corps infini $F$.
Pour tout entier $i\geq 1$, et tout entier $n \geq 0$,
on a une suite exacte
 de faisceaux   pour la topologie de Zariski sur $X$
 $$0\rightarrow \mathcal{K}^M_{i,X}/n
\rightarrow i_{X*}K^M_{i}(F(X))/n
\stackrel{\partial}{\rightarrow}
\bigoplus_{   D \in X^{(1)}  }
i_{D*}K^M_{i-1}(F(D))/n
\stackrel{\partial}{\rightarrow}\ldots
\stackrel{\partial}{\rightarrow}
\bigoplus_{D   \in X^{(i)}  }
 i_{D*}\Z /n \rightarrow 0.$$
\end{theo}

Soit $M$ un module continu discret sur le groupe de Galois absolu $G$  de $F$,
et $i\geq 0$ un entier.
On note $H^{i}(F,M)=H^{i}(G,M)$ le $i$-i\`{e}me groupe de  cohomologie galoisienne de $M$.

\begin{theo}\label{boresolutionfinie}(Bloch--Ogus \cite{blochogus})
Soit $X$ une   vari\'{e}t\'{e} lisse connexe sur un corps $F$.
Soit $n$ un entier positif premier \`{a} la caract\'{e}ristique de $F$,
et soit $\mu_{n}$ le faisceau \'{e}tale sur $X$ d\'{e}fini par les racines $n$-i\`{e}mes de l'unit\'{e}.
Soit $ \mathcal{H}^{i}_{X}(\mu_{n}^{\otimes i})$ le faisceau pour la topologie
de Zariski sur $X$ associ\'{e} au pr\'{e}faisceau $U \mapsto H^{i}_{\et}(U,\mu_{n}^{\otimes i})$.
Pour tout entier $i\geq 1$,  tout entier $j \in \Z$
et tout entier $n \geq 0$,
on a une suite exacte
 de faisceaux   pour la topologie de Zariski sur $X$
 $$0\rightarrow \mathcal{H}^{i}_{X}(\mu_{n}^{\otimes j})
\rightarrow i_{X*}H^{i}(F(X),\mu_{n}^{\otimes j})
\stackrel{\partial}{\rightarrow}
\bigoplus_{   D \in X^{(1)}  }
i_{D*} H^{i-1}(F(D),\mu_{n}^{\otimes j-1})
\stackrel{\partial}{\rightarrow}\ldots$$$$
\stackrel{\partial}{\rightarrow}
\bigoplus_{D   \in X^{(i)}  }
 i_{D*}H^0(F(D), \mu_{n}^{\otimes j-i})    \rightarrow 0.$$
\end{theo}
En d'autres termes, quand on prend les sections de cette suite sur un anneau local de $X$,
on obtient  une suite exacte de groupes ab\'{e}liens.
Ceci vaut encore quand on prend
les sections sur un anneau semilocal de $X$ (voir \cite{ctkh}).

\medskip

Pour toute vari\'{e}t\'{e} $X$ sur $\C$, tout groupe ab\'{e}lien {\it de torsion} $A$, et tout entier $i \geq 0$,
 les th\'{e}or\`{e}mes de comparaison de SGA4 \cite{artin} donnent des isomorphismes canoniques \break
$H^{i}_{\et}(X,A) \oi H^{i}(X(\C),A)$ et  $H^{i}(\C(X),A) \oi H^{i}_{B}(\C(X),A)$. On en d\'{e}duit un
isomorphisme canonique    entre  le faisceau Zariski
$\H^{i}_{X}(A)$ obtenu par faisceautisation de $U \mapsto H^{i}_{\et}(U,A)$ et
 le faisceau Zariski
$\H^{i}_{X}(A)$ obtenu par faisceautisation de $U \mapsto H^{i}(U(\C),A)$ ,
ce qui nous autorise \`{a} utiliser une notation unique pour ces deux faisceaux.
Pour $X$ lisse sur $\C$, et $A=\mu_{n}^{\otimes j}$,
cela permet
d'identifier les suites exactes de faisceaux des th\'{e}or\`{e}mes \ref{boresolution} et \ref{boresolutionfinie}.

\medskip

 Notons $CH^{i}(X)$ le groupe de Chow des cycles de codimension $i$ modulo
 l'\'{e}quivalence rationnelle et $H^{i}_{nr}(X,\mu_{n}^{\otimes j})=H^0(X,\mathcal{H}^{i}_{X}(\mu_{n}^{\otimes j}))$.
  Du th\'{e}or\`{e}me \ref{boresolutionfinie}, on tire comme au paragraphe pr\'{e}c\'{e}dent
 des isomorphismes $$CH^{i}(X)/n \oi H^{i}(X,\mathcal{H}^{i}_{X}(\mu_{n}^{\otimes i}))$$
et des suites exactes
$$ H^3_{\et}(X,\mu_{n}^{\otimes 2}) \to H^3_{nr}(X,\mu_{n}^{\otimes 2}) \to CH^2(X)/n \to H^4_{\et}(X,\mu_{n}^{\otimes 2}).$$

\medskip

Rappelons l'\'{e}nonc\'{e} de la conjecture  de Bloch--Kato en K-th\'{e}orie alg\'{e}brique.

 ($BK_{i,n}$)   Soient $i\geq 1$ et $n \geq 2$ des entiers. Pour tout corps $F$
 de caract\'{e}ristique premi\`{e}re \`{a} $n$, l'application de r\'{e}ciprocit\'{e} (d\'{e}finie par Tate)
 $$K_{i}^M(F)/n \to H^{i}(F,\mu_{n}^{\otimes i}),$$
 qui envoie les   groupes  de $K$-th\'{e}orie de Milnor
 modulo $n$
 vers les groupes de cohomologie galoisienne,
 est un isomorphisme.

Des cas particuliers importants de cette conjecture ont \'{e}t\'{e} d\'{e}montr\'{e}s par
Merkur'ev et Suslin (\cite{MS1}, \cite{MS2}),  Rost \cite{R},  et Voevodsky
\cite{voevodskymod2}. Le cas g\'{e}n\'{e}ral vient d'\^{e}tre \'{e}tabli par Voevodsky
 \cite{voevodskymodl} et Rost   (cf. \cite{suslinjoukhovitski}), voir aussi le texte de Weibel  \cite{weibel}.

Il est facile, dans chacun des  \'{e}nonc\'{e}s du pr\'{e}sent article,
de d\'{e}terminer lesquels parmi les \'{e}nonc\'{e}s $BK_{i,n}$ sont requis
pour les d\'{e}monstrations.

\medskip

Pour  $X$ une vari\'{e}t\'{e} lisse sur un corps infini $F$ et $n>0$ un entier premier \`{a} la caract\'{e}ristique de $F$,
les applications de r\'{e}ciprocit\'{e} d\'{e}finissent un morphisme de la suite exacte du th\'{e}or\`{e}me
\ref{kerzresolution} dans la suite exacte du th\'{e}or\`{e}me \ref{boresolutionfinie} (cf. \cite{kerz}).
Sous la conjecture de Bloch--Kato,
on d\'{e}duit  de ces deux th\'{e}or\`{e}mes
 un isomorphisme de faisceaux pour la topologie de Zariski
$ {\mathcal K}^M_{i,X}/n \oi {\mathcal H}^{i}_{X}(\mu_{n}^{\otimes i}).$
Comme expliqu\'{e} dans \cite[Thm. 7.8]{kerz}, une m\'{e}thode due \`{a} Hoobler permet
d'\'{e}tablir un \'{e}nonc\'{e} plus g\'{e}n\'{e}ral : {\it Sous la conjecture de Bloch--Kato,
pour tout anneau semilocal $A$ contenant un corps $F$  infini et tout entier $n$
premier \`{a} la caract\'{e}ristique de $F$,  pour tout entier $i \geq 1$,
l'application de r\'{e}ciprocit\'{e} $K_{i}^M(A)/n \to H^{i}_{\et}(A,\mu_{n}^{\otimes i})$
est un isomorphisme.}

\medskip

Pour $F$ un corps, on a des homomorphismes $K_{i}^M(F) \to K_{i}(F)$, qui sont des isomorphismes
pour $i=0,1,2$. Pour $i\leq 2$, on a donc $ {\mathcal K}_{i,X}/n \oi {\mathcal H}^{i}_{X}(\mu_{n}^{\otimes i}).$

\section{Cohomologie non ramifi\'{e}e et  conjecture de Hodge enti\`{e}re}

\subsection{G\'{e}n\'{e}ralisation d'un argument de Bloch et Srinivas}

Dans cette partie on d\'{e}taille et g\'{e}n\'{e}ralise
un argument de Bloch et Srinivas \cite{blochsrinivas}.
Soit $X$ une vari\'{e}t\'{e} alg\'{e}brique sur $\C$.  Les notations sont celles du paragraphe \ref{secrappelsbo}.

\begin{theo}\label{principal0}
 Soit $X$ une vari\'{e}t\'{e} alg\'{e}brique connexe sur $\C$.
Pour tout entier $i$,
la multiplication par un entier $n>0$  sur les faisceaux $ \H^{p}_{X}(\Z(i))$
induit des suites exactes  courtes de faisceaux Zariski sur $X$
$$ 0 \to \H^{p}_{X}(\Z(i)) \buildrel{\times n}\over{\to} \H^{p}_{X}(\Z(i))  \to \H^{p}_{X}(\mu_{n}^{\otimes i})  \to 0.$$
En particulier les faisceaux $\H^{p}_{X}(\Z(i))$ sont  sans torsion, et donc
 leurs groupes de sections globales $H^0(X,\H^{p}_{X}(\Z(i)))=H^p_{nr }(X,\Z(i))$ sont sans torsion.
\end{theo}

{\bf D\'{e}monstration.}
Les suites exactes de groupes ab\'{e}liens
$$ 0 \to \Z(i) \buildrel{\times n}\over \longrightarrow \Z(i) \to \mu_{n}^{\otimes i} \to 0$$
 donnent
naissance \`{a} de longues suites exactes de groupes de cohomologie
sur tout $U(\C)$ pour $U \subset X$ ouvert Zariski, et donc
\`{a} de longues suites de
 faisceaux sur $X_{Zar}$
$$\cdots \to   \H^{p-1}_{X}(\mu_{n}^{\otimes i})   \to \H^{p}_{X}(\Z(i))  \buildrel{\times n}\over \longrightarrow \H^{p}_{X}(\Z(i)) \to  \H^{p}_{X}(\mu_{n}^{\otimes i})
\to \H^{p+1}_{X}(\Z(i)) \to \cdots.$$

Il s'agit donc de montrer que les fl\`{e}ches $\H^{p}_{X}(\Z(i)) \to  \H^{p}_{X}(\mu_{n}^{\otimes i})$
sont surjectives. Il suffit pour cela d'\'{e}tablir que les fl\`{e}ches
$\H^{p}_{X}(\Z(p)) \to  \H^{p}_{X}(\mu_{n}^{\otimes p})$
le sont.

On a le diagramme commutatif de faisceaux sur $X_{cl}$

$$\begin{array}{ccccccccccccccccccccccccccccccccccccccccccccccccccc}
0 &  \longrightarrow & \Z(1)                                            & \longrightarrow &  \O_{X,cl}                                         & \buildrel{exp}\over \longrightarrow            & \O_{X,cl}^{\times}              &    \longrightarrow  & 1  \cr
  &        &       \downarrow{} &      &   \downarrow{}      &                                                                     &   \downarrow{id}      &    &    \cr
    0 &  \longrightarrow &     \mu_{n}                                    & \longrightarrow &  \O_{X,cl}^{\times}                                    &  \buildrel{z \mapsto z^n}\over \longrightarrow  &     \O_{X,cl}^{\times}          &      \longrightarrow  &  1\cr
\end{array} $$
o\`u les deux fl\`{e}ches  verticales de gauche sont d\'{e}finies par $z \mapsto exp(z/n)$.

Dans ce diagramme, le faisceau $\O_{X,cl}$ d\'{e}signe le faisceau
des fonctions continues sur $X(\C)$ \`{a} valeurs complexes
sur $X(\C)$, vu sur $X_{cl}$.

De ce diagramme on d\'{e}duit un diagramme commutatif de faisceaux sur $X_{Zar}$ :
$$\begin{array}{ccccccccccccccccccccccccccccccccccccccccccccccccccc}
\pi_{*}\O_{X,cl}^{\times}   & \to & R^1\pi_{*}(\Z(1)  ) \cr
  \downarrow{}{id}  & &   \downarrow{}{}  \cr
    \pi_{*}\O_{X,cl}^{\times}   & \to  & R^1\pi_{*}(\mu_{n})
    \end{array} $$

   On cherche \`{a} identifier la fl\`{e}che compos\'{e}e $\O_{X,Zar}^{\times} \to   \pi_{*}\O_{X,cl}^{\times}    \to  R^1\pi_{*}(\mu_{n})$.

 D'apr\`{e}s \cite[Thm. 4.4]{artin},
on a  $R^1f_{*}\mu_{n}=0$.
 La suite
 exacte de faisceaux sur $X_{cl}$
   $$ 1 \to \mu_{n} \to  \O_{X,cl}^{\times} \buildrel{z \mapsto z^n}\over \longrightarrow \O_{X,cl}^{\times} \to 1$$
 induit donc
 une  suite exacte de faisceaux sur $X_{\et}$
 $$ 1 \to f_{*}\mu_{n} \to f_{*}\O_{X,cl}^{\times}   \buildrel{z \mapsto z^n}\over \longrightarrow  f_{*}\O_{X,cl}^{\times} \to 1.$$
 On a le diagramme commutatif de faisceaux sur $X_{\et}$
 $$\begin{array}{ccccccccccccccccccccccccccccccccccccccccccccccccccc}
 1 & \to  &  f_{*}\mu_{n}  &  \to  &  f_{*}\O_{X,cl}^{\times}    &  \buildrel{z \mapsto z^n}\over \longrightarrow   &  f_{*}\O_{X,cl}^{\times}  &  \to   & 1 \cr
 && \uparrow{}{\simeq } &&    \uparrow&  &  \uparrow && \cr
 1  &  \to  &  \mu_{n}  &  \to  &  \O_{X,\et}^{\times}    & \buildrel{z \mapsto z^n}\over \longrightarrow   &  \O_{X,\et}^{\times}  &  \to  &  1
  \end{array} $$
  d'o\`u l'on d\'{e}duit un diagramme commutatif sur $X_{Zar}$
  $$\begin{array}{ccccccccccccccccccccccccccccccccccccccccccccccccccc}
g_{*} f_{*}\O_{X,cl}^{\times} & \to & R^1g_{*}(f_{*}\mu_{n} )\cr
\uparrow{}{ } &  &  \uparrow{}{\simeq} \cr
g_{*}\O_{X,\et}^{\times} & \to &  R^1g_{*}\mu_{n}
  \end{array} $$

  Comme $R^1f_{*}\mu_{n}=0$, la
  fl\`{e}che naturelle $R^1g_{*}(f_{*}\mu_{n}) \to R^1\pi_{*}\mu_{n}$
  est   un isomorphisme.

  Le diagramme ci-dessus se r\'{e}\'{e}crit
    $$\begin{array}{ccccccccccccccccccccccccccccccccccccccccccccccccccc}
\pi_{*}\O_{X,cl}^{\times} & \to & R^1\pi_{*}\mu_{n} \cr
\uparrow{}{ } &  &  \uparrow{}{\simeq} \cr
g_{*}\O_{X,\et}^{\times} & \to &  R^1g_{*}\mu_{n}
  \end{array} $$
  On obtient    le diagramme commutatif
  $$\begin{array}{ccccccccccccccccccccccccccccccccccccccccccccccccccc}
  \pi_{*}\O_{X,cl}^{\times}   & \to & R^1\pi_{*}(\Z(1)  ) \cr
  \downarrow{}{id}  & &   \downarrow{}{}  \cr
    \pi_{*}\O_{X,cl}^{\times}   & \to  & R^1\pi_{*}\mu_{n} \cr
    \uparrow{}{} &  & \uparrow{}{\simeq} \cr
    g_{*}\O_{X,\et}^{\times} & \to &  R^1g_{*}\mu_{n}
\end{array} $$
Le th\'{e}or\`{e}me de Hilbert 90, sous la forme de Grothendieck,
assure $R^1g_{*}\O_{X,\et}^{\times}=0$.
Ainsi la fl\`{e}che  $g_{*}\O_{X,\et}^{\times}   \to    R^1g_{*}\mu_{n}$
est surjective. Du diagramme ci-dessus il r\'{e}sulte que la fl\`{e}che
  $$R^1\pi_{*}(\Z(1)  ) \to R^1\pi_{*} \mu_{n}$$
  est surjective, ce qui se traduit encore ainsi :
  le faisceau $R^2\pi_{*}(\Z(1)  )$ est sans torsion
  (\cite[Lemma 3.2]{barbieriviale}).

  Observons que la fl\`{e}che $\O_{X,Zar}^{\times}  \to g_{*}\O_{X,\et}^{\times}$
  est un isomorphisme. On  a donc \'{e}tabli le fait suivant : la fl\`{e}che compos\'{e}e

$$\O_{X,Zar}^{\times}  \to  \pi_{*}\O_{X,cl}^{\times}    \to  R^1\pi_{*}(  \Z(1)  ) \to
R^1\pi_{*}\mu_{n} \to R^1g_{*}\mu_{n},$$
o\`u la derni\`{e}re fl\`{e}che est l'isomorphisme r\'{e}ciproque du compos\'{e} d'isomorphismes
$$R^1g_{*}\mu_{n} \oi R^1g_{*}(f_{*}\mu_{n}) \oi R^1(g\circ f)_{*}(\mu_{n})=R^1\pi_{*}\mu_{n},$$
est la surjection  naturelle de faisceaux induite par la suite de Kummer.

On prend maintenant des produits tensoriels et des cup-produits.
On obtient ainsi une suite d'homomorphismes
$$(\O_{X,Zar}^{\times})^{\otimes p}      \to  R^p\pi_{*}((\Z(p)  ) \to
R^p\pi_{*}(\mu_{n}^{\otimes p}) \to R^pg_{*}(\mu_{n}^{\otimes p}),$$
  le compos\'{e} $(\O_{X,Zar}^{\times})^{\otimes p} \to
R^pg_{*}(\mu_{n}^{\otimes p})$ est simplement le compos\'{e}
$$(\O_{X,Zar}^{\times})^{\otimes p} \to ( R^1g_{*}\mu_{n})^{\otimes p} \to R^pg_{*}(\mu_{n}^{\otimes p}),$$
la premi\`{e}re fl\`{e}che \'{e}tant donn\'{e}e par  la suite de Kummer.

Comme   $R^{i}f_{*}(\mu_{n}^{\otimes p})=0$ pour tout $i>0$ (\cite{artin}, loc. cit.),
la fl\`{e}che
$ R^p\pi_{*}(\mu_{n}^{\otimes p}) \to R^pg_{*}(\mu_{n}^{\otimes p})$
est un isomorphisme.

Comme rappel\'{e} \`{a} la fin de la section \ref{secrappelblochkato},
la conjecture de Bloch--Kato
implique
que la fl\`{e}che de faisceaux Zariski $$(\O_{X,Zar}^{\times})^{\otimes p} \to
R^pg_{*}(\mu_{n}^{\otimes p})$$ est surjective.
Notons que l'on ne
fait pas ici d'hypoth\`{e}se de lissit\'{e} sur $X$.

La combinaison de ces r\'{e}sultats implique
que la fl\`{e}che
$$R^p\pi_{*}(\Z(p)  ) \to
R^p\pi_{*}(\mu_{n}^{\otimes p})$$ est surjective.
Ceci \'{e}tablit le th\'{e}or\`{e}me.
\cqfd

{\begin{rema}  {\rm  Des cas particuliers du  th\'{e}or\`{e}me \ref{principal0}
avaient \'{e}t\'{e} \'{e}tablis.

Pour $X$ lisse et $p=2$, Bloch et Srinivas  \cite{blochsrinivas} (voir la d\'{e}monstration du th\'{e}or\`{e}me 1 (ii), p.~1240)
 utilisent
le th\'{e}or\`{e}me de Merkur'ev-Suslin
pour montrer que  la suite ci-dessus est exacte, et que donc le faisceau $\H^3_{X}(\Z)$
est sans torsion.

Barbieri-Viale  \cite{barbieriviale}  observe que lorsque la dimension de $X$ est 3, pour tout ouvert affine $U$ de $X$,
on a $H^{r}(U(\C),\Z)=0$ et $H^{r}(U(\C),\Z/n)=0$ pour $r\geq 4$ (th\'{e}or\`{e}me de Lefschetz faible, cf. \cite[Thm. 13.2.22]{voisinlivre}).
Ainsi les faisceaux  $\H^r_{X}(\Z)$ et $\H^r_{X}(\Z/n)$ sont nuls pour $r \geq 4$.
Le th\'{e}or\`{e}me de Merkur'ev--Suslin suffit alors  \`{a}
\'{e}tablir le th\'{e}or\`{e}me ci-dessus pour les vari\'{e}t\'{e}s de dimension 3.

  L. Barbieri-Viale nous signale une d\'{e}monstration du th\'{e}or\`{e}me  \ref{principal0} dans sa  pr\'{e}publication   \cite{barbieriviale2}.
  }

\end{rema}

\subsection{Action des correspondances et  cohomologie des faisceaux $\H^n(\Z)$}

 \begin{prop}\label{lemmecorresp} Soit $X$ une vari\'{e}t\'{e} projective, lisse, connexe sur $\C$, de dimension $d$.

(i) Si   le groupe $CH_0(X) $ est support\'{e} sur un ferm\'{e} alg\'{e}brique
 $Y\subset X$ avec ${\rm dim} \,Y = r$,
  alors $H^0(X,\mathcal{H}^p_X(A))$ est annul\'{e} par un entier $N\not=0$ pour $p>r$. En particulier
$H^0(X,\mathcal{H}^p_X(\mathbb{Z}))=0$  pour $p>r$.

 (ii) Si   le groupe $CH_0(X) $ est support\'{e} sur un ferm\'{e} alg\'{e}brique
 $Y\subset X$ avec ${\rm codim} \,Y = r$,
alors pour $p<r$ le groupe
$H^p(X,\mathcal{H}^d_X(A))$ est annul\'{e} par un entier $N>0$.

(iii) Pour l'espace projectif $X=\mathbb{P}^d_{\C}$, on a $H^p(X,\mathcal{H}^q_X(A))=0$, pour $p\not=q$
 et $H^p(X,\mathcal{H}^p_X(A))=A$ pour tout $p\leq d$.

 \end{prop}
 {\bf D\'{e}monstration.}
  Les correspondances modulo \'{e}quivalence alg\'{e}brique agissent
  sur la cohomologie des faisceaux $\H^{i}(A)$, pour tout
  groupe ab\'{e}lien $A$. Pour la commodit\'{e} du lecteur, nous avons inclus en appendice une preuve de ce fait dans
  le contexte de la cohomologie de Betti.
  Pour la cohomologie  \'{e}tale \`{a} coefficients des groupes $\Z/n$,
  dans le cadre des modules de cycles de Rost \cite{rostchowgroupswith},
 une telle action est construite
  par  Merkur'ev dans \cite{merkurevcorrespondances}.

 \medskip

 (i)
La d\'{e}composition de la diagonale (cf. \cite{blochsrinivas}) dit que sous l'hypoth\`{e}se
de la proposition, on a  pour un entier $N\not=0$ une \'{e}galit\'{e}
de cycles modulo \'{e}quivalence rationnelle sur $X\times X$:
$$N\Delta_X=\Gamma_1+\Gamma_2\,\,\,{\rm dans} \,\,\, CH^d(X\times X),$$
o\`{u} $\Gamma_1$ est support\'{e} sur $Y\times X$, o\`{u} l'on peut supposer que $Y$ est de pure  dimension $r$,
et $\Gamma_2$ est support\'{e} sur
$X\times D$,  avec $D\subsetneqq X$.
Prenant les actions de ces correspondances sur $H^p(X,\mathcal{H}^d_X(A))$, on obtient
l'\'{e}galit\'{e}:
$$N Id=\Gamma_{1*}+\Gamma_{2*}:H^0(X,\mathcal{H}_X^p(A))\rightarrow H^0(X,\mathcal{H}_X^p(A)).$$
Introduisant des d\'{e}singularisations $\widetilde{D}$
de $D$ et $\widetilde{Y}$ de $Y$ et des rel\`{e}vements de $\Gamma_1$, $\Gamma_2$,
on constate que $\Gamma_{1*}$ se factorise par la fl\`{e}che de restriction
 $$H^0(X,\mathcal{H}_X^p(A)) \to  H^{0}(\widetilde{Y},\mathcal{H}_{\widetilde{Y}}^{p}(A)).$$
 Or le faisceau   $\mathcal{H}_{\widetilde{Y}}^{p}(A)$ est nul car $\dim({\widetilde{Y}})<p$.
 Par ailleurs $\Gamma_{2*}$ a   son image  dans $H^0(X,\mathcal{H}_X^p(A))$ constitu\'{e}e de classes  \`{a} support
 dans le ferm\'{e} $D \subsetneqq X$, donc nulles (Th\'{e}or\`{e}me \ref{boresolution}).
 Ainsi $N H^0(X,\mathcal{H}_X^p(A))=0$.   Si  $A=\mathbb{Z}$,  le th\'{e}or\`{e}me \ref{principal0} dit que $ H^0(X,\mathcal{H}_X^p(\mathbb{Z}))$
 est sans torsion et donc
 $H^0(X,\mathcal{H}_X^p(\mathbb{Z}))=0$.

 (ii)
 En permutant les facteurs, on trouve  une \'{e}galit\'{e}
de cycles modulo \'{e}quivalence rationnelle sur $X\times X$:
$$N\Delta_X=\Gamma_1+\Gamma_2\,\,\,{\rm dans} \,\,\, CH^d(X\times X),$$
o\`{u} $\Gamma_1$ est support\'{e} sur $D\times X,\,D\subsetneqq X$ et $\Gamma_2$ est support\'{e} sur
$X\times Y$, o\`{u} l'on peut supposer que $Y$ est de pure codimension $r$.
Prenant les actions de ces correspondances sur $H^p(X,\mathcal{H}^d_X(A))$, on obtient
l'\'{e}galit\'{e}:
$$N Id=\Gamma_{1*}+\Gamma_{2*}:H^p(X,\mathcal{H}_X^d(A))\rightarrow H^p(X,\mathcal{H}_X^d(A)).$$
Introduisant des d\'{e}singularisations $\widetilde{D}$
de $D$ et $\widetilde{Y}$ de $Y$ et des rel\`{e}vements de $\Gamma_1$, $\Gamma_2$,
on constate que $\Gamma_{1*}$ se factorise par la fl\`{e}che de restriction
$$H^p(X,\mathcal{H}_X^d(A))\rightarrow H^p(\widetilde{D},\mathcal{H}_{\widetilde{D}}^d(A))$$
qui est nulle car $\dim\,D<d$
 et donc $\mathcal{H}_{\widetilde{D}}^d(A))=0$.
 Enfin
$\Gamma_{2*}$ se factorise par la fl\`{e}che
$$\widetilde{j}_*:H^{p-r}(\widetilde{Y},\mathcal{H}_{\widetilde{Y}}^{d-r}(A))\rightarrow H^p(X,\mathcal{H}_X^d(A))$$
induite par $j : \widetilde{Y} \to X$,
qui est nulle car $p<r$. Donc $NId=0$ sur $H^p(X,\mathcal{H}_X^d(A))$.

(iii) Pour $X=\mathbb{P}^d_{\C}$, on a une d\'{e}composition compl\`{e}te de la diagonale dans $CH^d(X\times X)$
$$\Delta_X=\sum_{i=0}^{d} h_1^ih_2^{d-i},$$
o\`{u} $h=c_1(\mathcal{O}_{\mathbb{P}^d}(1))\in CH^1(\mathbb{P}^d)$ et $h_i:=pr_i^*h,\,i=1,\,2$.
On en d\'{e}duit que pour $\alpha\in H^p(X,\mathcal{H}^q_X(A))$, on a
$\Delta_{X*}\alpha=\sum_{i=0}^{d}  (h_1^ih_2^{d-i})_*\alpha$.
Soient $X_i=\mathbb{P}^i\stackrel{j_i}{\hookrightarrow }\mathbb{P}^d$ et $\pi_i:\mathbb{P}^i\rightarrow pt$ l'application constante. Alors $h_1^ih_2^{d-i}$ est la classe de
$X_{d-i}\times X_{i}$ dans $X\times X$ et on  a 
$$(h_1^ih_2^{d-i})_*\alpha= {j_{i}}_*(\pi_{i}^*({\pi_{d-i}}_*(j_{d-i}^*\alpha))),$$
ce qui montre  que $H^p(X,\mathcal{H}^q_X(A))=0$ pour $p\not=q$ car ${\pi_{d-i}}_*\circ j_{d-i}^*$ s'annule sur
$H^p(X,\mathcal{H}^q_X(A))$ pour $p\not=q$. L'assertion restante r\'{e}sulte du m\^{e}me argument ou, si $A=\mathbb{Z}$,
du corollaire  (\ref{formuledebo}).

 \cqfd

\bigskip

 \begin{prop}\label{invbir}
 Pour tout groupe ab\'{e}lien
 $A$ et tout entier $i\geq 0$, les groupes $H^i(X, \mathcal{H}_X^d(A))$ sont des invariants birationnels des
 vari\'{e}t\'{e}s connexes, projectives et lisses sur $\C$ de dimension~$d$.
 \end{prop}
 {\bf D\'{e}monstration.}  Soient $X$ et $Y$ deux vari\'{e}t\'{e}s projectives lisses sur $\mathbb{C}$ de dimension
 $d$
 et $\phi:X\dashrightarrow Y$ une application birationnelle. Soit $\Gamma_\phi\subset X\times Y$
 le graphe de $\phi$ et $\Gamma_{\phi^{-1}}\subset Y\times X$ le graphe de $\phi^{-1}$.
 Regardons ces correspondances modulo \'{e}quivalence alg\'{e}brique, de fa\c{c}on \`{a} pouvoir les composer.
 \begin{lemm}\label{lemme9avril} Le compos\'{e} $\Gamma_{\phi^{-1}}\circ\Gamma_\phi\in CH^d(X\times X)/alg$
 se d\'{e}compose de la fa\c{c}on suivante:
 \begin{eqnarray} \label{decomp9avril} \Gamma_{\phi^{-1}}\circ\Gamma_\phi=\Delta_X+Z\,\,{\rm dans}\,\,CH^d(X\times X)/alg,
 \end{eqnarray}
o\`u  $\Delta_X$ est la diagonale de $X$ et $Z$ est un cycle de $X\times X$ support\'{e} sur $D\times X$,
 $D\subsetneqq X$ \'{e}tant un ferm\'{e} propre de $X$.
 \end{lemm}
 \`A l'aide de ce lemme, on conclut de la fa\c{c}on suivante:
 Comme dans la preuve de la proposition
  \ref{lemmecorresp},
 on note qu'un cycle $Z$ comme ci-dessus agit trivialement sur $H^i(X,\mathcal{H}_X^d(A))$ pour tout $i$.
 On conclut donc de la d\'{e}composition (\ref{decomp9avril}) que l'on a :
 $$\Gamma_{\phi^{-1}\,*}\circ\Gamma_{\phi\,*}=Id:H^i(X, \mathcal{H}_X^d(A))\rightarrow H^i(X, \mathcal{H}_X^d(A)).$$
 Mais, pour les m\^{e}mes raisons, on a aussi
 $$\Gamma_{\phi\,*}\circ\Gamma_{\phi^{-1}\,*}=Id:H^i(Y, \mathcal{H}_Y^d(A))\rightarrow H^i(Y, \mathcal{H}_Y^d(A)),$$
 d'o\`u l'on conclut que $\Gamma_{\phi\,*}:H^i(X, \mathcal{H}_X^d(A))\rightarrow H^i(Y, \mathcal{H}_Y^d(A))$
 est un isomorphisme.
 \cqfd

\medskip

{\bf D\'{e}monstration du lemme \ref{lemme9avril}.} D'apr\`{e}s la th\'{e}orie de l'intersection raffin\'{e}e
de Fulton (voir \cite[Chap. 6]{fulton}),
et compte tenu de la d\'{e}finition de la composition des correspondances,
il suffit de montrer qu'il existe un ferm\'{e} alg\'{e}brique propre $D\subsetneqq X$ tel que
$p_{13}(p_{12}^{-1}(\Gamma_\phi)\cap p_{23}^{-1}(\Gamma_{\phi^{-1}})),$
o\`u les $p_{ij}$ sont les diff\'{e}rentes projections de $X\times Y\times X$ sur les
produits de deux de ses facteurs,
 co\"\i ncide (multiplicit\'{e}s comprises) avec $\Delta_X$ au-dessus de $U\times X$, o\`u $U=X\setminus D$.

Notons   $V$ l'ouvert de d\'{e}finition de $\phi$ et prenons pour ouvert $U\subset X$ l'ouvert
de $V$ d\'{e}fini comme $\phi^{-1}(W)$, o\`u $W\subset Y$ est l'ouvert de d\'{e}finition
de $\phi^{-1}$.
Alors au-dessus de $U$, la premi\`{e}re projection $p_1:X\times Y\rightarrow X$
induit un isomorphisme local
entre le  graphe de $\phi$  et  $X$. La vari\'{e}t\'{e} $U$ peut alors \^{e}tre vue
comme un ouvert $U'$ de ce graphe. Via la seconde projection $p_2:X\times Y\rightarrow Y$, $U'$ est par construction envoy\'{e}
dans l'ouvert  de d\'{e}finition de $\phi^{-1}$, et donc l'intersection
$p_{12}^{-1}(\Gamma_\phi)\cap p_{23}^{-1}(\Gamma_{\phi^{-1}})$ s'identifie au-dessus de $U$
\`{a} la sous-vari\'{e}t\'{e} lisse
$$\{(x,\phi(x),\phi^{-1}(\phi(x))=x),\,x\in V,\,\phi(x)\in W\}$$ de $X\times Y\times X$, ce qui donne le r\'{e}sultat.
\cqfd

Certains des invariants birationnels donn\'{e}s par la proposition \ref{invbir} \'{e}taient connus. Pour ce qui est de la non-trivialit\'{e}
de ces invariants,
 les propri\'{e}t\'{e}s de la suite spectrale de Bloch--Ogus (\S \ref{secrappelsbo}) donnent les informations suivantes.
Soit  $X$ une vari\'{e}t\'{e} connexe,  projective et lisse sur $\C$ de dimension $d$.

Pour $i=0$, on trouve la cohomologie non ramifi\'{e}e  de $X$ en degr\'{e} $d$.

 Pour $i=d$, on a
$H^d(X,\mathcal{H}_X^d(\Z))=CH^d(X)/alg=\Z$ et $H^d(X,\mathcal{H}_X^d(\Z/n))=CH^d(X)/n=\Z/n.$

 Pour $i=d-1$, on a  $H^{d-1}(X,\mathcal{H}_X^d(\Z))=H^{2d-1}(X(\C),\Z) \simeq H_{1}(X(\C),\Z)$ et pour un entier $n>1$, on a
$H^{d-1}(X,\mathcal{H}_X^d(\Z/n))=H^{2d-1}(X(\C),\Z/n)$, groupe dual de $H^1(X(\C),\Z/n)$,
qui est un invariant birationnel  bien connu, co\"{\i}ncidant avec la $n$-torsion du groupe de Picard de $X$.
Ces groupes sont non nuls si $H^1(X,\mathcal{O}_X)\neq 0$, soit encore $b_{1}\neq 0$.
La torsion du groupe de N\'{e}ron-Severi $NS(X)$, \'{e}gale \`{a} celle du groupe $H^2(X(\C),\Z)$, peut aussi
apporter une contribution.

Pour ce qui est des coefficients rationnels, ces r\'{e}sultats se g\'{e}n\'{e}ralisent de la fa\c{c}on suivante:
\begin{prop} \label{prop11juin}Si $H^i(X,\mathcal{O}_X)\not=0$, on a $H^{d-i}(X,\mathcal{H}^d_X(\mathbb{Q}))\not=0$.
\end{prop}
{\bf D\'{e}monstration.} Cela r\'{e}sulte du fait 
  (cf. \cite[p. 194]{blochogus}, \cite{gilletsoule}, \cite{paranjape})
 que l'aboutissement de la filtration de Leray
sur $H^{2d-i}_B(X(\C),\mathbb{Q})$ relative  \`{a} $\pi: X_{cl}\rightarrow X_{Zar}$ est la filtration par le coniveau $N$.
Si le terme extr\^{e}me
$E_{2}^{d-i,d}=H^{d-i}(X(\C),\mathcal{H}^d_X(\mathbb{Q}))$
est nul, on conclut donc que $H^{2d-i}_B(X(\C),\mathbb{Q})=N^{d-i+1}H^{2d-i}_B(X(\C),\mathbb{Q})$.
Ceci entra\^{\i}ne classiquement (cf. \cite{grothhodge}) que le coniveau de Hodge de la structure de Hodge sur $H^{2d-i}_B(X(\C),\mathbb{Q})$
est $\geq d-i+1$, et donc que $H^{d,d-i}(X):=H^{d-i}(X,\Omega_X^d)=H^{d-i}(X,K_X)=0$. Par la dualit\'{e} de Serre,
ceci contredit l'hypoth\`{e}se
$H^i(X,\mathcal{O}_X)\not=0$ .
\cqfd

 \subsection{$H^3$ non ramifi\'{e} \`{a} coefficients finis et  d\'{e}faut de la conjecture de Hodge enti\`{e}re en degr\'{e} 4}\label{sectionprincipale}

\begin{theo}\label{principal1}
Soit $X$ une vari\'{e}t\'{e}  connexe, projective et  lisse sur $\C$.
On a les suites exactes
$$ 0 \to H^3_{nr}(X, \Z(2))   /n  \to  H^3_{nr}(X,\mu_{n}^{\otimes 2}) \to Z^4(X)[n] \to 0,$$
$$ 0 \to H^3_{nr}(X, \Z(2))   \otimes \Q/\Z  \to  H^3_{nr}(X,\Q/\Z(2)) \to Z^4(X)_{\rm tors}\to 0.$$
\end{theo}

{\bf D\'{e}monstration.}
 On consid\`{e}re la suite spectrale
$$E_{2}^{pq} = H^p(X,\H^q_{X}(\Z(2))) \Longrightarrow H^*(X(\C),\Z(2))$$
associ\'{e}e  au morphisme de sites $X(\C) \to X_{Zar}$.
Cette suite est concentr\'{e}e dans le premier quadrant.
Comme $X$ est lisse, la conjecture de Gersten, connue   pour
cette cohomologie (\cite{blochogus}, voir  \cite{ctkh} et le \S \ref{secrappelsbo} ci-dessus),
implique que cette suite spectrale est  en outre  concentr\'{e}e dans le demi-quadrant  $p \leq q$.
La fl\`{e}che qui va de $H^2(X,\H^2_{X}(\Z(2)))$ dans $H^4(X(\C),\Z(2))$ a pour image le groupe
$H^4_{alg}(X(\C),\Z(2))$ des cycles alg\'{e}\-briques.
Cela r\'{e}sulte de l'identification
de la filtration qu'on obtient ainsi avec la filtration par la codimension du support
(\cite{paranjape}, \cite{gilletsoule}).

De la suite spectrale on tire alors une suite exacte
$$0 \to H^1(X,\H^3_{X}(\Z(2))) \to [H^4(X(\C),\Z(2))/H^4_{alg}(X(\C),\Z(2))]  \to H^0(X,\H^4_{X}(\Z(2))).$$

Pour $X$ de dimension quelconque,
le faisceau $\H^4_{X}(\Z(2))$ est sans torsion
(th\'{e}or\`{e}me \ref{principal0}).
On en conclut
$$H^1(X,\H^3_{X}(\Z(2)))[n]  \simeq [H^4(X(\C),\Z(2))/H^4_{alg}(X(\C),\Z(2))][n].$$
Comme on l'a rappel\'{e} dans l'introduction, $Z^4(X)[n]    = [H^4(X(\C),\Z(2))/H^4_{alg}(X(\C),\Z(2))][n].$
On a montr\'{e} au  th\'{e}or\`{e}me~\ref{principal0}
que la suite de faisceaux
$$  0 \to \H^3_{X}(\Z(2)) \buildrel{\times n}\over{\to} \H^3_{X}(\Z(2)) \to \H^3_{X}(\mu_{n}^{\otimes 2}) \to 0$$
est exacte. Elle
donne naissance \`{a} la suite exacte
$$  0 \to H^0(X,\H^3_{X}(\Z(2)))/n  \to H^0(X,\H^3_{X}(\mu_{n}^{\otimes 2})) \to H^1(X,\H^3_{X}(\Z(2)))[n] \to 0.$$
Ceci \'{e}tablit le premier \'{e}nonc\'{e} du th\'{e}or\`{e}me.

Le second s'obtient par passage \`{a} la limite directe dans les suites exactes de (i)
pour $n$ variable.
C'est une cons\'{e}quence du th\'{e}or\`{e}me de Merkur'ev et Suslin
 que les fl\`{e}ches de passage de la suite
pour $n$ \`{a} la suite pour $n.n'$ sont toutes injectives. \cqfd

\medskip

Rappelons que la conjecture de Hodge rationnelle vaut en degr\'{e} $2$ et $2d-2$
pour toute vari\'{e}t\'{e} projective et lisse de dimension $d$, donc en tout degr\'{e}  $2i$
pour une vari\'{e}t\'{e}
de dimension au plus~3. Plus g\'{e}n\'{e}ralement, comme rappel\'{e} dans l'introduction,
Bloch et Srinivas \cite[Thm. 1 (iv)]{blochsrinivas}  ont montr\'{e} que si le groupe $CH_{0}(X)$
est  \og repr\'{e}sent\'{e} \fg \, par le groupe $CH_{0}$ d'une vari\'{e}t\'{e} de dimension au plus~3
alors la conjecture de Hodge rationnelle vaut en  degr\'{e}~4 pour $X$.

On a donc l'\'{e}nonc\'{e} suivant :

\begin{theo}\label{principal2}
Soit $X$ connexe, projective et lisse. S'il existe une vari\'{e}t\'{e} $Y$ connexe, projective et lisse
de dimension au plus 3 et un morphisme $f :Y \to X$ tels que l'application induite
$f_{*} : CH_{0}(Y)\to CH_{0}(X)$ soit surjective, alors
le groupe $Z^4(X)$
est un groupe fini, et l'on a la suite exacte
$$ 0 \to H^3_{nr}(X, \Z(2))   \otimes \Q/\Z  \to  H^3_{nr}(X,\Q/\Z(2))   \to Z^4(X)\to 0.$$
\end{theo}

Le th\'{e}or\`{e}me suivant s'applique en particulier aux vari\'{e}t\'{e}s unirationnelles
et aux solides unir\'{e}gl\'{e}s.

\begin{theo} \label{blochsrinivasgriffnul}
Soit $X$ une vari\'{e}t\'{e} connexe, projective et lisse sur $\C$.
S'il existe une vari\'{e}t\'{e} $Y$, projective et lisse
de dimension  2
et un morphisme $f :Y \to X$ tels que l'application induite
$f_{*} : CH_{0}(Y)\to CH_{0}(X)$ soit surjective, alors
le groupe $Z^4(X)$ est fini, le groupe
  $H^3_{nr}(X,\Z(2))$ est nul  et
l'on  a un isomorphisme de groupes finis
$H^3_{nr}(X, \Q/\Z(2))  \oi Z^4(X).$
\end{theo}
{\bf D\'{e}monstration.}
D'apr\`{e}s  les th\'{e}or\`{e}mes \ref{principal1} et \ref{principal2}, il suffit de montrer l'annulation
de $H^3_{nr}(X,\Z(2))$. Sous l'hypoth\`{e}se ci-dessus, celle-ci est un cas particulier
de la proposition \ref{lemmecorresp}  (i).
\cqfd

\begin{rema} {\rm  Dans \cite{barbieriviale},  Barbieri-Viale consid\`{e}re
les vari\'{e}t\'{e}s projectives et  lisses de dimension~3.
 Il montre (Teorema 5.2)  que
si $H^0(X,\H^3_{X}(\Z(2)))=0$ (ce qui est en particulier satisfait si
$X$ est unirationnelle) alors pour tout entier $n>0$,
$$H^0(X,\H^3_{X}(\mu_{n}^{\otimes 2})) \oi Z^4(X)[n].$$
Les  d\'{e}monstrations de ce paragraphe
\'{e}tendent  directement sa d\'{e}monstration. }
\end{rema}

\subsection{D\'{e}faut de la conjecture de Hodge enti\`{e}re  en degr\'{e} $2d-2$.}

Le cas particulier $\dim(X)=3$ du th\'{e}or\`{e}me \ref{principal1} porte sur les
cycles de dimension 1. On peut aussi \'{e}tablir un \'{e}nonc\'{e} sur les cycles de dimension 1
sur les vari\'{e}t\'{e}s de dimension quelconque.

\begin{theo}\label{principal3}
Soit $X$ une vari\'{e}t\'{e} connexe, projective et lisse sur $\C$,   de dimension $d$.
On a une suite exacte
$$0 \to   H^{d-3}(X, \H^d_{X}(\Z(d-1)))\otimes \Q/\Z\to    H^{d-3}(X, \H^d_{X}(\Q/\Z(d-1)))  \to Z^{2d-2}(X) \to 0,$$
o\`u le groupe $Z^{2d-2}(X)$ est fini. \cqfd
\end{theo}

{\bf D\'{e}monstration.}
Comme montr\'{e} au  th\'{e}or\`{e}me \ref{principal0},
on a la suite exacte courte de faisceaux
$$  0 \to \H^d_{X}(\Z(d-1)) \buildrel{\times n}\over{\to} \H^d_{X}(\Z(d-1)) \to \H^d_{X}(\mu_{n}^{\otimes d-1}) \to 0.$$
Pour tout entier $n>0$, on a donc une suite exacte
$$0 \to   H^{d-3}(X, \H^d_{X}(\Z(d-1)))/n\to    H^{d-3}(X,  \H^d_{X}(\mu_{n} ^{\otimes d-1  }))  \to H^{d-2}(X,\H^d_{X}(\Z(d-1)))[n] \to 0.$$
Par passage \`{a} la limite sur $n$, on a une suite exacte
$$0 \to   H^{d-3}(X, \H^d_{X}(\Z(d-1))) \otimes  \Q/\Z \to    H^{d-3}(X, \H^d_{X}(\Q/\Z(d-1)))  \to H^{d-2}(X,\H^d_{X}(\Z(d-1)))_{\rm tors} \to 0.$$

Les termes $E_{2}^{pq}$ de la suite spectrale de Bloch--Ogus pour le faisceau $\Z(d-1)$  sont nuls
en dehors du triangle $0\leq p \leq q \leq d$ (voir le \S \ref{secrappelsbo}). On obtient ainsi la suite exacte
$$H^{d-3}(X,\H^d_{X}(\Z(d-1))) \to CH^{d-1}(X)/{alg} \to H^{2d-2}(X(\C),\Z(d-1)) \to H^{d-2}(X,\H^d_{X}(\Z(d-1))) \to 0.$$
Les arguments donn\'{e}s dans l'introduction montrent que le groupe $Z^{2d-2}(X)$
 est fini,
et qu'il s'identifie \`{a} la torsion du conoyau de  l'application $CH^{d-1}(X)/{alg} \to H^{2d-2}(X(\C),\Z(d-1))$,
et donc,  d'apr\`{e}s ce qui pr\'{e}c\`{e}de, \`{a} la torsion, finie, du groupe de type fini $H^{d-2}(X,\H^d_{X}(\Z(d-1))) $.
\cqfd

\begin{coro}\label{coronomme9avril}
Soit $X$ une vari\'{e}t\'{e} connexe, projective et lisse sur $\C$,  de dimension $d$.
Si le groupe $CH_{0}(X)$ est support\'{e} sur une surface, alors on a un isomorphisme
de groupes finis
 $$  H^{d-3}(X, \H^d_{X}(\Q/\Z(d-1)))  \oi Z^{2d-2}(X).$$
 \end{coro}

 {\bf D\'{e}monstration.}
D'apr\`{e}s la proposition \ref{lemmecorresp}, le groupe
  $H^{d-3}(X, \H^d_{X}(\Z(d-1)))$ est, sous l'hypoth\`{e}se ci-dessus,  de torsion.
On conclut par une application du
 th\'{e}or\`{e}me
 \ref{principal3}.
\cqfd

\begin{rema} \label{rema11juin} {\rm
  Comme il est remarqu\'{e} dans \cite[Lemme 1]{soulevoisin}, le groupe
 $Z^{2d-2}(X)$ est un invariant birationnel des vari\'{e}t\'{e}s projectives, lisses, connexes sur $\C$.
 Le th\'{e}or\`{e}me \ref{principal3} d\'{e}crit ce groupe \`{a} l'aide de la cohomologie
 de $X$ \`{a} valeurs dans les faisceaux $\mathcal{H}_X^d$, $d=\dim\,X$, pour des coefficients ad\'{e}quats.
 La proposition \ref{invbir} montre que ces groupes sont aussi des invariants birationnels
 et la Proposition \ref{prop11juin} donne
 une condition suffisante 
 pour leur non trivialit\'{e} \`{a} coefficients dans $\mathbb{Q}$.

Quant aux coefficients de torsion, pour tout entier $d\geq 3$, des exemples de Koll\'ar (voir \S  \ref{subseckollar} ci-apr\`{e}s,
 \cite{kollar}, \cite[\S 2]{soulevoisin}, \cite[Thm. 1]{voisinuniruled}) montrent que
 le groupe $Z^{2d-2}(X)$ n'est pas forc\'{e}ment nul.
 
 Pour $d\geq 4$, ceci donne des exemples o\`u les
invariants  birationnels $H^{d-2}(X,\H^d_{X}(\Z))$ et  $H^{d-3}(X,\H^d_{X}(\Z/n))$ sont non nuls,   mais o\`{u} cependant
les groupes $H^i(X,\mathcal{O}_X)$  intervenant dans la Proposition \ref{prop11juin} sont nuls.}

 \end{rema}

\section{Cohomologie non ramifi\'{e}e \`{a} coefficients entiers en degr\'{e} 3,
cohomologie coh\'{e}rente en degr\'{e} 3,
groupe de Griffiths en degr\'{e}~2}\label{griffiths}

Dans toute cette section on consid\`{e}re une vari\'{e}t\'{e} $X$ connexe, projective et  lisse  sur $\C$.
On se propose  de discuter un lien  conjectural entre $H^3_{nr}(X,\Z)$
 et $H^3(X,\mathcal{O}_X)$.

 \subsection{La situation en degr\'{e}s 1 et 2}\label{groupebrauer}

 Les groupes  $H^1_{nr}$ et $H^2_{nr}$ sont compris depuis longtemps
(\cite[ II.3 et  III.8]{grothendieck}).
Le lecteur v\'{e}rifiera  que l'on
a $H^1_{nr}(X,\Z)=0$ si et seulement si $H^1(X,\O_{X})=0$.

\medskip

Notons ${\rm Br}(X)=H^2_{\et}(X,\mathbb{G}_{m})$ le groupe de Brauer--Grothendieck.
Pour tout entier $n>0$, on~a
$$H^2_{nr}(X,\mu_{n} ) \oi {\rm Br}(X)[n].$$
On a la suite exacte
$$0 \to NS(X) \to H^2(X(\C),\Z(1)) \to H^2_{nr}(X,\Z(1)) \to 0,$$
o\`u $NS(X)$ est le groupe de N\'{e}ron--Severi de $X$. La
fl\`{e}che $NS(X) \to H^2(X(\C),\Z(1))$ induit un isomorphisme sur les
groupes de torsion.  Si l'on note $\rho$ le rang du $\Q$-vectoriel $NS(X)\otimes \Q$
et $b_{2}$ le deuxi\`{e}me nombre de Betti de $X(\C)$,
on a un isomorphisme de groupes ab\'{e}liens
$$ H^2_{nr}(X,\Z(1)) \simeq \Z^{b_{2}-\rho}.$$
Enfin on a une  suite exacte
$$ 0 \to (\Q/\Z)^{b_{2}-\rho} \to H^2_{nr}(X,\Q/\Z(1)) \to H^3(X(\C),\Z(1))_{\rm tors} \to 0.$$

Par ailleurs la th\'{e}orie de Hodge montre que l'on a $b_{2}-\rho=0$ si et seulement si $H^2(X,\O_{X})=0$.
En d'autres termes, $H^2_{nr}(X,\Z(1))=0$ si et seulement si $H^2(X,\O_{X})=0$.

\subsection{La situation en degr\'{e} 3}

Au paragraphe~\ref{secrappelsbo}, on vu  l'isomorphisme
 $$ H^3_{nr}(X, \Z(2))/ \Im[ H^3(X(\C),\Z(2))] \oi  {\rm Griff}^2(X).$$

 Le th\'{e}or\`{e}me \ref{principal0} assure
  que le groupe $ H^3_{nr}(X, \Z(2))$ est sans torsion.
 Comme le groupe $H^3(X(\C),\Z(2))$ est de type fini, on en d\'{e}duit imm\'{e}diatement :
  \begin{prop}
 Pour tout $l$ premier, le sous-groupe de torsion $l$-primaire  ${\rm Griff}^2(X)\{l\}$
 est de cotype fini, son sous-groupe divisible maximal \'{e}tant de corang au plus le troisi\`{e}me nombre de Betti $b_{3}$.
  \end{prop}

 Le r\'{e}sultat suivant  montre  que le groupe
 $H^3_{nr}(X,\Z(2)    )/\Im[ H^3(X(\C),\Z(2))]$ n'est pas en g\'{e}n\'{e}ral de torsion.
 \begin{theo} (Griffiths  \cite{griffiths}) Soit $X\subset \mathbb{P}^4$ une hypersurface tr\`{e}s g\'{e}n\'{e}rale   de degr\'{e}
 $5$. Alors le groupe ${\rm Griff}^2(X)$ n'est pas de torsion.
 \end{theo}
 Ce r\'{e}sultat a \'{e}t\'{e} am\'{e}lior\'{e} de fa\c{c}on frappante par Clemens \cite{clemens} (voir aussi \cite{voisinduke}
 pour le cas des vari\'{e}t\'{e}s de Calabi-Yau de dimension $3$ arbitraires) qui montre que dans le cas
 consid\'{e}r\'{e} par Griffiths,  le
  $\mathbb{Q}$-espace vectoriel ${\rm Griff}^2(X)\otimes\mathbb{Q} \simeq
  H^3_{nr}(X,\mathbb{Q}(2))/\Im[H^3(X(\C),\mathbb{Q}(2))]$
  n'est pas de dimension finie.

  \medskip

  Par ailleurs   le groupe
  $ H^3_{nr}(X, \Z(2))/ \Im[ H^3(X(\C),\Z(2))] \simeq  {\rm Griff}^2(X)$ n'est pas en g\'{e}n\'{e}ral divisible :
  \begin{theo}(Bloch--Esnault \cite{blochesnault})\label{blochesnault}
  Il existe des intersections compl\`{e}tes lisses $X \subset {\mathbb P}^5_{\C}$
  de dimension 3
  pour lesquelles le groupe de Griffiths  ${\rm Griff}^2(X)$ n'est pas divisible.
  \end{theo}

 Ceci a \'{e}t\'{e} rafffin\'{e} par C. Schoen \cite{schoen}), qui montre qu'il existe
  des vari\'{e}t\'{e}s $X$ projectives lisses connexes
  de dimension 3 sur $\C$ et des entiers $n>1$ tels que les quotients
  ${\rm Griff}^2(X)/n$ et
  $CH^2(X)/n= [CH^2(X)/{alg}]/n$ soient  infinis.

 Les d\'{e}monstrations de ces r\'{e}sultats sont arithm\'{e}tiques, les vari\'{e}t\'{e}s
 consid\'{e}r\'{e}es sont d\'{e}finies sur un corps de nombres.

  \medskip

Supposons maintenant $H^3(X,\mathcal{O}_X)=0$. La conjecture de Hodge g\'{e}n\'{e}ralis\'{e}e
(par  Grothendieck \cite{grothhodge})
pr\'{e}dit alors
que la cohomologie $H^3(X,\mathbb{Q})$ est de coniveau g\'{e}om\'{e}trique $1$, de sorte qu'il devrait exister
une vari\'{e}t\'{e} projective complexe lisse $Y$ de dimension $n-1$, o\`{u} $n=\dim\,X$, et
un morphisme $j: Y\rightarrow X$, tel que
$j_*:H^1(Y,\mathbb{Q})\rightarrow H^3(X,\mathbb{Q})$ soit surjectif.
Il en r\'{e}sulte imm\'{e}diatement que
le morphisme induit par $j_*$ entre les jacobiennes interm\'{e}diaires
$J^{1}(Y)=\Pic^0(Y)$ et $J^2(X)=H^3(X,\C)/(F^2H^3(X,\C)\oplus H^3(X,\Z)/{\rm tors})$
est surjectif.
Comme l'application d'Abel-Jacobi $AJ_Y^1:CH^1(Y)_{hom}\rightarrow J^1(Y)$ est surjective, il en r\'{e}sulte aussi
que l'application d'Abel-Jacobi $AJ_X^2:CH^2(X)_{hom}\rightarrow J^2(X)$  de $X$ envoie
surjectivement $j_*CH^1(Y)_{hom}$ vers $J^2(X)$.
Or $j_*CH^1_{hom}(Y)$ est contenu dans le groupe $CH^2(X)_{alg}$
des cycles alg\'{e}briquement \'{e}quivalents \`{a} $0$
de $X$.

Soit maintenant $z\in CH^2(X)_{hom}$. D'apr\`{e}s ce qui pr\'{e}c\`{e}de, il existe $z'\in CH^2(X)_{alg}$ tel que
$AJ_X^2(z-z')=0$. Or on a la conjecture suivante, due \`{a} Nori \cite{nori}:
\begin{conj} \label{noriconj}Les cycles de codimension $2$ homologues \`{a} $0$ et annul\'{e}s par l'application d'Abel-Jacobi
sont de torsion modulo l'\'{e}quivalence alg\'{e}brique.
\end{conj}

On conclut donc que le groupe $H^3_{nr}(X,\Z(2))/\Im[ H^3(X(\C),\Z (2)) ] \simeq {\rm Griff}^2(X)$ est conjecturalement  de torsion si
$H^3(X,\mathcal{O}_X)=0$.

 Par ailleurs
la conjecture de Hodge g\'{e}n\'{e}ralis\'{e}e
 (\cite{grothhodge}) et l'hypoth\`{e}se $H^3(X,\mathcal{O}_X)=0$
 impliquent l'existence d'un ouvert de Zariski non vide $U \subset X$ tel que la restriction
 $H^3(X(\C),\Q(2)) \to H^3(U(\C),\Q(2))$ soit nulle. Il existe donc un entier $n>0$
 tel que la restriction $H^3(X(\C),\Z(2)) \to H^3(U(\C),\Z(2))$ soit annul\'{e}e par $n$.
 A fortiori la fl\`{e}che $H^3(X(\C),\Z(2))  \to H^3_{B}(\C(X),\Z(2))$ est-elle annul\'{e}e par $n$.
   La restriction  $H^3(X(\C),\Z(2))  \to H^3_{B}(\C(X),\Z(2))$
 se factorise par la fl\`{e}che  $H^0(X,\H^3(\Z(2))) \rightarrow H^3_{B}(\C(X),\Z(2))$
 qui est injective par la th\'{e}orie de Bloch--Ogus
 (th\'{e}or\`{e}me  \ref{boresolution}).
 Ainsi
  la fl\`{e}che  $H^3(X(\C),\Z(2)) \to H^0(X,\H^3(\Z(2))) $ a son image annul\'{e}e par $n>0$.
Mais  $H^0(X,\H^3(\Z(2)))=H^3_{nr}(X,\Z(2))$ est sans torsion
  (\cite{blochsrinivas}; thm. \ref{principal0} ci-dessus).
La fl\`{e}che $H^3(X(\C),\Z(2)) \to H^3_{nr}(X,\Z(2))$ est donc nulle.
Ainsi  $H^3_{nr}(X,\Z(2)) \simeq {\rm Griff}^2(X)$, mais le groupe de gauche est sans torsion
et celui de droite  conjecturalement  de torsion.

La conjecture de Nori, combin\'{e}e avec la conjecture de Hodge g\'{e}n\'{e}ralis\'{e}e,
et le th\'{e}or\`{e}me \ref{principal0}
  permettraient  donc d'\'{e}tablir la conjecture suivante.

\begin{conj} Si $H^3(X,\mathcal{O}_X)=0$, le groupe $H^3_{nr}(X,\Z (2))$ est nul,
 ainsi que le groupe
${\rm Griff}^2(X)$, et l'on a un isomorphisme de groupes finis
$H^3_{nr}(X, \Q/\Z(2))
  \oi Z^4(X).$
\end{conj}

Cet \'{e}nonc\'{e} conjectural, faisant suite aux \'{e}nonc\'{e}s du paragraphe \ref{groupebrauer},
sugg\`{e}re le probl\`{e}me intriguant suivant:

Pour $ i\geq 4$,  existe-t-il des vari\'{e}t\'{e}s $X$ avec $ H^i(X,\mathcal{O}_X)=0$   pour lesquelles  $H^i_{nr}(X,\Z(i-1))\neq0$ ?

Notons qu'on a par ailleurs la conjecture suivante  (pour laquelle les   hypoth\`{e}ses   ont une structure diff\'{e}rente) :
\begin{conj} Si $H^j(X,\mathcal{O}_X)=0$ pour $j\geq i$, le groupe $H^i_{nr}(X,\Z (i-1))$ est nul.
\end{conj}
Cette derni\`{e}re conjecture est en effet entra\^{\i}n\'{e}e par la conjecture de Bloch g\'{e}n\'{e}ralis\'{e}e
(c'est-\`{a}-dire, la r\'{e}ciproque  du th\'{e}or\`{e}me de Mumford g\'{e}n\'{e}ralis\'{e} \cite[Th\'{e}or\`{e}me 22.17]{voisinlivre}), qui dit que sous ces hypoth\`{e}ses, il existe un ferm\'{e} alg\'{e}brique
$Y\subset X$ de dimension $\leq i-1$ tel que la fl\`{e}che $CH_0(Y)\rightarrow CH_0(X)$ soit surjective, combin\'{e}e avec
la proposition \ref{lemmecorresp} (i).

\section{Vari\'{e}t\'{e}s avec $ H^3_{nr}(X,\mu_{n}^{\otimes 2}) \neq 0$ ou $Z^4(X)_{\rm tors}\neq 0$.}

Soit $X$ une vari\'{e}t\'{e} connexe, projective et lisse sur $\C$ et $n>0$ un entier.  On a le diagramme
commutatif de suites exactes :

$$\begin{array}{ccccccccccccccccccccccccccccccccccccccccccccccccccc}
0 & \to  &  H^3(X(\C),\Z(2))/n   &  \to  &  H^3_{\et}(X,\mu_{n}^{\otimes 2})   &  \to   &  H^4(X(\C),\Z(2))[n] &  \to   & 0 \cr
 && \downarrow  &&    \downarrow &  &  \downarrow && \cr
 0  &  \to  &  H^3_{nr}(X, \Z(2)))/n  &  \to  & H^3_{nr}(X,\mu_{n}^{\otimes 2})   & \to   & Z^4(X)[n]  &  \to  &  0 \cr
  && \downarrow &&    \downarrow &  &    && \cr
    &   &  {\rm Griff}^2(X)/n  &  \to  &  \Ker[CH^2(X)/n \to H^4_{\et}(X,\mu_{n}^{\otimes 2})]   &     &   &     & \cr
    && \downarrow &&    \downarrow &  &    && \cr
     &   &  0  &    &  0  &     &   &     & \cr
  \end{array} $$

La premi\`{e}re suite exacte horizontale vient de la longue suite exacte de groupes de cohomologie de Betti associ\'{e}e \`{a}
la multiplication par $n$ sur le faisceau constant $\Z(2)$; on a utilis\'{e} l'identification $H^3_{\et}(X,\mu_{n}^{\otimes 2}) \oi
H^3(X(\C),\mu_{n}^{\otimes 2})$.
La suite exacte horizontale m\'{e}diane est donn\'{e}e par le th\'{e}or\`{e}me \ref{principal1}. Les deux  fl\`{e}ches verticales sup\'{e}rieures de gauche sont
les fl\`{e}ches \'{e}videntes, la fl\`{e}che verticale sup\'{e}rieure droite  est induite par le reste du diagramme.
La suite verticale de gauche est induite par le th\'{e}or\`{e}me \ref{theogriffblochogus}.
La suite verticale m\'{e}diane est celle donn\'{e}e au paragraphe \ref{secrappelblochkato}, apr\`{e}s le th\'{e}or\`{e}me \ref{boresolutionfinie}.

Dans ce paragraphe, on passe en revue un certain nombre d'exemples de vari\'{e}t\'{e}s $X$ construites dans la litt\'{e}rature, pour lesquelles
certains des groupes apparaissant dans les deux lignes inf\'{e}\-rieures du diagramme ci-dessus sont non nuls.

 Comme on le voit, deux groupes de types tr\`{e}s diff\'{e}rents peuvent
contribuer \`{a} la non nullit\'{e} de  $H^3_{nr}(X,\mu_{n}^{\otimes 2})$ : d'une part le groupe $H^3_{nr}(X, \Z(2))/n $, d'autre part le groupe
$ Z^4(X)[n] $. Nous commentons \'{e}galement ces exemples du
point de vue du probl\`{e}me  de d\'{e}formation ou de  sp\'{e}cialisation qui est d\'{e}crit pr\'{e}cis\'{e}ment dans
la sous-section \ref{subsecspecialisation}.

\subsection{\label{subsecspecialisation}
D\'{e}formation et sp\'{e}cialisation du groupe $Z^{2i}$}

Lorsque l'on a un exemple de vari\'{e}t\'{e} projective
et lisse $X$ avec $Z^{2i}(X) \neq 0$, on s'int\'{e}resse \`{a} la question de la stabilit\'{e} de cet exemple par
 d\'{e}formation, petite ou globale. Voici comment on peut
 donner un sens pr\'{e}cis \`{a} ce probl\`{e}me
 (nous renvoyons \`{a} \cite{voisinhodgeloci} pour plus de d\'{e}tails).

 On a les inclusions $H^{2i}_{alg}(X(\C),\Z(i)) \subset Hdg^{2i}(X,\Z) \subset H^{2i}(X(\C),\Z(i))$.
Une famille de d\'{e}formations  de $X$ est une  famille projective et lisse $\mathcal{X} \to T$
au-dessus d'un $\C$-sch\'{e}ma $T$    connexe \'{e}quip\'{e} d'un point marqu\'{e} $0$, de fibre sp\'{e}ciale $\mathcal{X}_{0} = X$.
Soit $f:\mathcal{X}(\mathbb{C}) \to T(\mathbb{C})$ le morphisme correspondant de vari\'{e}t\'{e}s complexes.
D'apr\`{e}s le lemme d'Ehresmann \cite[9.1]{voisinlivre}, $f$ est une fibration topologique au-dessus de $T(\mathbb{C})$.
Il en r\'{e}sulte que $R^{2i}f_*\mathbb{Z}(i)$ est un syst\`{e}me local sur $T(\mathbb{C})$.

Pour  $t,\,s \in  T(\C)$, on a alors une identification
 $\rho_{t,s}  : H^{2i}(\mathcal{X}_t(\C),\Z(i)) \oi H^{2i}(\mathcal{X}_{s}(\C),\Z(i))$, canoniquement
 d\'{e}termin\'{e}e par le choix d'un chemin continu $[0,1]\rightarrow T(\C)$ de $t$ \`{a} $s$.

 Partant d'une classe $\alpha$ dans $H^{2i}(\mathcal{X}_t(\mathbb{C}),\Z(i))$, on se pose deux questions :

(a) pour quels $s \in T(\C)$ a-t-on $\rho_{t,s}(\alpha) \in Hdg^{2i}(\mathcal{X}_{s},\Z)$ (pour un choix ad\'{e}quat de chemin de
$t$ \`{a} $s$)?

(b)  pour quels $s \in T(\C)$ a-t-on $\rho_{t,s}(\alpha) \in H^{2i}_{alg}(\mathcal{X}_{s}(\mathbb{C}),\Z(i))$ (pour un choix ad\'{e}quat de chemin de
$t$ \`{a} $s$)?

Une cons\'{e}quence tr\`{e}s simple de l'existence, de la projectivit\'{e}, et de la d\'{e}nombrabilit\'{e} des sch\'{e}mas de Hilbert relatifs
est  l'\'{e}nonc\'{e} suivant (cf. \cite[1.2]{voisinhodgeloci}) :

\begin{lemm} \label{lemmedefalg}  L'ensemble $T_{\alpha,alg}$ des $s \in T(\C)$
 tels que  $\rho_{t,s}(\alpha) \in H^{2i}_{alg}(\mathcal{X}_{s}(\mathbb{C}),\Z(i))$ pour un choix ad\'{e}quat de chemin de
$t$ \`{a} $s$ est une union d\'{e}nombrable de ferm\'{e}s alg\'{e}briques  de $T$.
\end{lemm}
 Soit $T_{alg}$ l'union (d\'{e}nombrable) des ensembles $T_{\alpha,alg}$, prise
 sur les classes $\alpha$ telles que  $T_{\alpha,alg}\not=T(\mathbb{C})$. Cet ensemble ne d\'{e}pend pas de $t$, comme il r\'{e}sulte des propri\'{e}t\'{e}s formelles du transport parall\`{e}le.

Infiniment plus difficile est l'\'{e}nonc\'{e} correspondant pour les classes de Hodge, montr\'{e}
par Cattani, Deligne et Kaplan \cite{CDK}:
\begin{theo} \label{CDK}   L'ensemble $T_{\alpha,Hdg}$ des $s \in T(\C)$ tels que  $\rho_{t,s}(\alpha) \in Hdg^{2i}(\mathcal{X}_{s},\Z)$ pour un choix ad\'{e}quat de chemin de
$t$ \`{a} $s$ est un ferm\'{e} alg\'{e}brique de $T$.
\end{theo}
Soit $T_{Hdg}$ l'union (d\'{e}nombrable) des ferm\'{e}s $T_{\alpha,Hdg}$, prise
 sur les classes $\alpha$ telles que  $T_{\alpha,Hdg}\not=T(\mathbb{C})$. Comme pr\'{e}c\'{e}demment, cet ensemble ne d\'{e}pend pas de $t$.

On en d\'{e}duit  imm\'{e}diatement
le r\'{e}sultat suivant:
\begin{coro} i) Si $t\in T(\C)$ est tr\`{e}s g\'{e}n\'{e}ral (et plus pr\'{e}cis\'{e}ment, en dehors de
$T_{alg}$), pour toute classe
$\alpha\in H^{2i}_{alg}(\mathcal{X}_t(\mathbb{C}),\mathbb{Z}(i))$, pour tout $s\in T(\mathbb{C})$, et
pour tout choix  de chemin de
$t$ \`{a} $s$, on a
$\rho_{t,s}(\alpha) \in H^{2i}_{alg}(\mathcal{X}_{s}(\mathbb{C}),\Z(i))$.

ii) Si $t\in T(\C)$ est tr\`{e}s g\'{e}n\'{e}ral (et plus pr\'{e}cis\'{e}ment, en dehors de
$T_{Hdg}$), pour toute classe
$\alpha\in Hdg^{2i}(\mathcal{X}_t,\mathbb{Z})$, pour tout $s\in T(\mathbb{C})$, et pour tout choix  de chemin de
$t$ \`{a} $s$, on a
$\rho_{t,s}(\alpha) \in Hdg^{2i}(\mathcal{X}_{s},\Z)$.
\end{coro}
Ceci montre que pour $t$ tr\`{e}s g\'{e}n\'{e}ral, et pour tout $s\in T(\mathbb{C})$, on dispose d'une application
\begin{eqnarray}\label{flechedespe}\rho_{t,s}^i:Z^{2i}(\mathcal{X}_t)\rightarrow Z^{2i}(\mathcal{X}_s)
\end{eqnarray}
d\'{e}termin\'{e}e par le choix d'un chemin de $t$ \`{a} $s$ et donc
bien d\'{e}finie \`{a} composition pr\`{e}s \`{a} droite et \`{a} gauche par l'action de monodromie (l'action \`{a} gauche pr\'{e}servant
$Hdg^{2i}(\mathcal{X}_{t},\Z)$ et $H^{2i}_{alg}(\mathcal{X}_{t}(\mathbb{C}),\Z(i))$ et donc agissant sur $Z^{2i}(\mathcal{X}_t)$).

 Une fa\c{c}on de formuler le probl\`{e}me de d\'{e}formation consiste \`{a} examiner les questions suivantes:

  (c)
Pour   $t$  tr\`{e}s g\'{e}n\'{e}ral dans $T(\mathbb{C})$, l'application
$\rho_{t,s}^i$ est-elle injective pour tout $s$ ?

  (d)
Pour   $t$   tr\`{e}s g\'{e}n\'{e}ral dans $T(\mathbb{C})$,
l'application
$\rho_{t,s}^i$ est-elle surjective pour tout $s$?

 Les questions (c) et (d) ne sont en fait int\'{e}ressantes que pour les entiers $i$ pour lesquels les groupes
$Z^{2i}$ sont des  invariants  birationnels non triviaux, c'est-\`{a}-dire $i=2$ et $i=n-1,\,n={\rm dim}\,\mathcal{X}_t$.
Elles sont motiv\'{e}es par l'application potentielle \`{a} l'\'{e}tude ou la caract\'{e}risation  des
vari\'{e}t\'{e}s rationnelles: on ignore en effet si, pour une famille $\mathcal{X}\rightarrow B$ de vari\'{e}t\'{e}s projectives lisses sur un corps alg\'{e}briquement clos $K$,
la propri\'{e}t\'{e} de rationalit\'{e} des fibres $\mathcal{X}_b$ est une propri\'{e}t\'{e} ouverte, ou ferm\'{e}e, sur $B(K)$ pour la topologie de Zariski.
Un seul fait est \'{e}vident, par consid\'{e}ration de la d\'{e}nombrabilit\'{e} et de la projectivit\'{e} des sch\'{e}mas de Hilbert relatifs de
$\mathbb{P}^n\times \mathcal{X}\rightarrow B$ :
 cette propri\'{e}t\'{e} est satisfaite sur une union d\'{e}nombrable
de sous-ensembles localement ferm\'{e}s pour la topologie de Zariski.

 La  question (c) concernant l'injectivit\'{e} peut aussi \^{e}tre formul\'{e}e pour une classe $\alpha\in Z^{2i}(\mathcal{X}_t)$ non nulle donn\'{e}e: a-t-on  $\rho^i_{t,s}(\alpha)\not=0$ dans $Z^{2i}(\mathcal{X}_s)$ pour tout $s$?
On verra au paragraphe \ref{subseckollar}  que la question (c) a en g\'{e}n\'{e}ral une r\'{e}ponse n\'{e}gative en degr\'{e} $i=2$ et pour des vari\'{e}t\'{e}s de dimension $3$.

En ce qui concerne la question (d), la r\'{e}ponse est n\'{e}gative pour  $i\not=2,\,n-1$, comme le montre l'exemple suivant:
Partons  de n'importe quelle vari\'{e}t\'{e} $Y$ pour laquelle
  $Z^4(Y)$  est un groupe cyclique d'ordre $d\geq 4$ (voir paragraphe \ref{subseckollar}), introduisons
la famille $\mathcal{S}\rightarrow B$ des surfaces lisses de degr\'{e} $d$ dans $\mathbb{P}^3$, et  consid\'{e}rons
la famille $\mathcal{X}=\mathcal{S}\times Y\rightarrow B$. Soit $\alpha\in Hdg^4(Y,\mathbb{Z} )$ une classe de Hodge enti\`{e}re non alg\'{e}brique.
Choisissons $s\in B$ tel que la surface $\mathcal{S}_s$ contienne une droite $\Delta$, et soit $\delta:=[\Delta]\in Hdg^2(\mathcal{S}_s,\mathbb{Z})\subset H^2(\mathcal{S}_s(\mathbb{C}),\mathbb{Z})$. Alors la classe $\alpha':=pr_1^*\delta\cup pr_2^*\alpha\in Hdg^6(\mathcal{X}_s,\mathbb{Z})$ n'est pas alg\'{e}brique,
car sinon $\alpha=pr_{2*}(pr_1^*c_1(\mathcal{O}_{\mathcal{S}_s}(1))\cup \alpha')$ le serait aussi.
Cette classe fournit donc une classe non nulle $\overline{\alpha'}$ dans $Z^6(\mathcal{X}_s)$, et on a :
\begin{lemm} La classe $\overline{\alpha'}$  n'est pas dans l'image de
 l'application
$\rho^3_{t,s}$ pour un point tr\`{e}s g\'{e}n\'{e}ral $t$ de $B$.
\end{lemm}
{\bf D\'{e}monstration.} En effet, supposons par l'absurde que $\overline{\alpha'}$ soit dans l'image
de $\rho^3_{t,s}$. De fa\c{c}on \'{e}quivalente, on peut \'{e}crire  $\alpha'=\alpha'_1+\alpha'_2$, o\`{u}
$\alpha'_1$ est une classe alg\'{e}brique sur $\mathcal{S}_s$ et $\alpha'_2\in Hdg^6(\mathcal{X}_s,\mathbb{Z})$ est une classe de Hodge qui est obtenue par transport parall\`{e}le d'une classe de Hodge
au point tr\`{e}s g\'{e}n\'{e}ral $t$. Le  th\'{e}or\`{e}me de Noether-Lefschetz
\cite[15.3]{voisinlivre}  dit que les classes de Hodge sur $\mathcal{S}_t$, pour $t$ tr\`{e}s g\'{e}n\'{e}ral dans $B$, sont r\'{e}duites \`{a} la classe
$c_1(\mathcal{O}_{\mathcal{S}_t}(1))$.  On  en d\'{e}duit que $\alpha'_2\in H^6(\mathcal{X}_s(\mathbb{C}),\mathbb{Z})$ a sa composante de K\"{u}nneth de type $(2,4)$ de la forme $pr_1^*c_1(\mathcal{O}_{\mathcal{S}_s}(1))\cup pr_2^*\beta$ pour une certaine classe
$\beta\in Hdg^4(Y,\mathbb{Z})$, puis que
$$\alpha=pr_{2*}(pr_1^*c_1(\mathcal{O}_{\mathcal{S}_s}(1))\cup \alpha')=pr_{2*}(pr_1^*c_1(\mathcal{O}_{\mathcal{S}_s}(1))\cup (\alpha'_1+\alpha'_2))$$
est \'{e}gal modulo les classes alg\'{e}briques sur $Y$ \`{a} ${\rm deg}\,(c_1(\mathcal{O}_{\mathcal{S}_s}(1))^2) \beta=d\beta$. Comme le groupe
$Z^4(Y)$ est annul\'{e} par $d$, le dernier terme est nul dans $Z^4(Y)$, ce qui contredit le fait que $\alpha$ est non nul dans $Z^4(Y)$.

\cqfd

Par contre pour $i=2$ ou $i=n-1$,
 on ne conna\^{\i}t pas
  de contre-exemple \`{a} la question (d).

\subsection{Exemples d'Atiyah--Hirzebruch}

On sait  depuis Atiyah et Hirzebruch \cite{atiyahhirzebruch}
que les groupes $Z^{2i}(X)$ ne sont pas forc\'{e}ment triviaux pour $i \neq 1$.
Pour tout premier $l$, leur m\'{e}thode permet en particulier de fabriquer une vari\'{e}t\'{e}
$X$ avec $Z^4(X)[l] \neq 0$. Pour ces vari\'{e}t\'{e}s, on a donc $H^3_{nr}(X,\Z/l) \neq 0$.
L'exemple de dimension minimale chez eux, pour $l=2$, est une vari\'{e}t\'{e} de dimension 7.
Les exemples de Atiyah et Hirzebruch sont de nature
topologique. Ce sont des classes de torsion dans $H^{4}(X(\C),\Z(2))$ sur une vari\'{e}t\'{e}
projective complexe $X(\C)$, qui ne peuvent pour des raisons de cobordisme complexe
(voir Totaro \cite{totaro}) \^{e}tre des classes de cycles analytiques
pour aucune structure  complexe  sur $X(\C)$.
 En particulier, dans le langage introduit
 au paragraphe \ref{subsecspecialisation},
  ces classes $\alpha\in H^{4}(\mathcal{X}_t(\C),\Z(2))$, avec $X\cong \mathcal{X}_t$, ont la propri\'{e}t\'{e}
que $\rho^2_{t,s}(\alpha)$ reste de Hodge en tout point $s$, et $\rho_{t,s}(\alpha)$ ne s'annule pas dans
$Z^{4}(\mathcal{X}_s)$, quel que soit $s\in T(\mathbb{C})$. L'espace de param\`{e}tres $T$ est ici
un ouvert non vide dans un espace projectif sur $\Q$. On peut donc choisir  $s\in T(\Q)$
et facilement donner des exemples  d'Atiyah--Hirzebruch d\'{e}finis sur $\Q$.
En cohomologie $l$-adique, on peut donner des exemples analogues
 sur des corps finis  (voir \cite{ctszamuely}).

Par ces m\'{e}thodes topologiques,
Totaro  (\cite[Thm. 7.1]{totaro})  a aussi construit des vari\'{e}t\'{e}s  projectives et lisses $X$
de dimension 7 sur $\C$, de type Godeaux--Serre,  telles que l'application $CH^2(X)/2 \to H^4_{\et}(X,\Z/2)$ n'est pas injective.
Ainsi   l'application $H^3(X(\C), \Z/2) \to H^3_{nr}(X,\Z/2)$ n'est pas surjective.

\medskip

\subsection{Exemples de Koll\'ar\label{subseckollar}}

Pour toute hypersurface lisse $X$ dans $\mathbb{P}^4_{\C}$, on a
  $H^4(X(\C),\Z(2))=\Z$ et donc $Hdg^4(X)=H^4(X(\C),\Z(2))=\Z$. Le sous-groupe $H^4_{alg}(X(\C),\Z(2)) \subset H^4(X(\C),\Z(2))=\Z$
  est l'id\'{e}al  $\Z.I \subset \Z$ (avec $I\in \N$) engendr\'{e} par les degr\'{e}s des courbes contenues dans $X$.
Dans \cite{kollar},
Koll\'ar montre que pour les hypersurfaces tr\`{e}s g\'{e}n\'{e}rales dans $\mathbb{P}^4$,
 de degr\'{e} suffisamment divisible, on a $I \neq 1$ et donc d'une part
 $Z^4(X)=\Z/I \neq 0$,  d'autre part $CH^2(X)/I \to H^4(X(\C),\mu_{I}^{\otimes 2})$ non injective
 (pour plus de d\'{e}tails, voir  \cite{soulevoisin}).

 Ces exemples sont tr\`{e}s instructifs du point de vue du probl\`{e}me de d\'{e}formation/sp\'{e}cialisation
 soulev\'{e} dans la section \ref{subsecspecialisation}. Tout d'abord,
 notons que la classe de Hodge $\alpha$ engendrant le groupe $H^4(X(\C),\Z(2))$ est de Hodge en tout point de la famille naturelle $T$
  de d\'{e}formations de $X$. Ainsi les fl\`{e}ches $\rho^2_{t,s}:Z^4(\mathcal{X}_t)\rightarrow Z^4(\mathcal{X}_s)$ de (\ref{flechedespe}), d\'{e}finies pour $t$ tr\`{e}s g\'{e}n\'{e}ral
 dans $T(\mathbb{C})$, sont surjectives pour tout $s$. Elles ne sont cependant pas injectives
 pour tout $s$: en effet, lorsque $\mathcal{X}_s$ contient une droite, la classe $\alpha$ devient alg\'{e}brique
 et donc $Z^4(\mathcal{X}_s)$=0.
 Il est m\^{e}me montr\'{e} dans \cite{soulevoisin}
 que, $t$ \'{e}tant fix\'{e} tr\`{e}s g\'{e}n\'{e}ral, la classe $\rho^2_{t,s}(\alpha)$ s'annule pour $s$ dans une union d\'{e}nombrable de ferm\'{e}s alg\'{e}briques de $T$, dense pour la topologie usuelle de $T(\mathbb{C})$.
 Ceci montre qu'il faut bien prendre une union d\'{e}nombrable d'ensembles alg\'{e}briques dans le lemme
 \ref{lemmedefalg}, alors que le th\'{e}or\`{e}me \ref{CDK} sugg\'{e}rerait qu'un ferm\'{e} alg\'{e}brique suffirait.

\begin{rema}{\rm Malgr\'{e} la n\'{e}cessit\'{e} d'\^{o}ter de $T(\mathbb{C})$ une union d\'{e}nombrable de sous-ensembles ferm\'{e}s alg\'{e}briques d\'{e}finis sur $\overline{\mathbb{Q}}$
pour construire des vari\'{e}t\'{e}s $X$ avec $Z^4(X)\not=0$, Hassett et Tschinkel ont montr\'{e} que les exemples de Koll\'ar peuvent \^{e}tre  choisis d\'{e}finis sur $\overline{\mathbb{Q}}$ (voir la remarque \ref{REHatschi}). Par contraste, sous l'hypoth\`{e}se que la conjecture de Tate sur les diviseurs sur les surfaces
sur un corps fini vaut, on  ne peut pas construire d'exemple \`{a} la Koll\'ar sur la cl\^{o}ture alg\'{e}brique d'un corps fini (\cite[\S 6]{ctszamuely}).

}
\end{rema}

\subsection{Exemples de Bloch--Esnault et Schoen}

Comme rappel\'{e} ci-dessus (Th\'{e}or\`{e}me \ref{blochesnault}),
 Bloch et Esnault  \cite{blochesnault}, resp. Schoen \cite{schoen},
 construisent des vari\'{e}t\'{e}s de dimension 3 pour lesquelles, pour
 un entier $n>1$ convenable, le quotient ${\rm Griff}^2(X)/n$ est non nul, resp. infini.
 Des suites exactes du d\'{e}but du paragraphe, il suit donc
que, dans les cas cit\'{e}s, $H^3_{nr}(X,\Z(2)) /n$ est non nul, resp. infini, et
qu'il en est de m\^{e}me de $H^3_{nr}(X,\mu_{n}^{\otimes 2})$.
Pour les exemples de Schoen, il est alors clair que le noyau de
$CH^2(X)/n \to H^4_{\et}(X,\mu_{n}^{\otimes 2})$ est infini.

  Dans les exemples de Bloch et Esnault, le groupe  $H^4(X(\mathbb{C}),\Z)$ est sans
torsion. Des suites exactes du d\'{e}but du paragraphe et
de la non nullit\'{e} du groupe ${\rm Griff}^2(X)/n$  on d\'{e}duit   alors que
l'application $CH^2(X)/n \to H^4_{\et}(X,\mu_{n}^{\otimes 2})$ n'est pas injective.
Ceci laisse ouverte la question de la trivialit\'{e} du groupe
  $Z^4(X)_{\rm tors}$ pour ces vari\'{e}t\'{e}s.

 \subsection{Exemples de Gabber}

Soient   $E_{i}, i=1,2,3$ des courbes elliptiques,
  $l$ un nombre premier,
et pour chaque $i$,  $\alpha_{i} \in H^1_{\et}(E_{i},\Z/l)$ non nul.
Un cas particulier d'un r\'{e}sultat de Gabber  \cite{gabber} dit que si les
$j$-invariants   des courbes $E_{i} $ sont alg\'{e}briquement ind\'{e}pendants sur $\Q$,
alors le cup-produit $\alpha_{1}\cup \alpha_{2} \cup \alpha_{3} \in H^3_{\et}(X,\Z/l)$
a une image non nulle dans $H^3(\C(X),\Z/l)$. On obtient ainsi des
classes non nulles dans $H^3_{nr}(X,\Z/l)$.  Mais celles-ci proviennent
de classes dans $H^3_{nr}(X,\Z)/l$, elles ont  donc une
 image nulle
dans $Z^4(X)[l]$.

\subsection{Certaines vari\'{e}t\'{e}s unirationnelles}

Les vari\'{e}t\'{e}s consid\'{e}r\'{e}es dans les exemples pr\'{e}c\'{e}dents ne sont pas
rationnellement connexes (les exemples d'Atiyah et Hirzebruch, par construction,
admettent des rev\^{e}tements non ramifi\'{e}s).

On a
 le r\'{e}sultat suivant, qui repose sur un th\'{e}or\`{e}me d'Arason sur la cohomologie galoisienne des corps de fonctions
 de certaines quadriques :
 \begin{theo}\label{theoctojan} (Colliot-Th\'{e}l\`{e}ne--Ojanguren \cite{colliotojanguren}) Il existe des  vari\'{e}t\'{e}s  $X$ projectives, lisses, con\-nexes,
 unirationnelles, de dimension $6$, telles que
 $H^3_{nr}(X,\Z/2)\not=0$.
 \end{theo}
 \'{E}tant unirationnelles, ces vari\'{e}t\'{e}s ont un groupe  $CH_0$ isomorphe \`{a} $\mathbb{Z}$.
 Le th\'{e}or\`{e}me  \ref{blochsrinivasgriffnul} du pr\'{e}sent article montre alors que l'on a $H^3_{nr}(X,\Z(2))=0$ et
 $Z^4(X)=Z^4(X)_{\rm tors} \neq 0$.
Ces exemples permettent donc, en  dimension $\geq 6$, de {\it r\'{e}pondre
  n\'{e}gativement  \`{a} une partie de la question 16
   pos\'{e}e  dans} \cite{voisinjapjmath}: la conjecture de Hodge enti\`{e}re est-elle
   satisfaite pour les classes de Hodge de degr\'{e} $4$ sur les vari\'{e}t\'{e}s rationnellement connexes?

 Les exemples de \cite{colliotojanguren} sont des vari\'{e}t\'{e}s  fibr\'{e}es en quadriques de dimension 3 au-dessus
 de l'espace projectif de dimension 3.
 D'autres exemples de vari\'{e}t\'{e}s unirationnelles satisfaisant $H^3_{nr}(X,\mu_{n}^{\otimes 2}) \neq 0$
pour $n$ convenable ont \'{e}t\'{e} obtenus par E. Peyre \cite{peyre1, peyre2}. Pour la m\^{e}me
raison que pr\'{e}c\'{e}demment,
   ces exemples satisfont  aussi $Z^4(X)=Z^4(X)_{\rm tors} \neq 0$.

La construction donn\'{e}e dans
\cite{colliotojanguren} de  classes non triviales
dans $H^3_{nr}(X,\Z/2)$ est explicite. Des mod\`{e}les lisses
des vari\'{e}t\'{e}s consid\'{e}r\'{e}es \'{e}tant difficiles \`{a} construire, les \'{e}l\'{e}ments
non triviaux correspondants du groupe $Z^4(X)$, c'est-\`{a}-dire des classes de Hodge enti\`{e}res
non alg\'{e}briques, sont difficiles \`{a} analyser.
Ces classes proviennent-elles  ou non comme dans les exemples d'Atiyah-Hirzebruch
de  classes de torsion dans $H^4(X(\C),\Z(2))$?
Ce probl\`{e}me est li\'{e} \`{a} la question de savoir si l'application $CH^2(X)/2 \to H^4_{\et}(X,\Z/2)$
est injective.

 La th\'{e}orie des d\'{e}for\-mations de ces vari\'{e}t\'{e}s
semble
tr\`{e}s difficile. Ces classes de Hodge
 sont-elles  stables par  d\'{e}formation, c'est-\`{a}-dire, sont-elles dans l'image de l'application
 de sp\'{e}cialisation $\rho^2_{t,s}$ de (\ref{flechedespe}) pour $t$ tr\`{e}s g\'{e}n\'{e}ral dans
 la base d'une famille compl\`{e}te de d\'{e}formations de $X$?
Il n'est pas clair non plus si, se pla\c{c}ant en un point tr\`{e}s g\'{e}n\'{e}ral $t$ de
la famille de d\'{e}formations $T$ de $X$ pr\'{e}servant ces classes de Hodge (cf. th\'{e}or\`{e}me
\ref{CDK}), ces classes peuvent s'annuler sous l'application de sp\'{e}cialisation
$\rho^2_{t,s}$, pour un $s\in T(\mathbb{C})$.

\subsection{Une vari\'{e}t\'{e} $X$ de dimension 3 avec  $H^i(X,\mathcal{O}_X)=0$ pour $i>0$}\label{nouveaunonnul}

On donne maintenant un exemple de telle vari\'{e}t\'{e} (projective et lisse)
avec $Z^4(X)\neq 0$.  Dans cet exemple, le groupe $Z^4(X)$
n'est pas localement constant par d\'{e}formations.

Soit $G=\Z/5$. On choisit une racine $5$-i\`{e}me de l'unit\'{e} $\zeta$ non triviale et un g\'{e}n\'{e}rateur $g$ de $G$, et on fait agir $G$ sur $\mathbb{P}^1=Proj \,\C[x,y]$ et sur
$\mathbb{P}^3=Proj \,\C[x_0,x_1,x_2,x_3]$ de la fa\c{c}on  suivante:

$$g^*x=x,\,g^*y=\zeta y,$$
$$g^*x_i=\zeta^ix_i,\,i=0,\ldots, 3.$$

Soit $X\subset \mathbb{P}^1\times \mathbb{P}^3$ une hypersurface de bidegr\'{e} $(3,4)$  d\'{e}finie par une \'{e}quation $f=0$, avec $f\in H^0(\mathbb{P}^1\times \mathbb{P}^3,\mathcal{O}_{\mathbb{P}^1\times \mathbb{P}^3}(3,4))$ invariante sous $G$.
Soit $\tau:\widetilde{W}\rightarrow W$ une d\'{e}singularisation  de $W=X/G$.
Si $f$ est g\'{e}n\'{e}riquement choisie, $X$ n'a que des singularit\'{e}s isol\'{e}es et  un th\'{e}or\`{e}me de Bertini implique que la premi\`{e}re projection
$pr_1:X\rightarrow \mathbb{P}^1$ a pour fibre g\'{e}n\'{e}rale une  surface $K3$  lisse de degr\'{e} $4$.
Le groupe $G$ agit librement sur $X$ et $\mathbb{P}^1$ en dehors  de points
fixes isol\'{e}s et $pr_1$ est  \'{e}quivariante. Il en r\'{e}sulte que l'application $W\rightarrow \mathbb{P}^1/G$ a aussi pour fibre g\'{e}n\'{e}rale une  surface $K3$  lisse de degr\'{e} $4$.
\begin{prop}\label{familleK3}
Soient $X,W, \widetilde{W}$ comme ci-dessus. Si $X$ est  tr\`{e}s  g\'{e}n\'{e}rale en modules, alors :

(i) On a
$H^i(\widetilde{W},\mathcal{O}_{\widetilde{W}})=0$ pour $i>0$.

(ii) La $2$-torsion du groupe $Z^4(\widetilde{W})$ est non nulle.

(iii) Le groupe $H^3_{nr}(\widetilde{W},\Z/2)$ est  non nul.
\end{prop}
{\bf D\'{e}monstration.} Notons d'abord que pour un choix g\'{e}n\'{e}rique de $f$,  la vari\'{e}t\'{e} $X$ a au pire des singularit\'{e}s
quadratiques ordinaires.
En effet, par le th\'{e}or\`{e}me de Bertini $X$ est g\'{e}n\'{e}riquement
lisse en dehors du lieu de base du syst\`{e}me lin\'{e}aire
$| \mathcal{O}_{\mathbb{P}^1\times \mathbb{P}^3}(3,4))^G|$, qui est constitu\'{e}
des points
$$O_1:y=0,x_0=0,x_2=0,x_1=0,$$
$$ O_2:y=0,x_0=0,x_2=0,x_3=0,$$
$$O_3:y=0,x_0=0,x_3=0,x_1=0,
$$
$$O'_1: x=0,x_3=0, x_0=0, x_1=0,$$
$$  O'_2: x=0,X_3=0, x_0=0, x_2=0,   $$
$$O'_3:x=0,x_3=0, x_1=0, x_2=0 .$$
 Un examen de l'\'{e}quation
g\'{e}n\'{e}rique de $X$ en ces points permet de conclure.

Comme les singularit\'{e}s de $X$ sont rationnelles,
pour toute d\'{e}singularisation $\widetilde{X}$ de $X$, on a
$H^i(\widetilde{X},\mathcal{O}_{\widetilde{X}})=H^i(X,\mathcal{O}_X)$.
Pour toute d\'{e}singularisation $\widetilde{W}$ de $W$, l'application rationnelle naturelle
$$Q:\widetilde{X}\dashrightarrow \widetilde{W}$$
est $G$-\'{e}quivariante, et identifie birationnellement $\widetilde{W}$ \`{a} $\widetilde{X}/G$.
 Il en r\'{e}sulte
que
\begin{eqnarray}\label{eqn3avril}H^i(\widetilde{W},\mathcal{O}_{\widetilde{W}})\subset H^i(\widetilde{X},\mathcal{O}_{\widetilde{X}})^G=
H^i({X},\mathcal{O}_{{X}})^G.
\end{eqnarray}
En effet, nous pouvons supposer que l'application quotient $Q$ est un morphisme. Il en r\'{e}sulte
 (cf. \cite[7.3.2]{voisinlivre}) que l'application $$Q^*: H^i(\widetilde{W},\mathcal{O}_{\widetilde{W}})\rightarrow H^i(\widetilde{X},\mathcal{O}_{\widetilde{X}})
 =H^i({X},\mathcal{O}_{{X}})$$ est injective, et  elle envoie
$H^i(\widetilde{W},\mathcal{O}_{\widetilde{W}})$ dans
$H^i(\widetilde{X},\mathcal{O}_{\widetilde{X}})^G$ par $G$-invariance de $Q$.

Pour \'{e}tablir l'\'{e}nonc\'{e} (i), il suffit donc de montrer que l'on a $H^i({X},\mathcal{O}_{{X}})^G=0$ pour $i>0$. Pour $i=1,\,2$, cela r\'{e}sulte du fait
que $X\subset \mathbb{P}^1\times \mathbb{P}^3$ est un diviseur ample. Pour $i=3$, on doit montrer de fa\c{c}on \'{e}quivalente que
$H^0(X,K_X)^G=0$. Or par adjonction $K_X=\mathcal{O}_X(1,0)$, et plus pr\'{e}cis\'{e}ment
les sections de $K_X$ sont donn\'{e}es par les
 $\eta_P:=Res_X \frac{P\Omega}{f}$, o\`{u} $P\in H^0(X,\mathcal{O}_X(1,0))$ et $\Omega$ est le g\'{e}n\'{e}rateur
 naturel de $H^0(\mathbb{P}^1\times \mathbb{P}^3,K_{\mathbb{P}^1\times \mathbb{P}^3}(2,4))$.

 Notant que $g^*\Omega=\zeta^2\Omega$, on conclut que $\eta_P$ est $G$-invariante
 si et seulement si $g^*P=\zeta^3P$. Comme $g$ n'agit pas avec la valeur propre $\zeta^3$
 sur $H^0(\mathbb{P}^1,\mathcal{O}_{\mathbb{P}^1}(1))$, on a bien montr\'{e} que $H^0(X,K_X)^G=0$.

 Notons $\pi: {W}\rightarrow \mathbb{P}^1/G\cong\mathbb{P}^1$ l'application naturelle, et
 $l:=(\pi\circ\tau)^* c_1(\mathcal{O}_{\mathbb{P}^1}(1))$ : c'est la classe des fibres de $\widetilde{W} \to \mathbb{P}^1$.
 Comme $H^2(\widetilde{W},\mathcal{O}_{\widetilde{W}})=0$, la structure de Hodge sur $H^4(\widetilde{W} (\mathbb{C}),\mathbb{Q})$ (qui est dual au sens de Poincar\'{e} de $H^2(\widetilde{W}(\mathbb{C}),\mathbb{Q})$)
 est triviale, c'est-\`{a}-dire enti\`{e}rement de type $(2,2)$.

 Toutes les classes enti\`{e}res dans
 $H^4(\widetilde{W}(\mathbb{C}),\Z)$ sont donc de Hodge.
Pour montrer que la $2$-torsion
 de $Z^4(\widetilde{W})$ est non nulle, il suffit de montrer:

(1)  Il existe une classe
 $\alpha\in H^4(\widetilde{W}(\mathbb{C}),\Z)=H_2(\widetilde{W}(\mathbb{C}),\Z)$ telle que
 $\deg_\alpha l=5$.

 (2) Pour $X$ tr\`{e}s g\'{e}n\'{e}rale, et pour toute classe alg\'{e}brique $[Z]\in H^4(\widetilde{W}
 (\mathbb{C}),\Z)$, $\deg_{[Z]}l$ est pair. De fa\c{c}on \'{e}quivalente, pour toute courbe
    $Z$  contenue dans $W$, le degr\'{e} de la restriction de  $\pi$ \`{a} $Z$ est pair.

     En effet, l'intersection avec $Z$ induit alors un homomorphisme de groupes finis
      $Z^4(\widetilde{W}) \to \Z/2$
     qui envoie $\alpha$ sur la classe $1 \in \Z/2$.

L'\'{e}nonc\'{e} (1) est obtenu en sp\'{e}cialisant $X$  en un $X_0$ g\'{e}n\'{e}rique contenant une droite $\mathbb{P}^1\times M$, o\`{u} $M$ est
 g\'{e}n\'{e}riquement choisi. On v\'{e}rifie que $X_0$ est alors lisse le long de cette droite, et que
 $W_0$ est lisse le long de son image $\Delta$ dans $W_0$.
 La courbe $\Delta $ se rel\`{e}ve naturellement en une courbe  $\tilde{\Delta}\subset\widetilde{W_0}$ qui  satisfait $\deg_{\tilde{\Delta}} l=5$. Comme $W_0$ est lisse
 le long de $\Delta$, la classe $\delta=[\tilde{\Delta}]\in H_2(\widetilde{W_0}(\mathbb{C}),\Z)$
 s'\'{e}tend par ailleurs en une classe dans $H_2(\widetilde{W}(\mathbb{C}),\Z)$
   pour toute petite d\'{e}formation
   $W$ de $W_0$, ce qui donne la classe $\alpha$ souhait\'{e}e.

 La preuve de (2)  se fait par sp\'{e}cialisation. En effet, $X$ et donc $W$ \'{e}tant choisies tr\`{e}s g\'{e}n\'{e}rales,
 tout $1$-cycle de $W$ se sp\'{e}cialise dans n'importe quelle sp\'{e}cialisation $W_0$ de $W$.
 Il suffit donc de sp\'{e}cialiser $W$ en $W_0$, de fa\c{c}on que
 $W_0$ ne contienne pas de $1$-cycle de degr\'{e} impair au-dessus de $\mathbb{P}^1$.
 On choisit la sp\'{e}cialisation suivante (cf. \cite{starr}):
 On consid\`{e}re le morphisme $G$-\'{e}quivariant
 $$\phi: \mathbb{P}^1\rightarrow \mathbb{P}^1$$
 d\'{e}fini par $\phi^*x=u^4,\,\phi^*y=v^4$, o\`{u} $u,v$ sont des coordonn\'{e}es homog\`{e}nes sur $\mathbb{P}^1$, avec l'action (lin\'{e}aris\'{e}e) suivante de $G$: $g^*u=u,\,g^*v=\zeta^4v$.
 On choisit un \'{e}l\'{e}ment  $G$-invariant
 g\'{e}n\'{e}rique $Q$ (relativement \`{a} cette derni\`{e}re action)
 de $H^0(\mathbb{P}^1\times \mathbb{P}^3,\mathcal{O}_{\mathbb{P}^1\times \mathbb{P}^3}(3,1))$.

 Un tel \'{e}l\'{e}ment s'\'{e}crit (modulo l'action de scalaires sur les $x_i$) sous la forme
 \begin{eqnarray}\label{equation}
 Q=u^3x_0+u^2vx_1+uv^2x_2+v^3x_3.
 \end{eqnarray}
  Soit $\Gamma\subset \mathbb{P}^1\times \mathbb{P}^3$ le diviseur de
$Q$, et soit 
$X_0:=(\phi,Id)(\Gamma)$.
Comme le degr\'{e} de $\phi$ est $4$,  on a $X_0\in \mid \mathcal{O}_{\mathbb{P}^1\times \mathbb{P}^3}(3,4)\mid $.
Par $G$-\'{e}quivariance (lin\'{e}aris\'{e}e) de $\phi$, $X_0$ est de plus d\'{e}fini par une \'{e}quation $G$-invariante.
Le lemme \ref{ledegenerescence} suivant conclut la d\'{e}monstration de l'\'{e}nonc\'{e} (ii).
Le  th\'{e}or\`{e}me \ref{principal1} donne alors l'\'{e}nonc\'{e} (iii).
\cqfd
\begin{lemm}\label{ledegenerescence} Pour toute courbe $Z\subset W_0$, le degr\'{e} de
$\pi_{\mid Z}: Z\rightarrow\mathbb{P}^1$ est pair.
\end{lemm}
{\bf D\'{e}monstration.} Comme le degr\'{e} de l'application quotient
$X_0\rightarrow W_0$ est impair, il suffit
de montrer que pour toute courbe $Z\subset X_0$,
le degr\'{e} de $pr_{1\mid Z}: Z\rightarrow\mathbb{P}^1$ est pair. La preuve est semblable \`{a} celle de Starr \cite{starr}.
La fibre de $pr_1:\Gamma\rightarrow \mathbb{P}^1$ au-dessus de $t$ est un hyperplan $H_t$, et donc la fibre de $pr_1:X_0\rightarrow\mathbb{P}^1$ est constitu\'{e}e de l'union
des quatre hyperplans $H_{t_i}$, o\`{u} $\{t_1,\ldots,t_4\}=\phi^{-1}(\{t\})$. Ces quatre hyperplans sont
 pour $t$ g\'{e}n\'{e}rique en position g\'{e}n\'{e}rale, comme le montre l'\'{e}quation (\ref{equation}). Soit $Z\subset X_0$
une courbe dominant $\mathbb{P}^1$. On peut supposer $Z$ irr\'{e}ductible. Au point g\'{e}n\'{e}rique $z$ de $Z$, envoy\'{e} par $pr_1$
sur le  point $t\in \mathbb{P}^1$, $z$ appartient \`{a} $m$ des hyperplans $H_{t_i}$, o\`{u}
$m\leq 3$. Donc il y a une factorisation naturelle de
$pr_{1\mid Z}$ \`{a} travers la courbe
$C_m\rightarrow \mathbb{P}^1$ param\'{e}trant $m$ points dans les fibres de $\phi:\mathbb{P}^1\rightarrow \mathbb{P}^1$. Pour $m=1$ ou $3$, $C_m\rightarrow \mathbb{P}^1$ s'identifie \`{a}
$\phi:\mathbb{P}^1\rightarrow \mathbb{P}^1$. Donc le degr\'{e} de $C_m$ sur $\mathbb{P}^1$ est divisible par $4$.
Pour $m=2$, comme $\phi$ est un rev\^{e}tement galoisien cyclique d'ordre $4$, $C_2\rightarrow \mathbb{P}^1$
a deux composantes au-dessus de $\mathbb{P}^1$, l'une de degr\'{e} $2$, l'autre de degr\'{e} $4$. Donc toutes
les composantes des courbes $C_m,\,m\leq 3$, sont de degr\'{e} pair au-dessus de $\mathbb{P}^1$.
Comme $pr_{1\mid Z}$ se factorise \`{a} travers l'une de ces composantes, on conclut que le
degr\'{e} de $pr_{1\mid Z}$ est pair.
\cqfd

\medskip

\begin{rema}\label{REHatschi}{\rm De tels exemples peuvent en fait \^{e}tre construits sur $\mathbb{Q}$ et c'est aussi le cas
pour les exemples de Koll\'ar (cf. section \ref{subseckollar}).
 Ceci nous a  \'{e}t\'{e} signal\'{e} par Hassett et Tschinkel.
En effet, la terminologie
\og tr\`{e}s g\'{e}n\'{e}rale \fg \, a \'{e}t\'{e} utilis\'{e}e dans la d\'{e}monstration ci-dessus pour s'assurer que
tout $1$-cycle de $W$ se sp\'{e}cialise avec $W$, ce qu'on obtient en demandant
que $W$ ne soit pas param\'{e}tr\'{e} par un point de l'espace de modules
qui appartient \`{a} l'image d'un sch\'{e}ma de Hilbert relatif   ne dominant pas l'espace de modules.

Mais si $W$ est d\'{e}finie sur
$\mathbb{Q}$,  les cycles de $W$ se sp\'{e}cialisent aussi bien aux r\'{e}ductions
de $W$ modulo $p$ (pour presque tout $p$). Dans l'argument pr\'{e}c\'{e}dent, il suffit donc de supposer que $W$ a une r\'{e}duction
modulo $p$, pour $p$ premier ad\'{e}quat, qui est de la forme $W_0$ comme ci-dessus. En effet, le lemme \ref{ledegenerescence}
est vrai en caract\'{e}ristique $p>2$.
}
\end{rema}
\begin{rema}{\rm
  Comme on le voit
  dans la preuve, les contre-exemples  \`{a} la conjecture de Hodge enti\`{e}re  construits dans   la proposition \ref{familleK3}  sont,
  comme ceux de Koll\'ar,
instables par d\'{e}formation :
la classe consid\'{e}r\'{e}e non nulle $\alpha\in Z^4(\widetilde{W}_t)$, o\`u $t$ est un point tr\`{e}s g\'{e}n\'{e}ral
de la base $T$ param\'{e}trant les d\'{e}formations des vari\'{e}t\'{e}s  $\widetilde{W}$,
 s'annule en certains points : $\rho_{t,s}(\alpha)=0$ dans $Z^4(\widetilde{\mathcal{W}}_s)$ pour
 certains $s\in T(\mathbb{C})$.
}
\end{rema}

\section{Cas de nullit\'{e} des groupes $Z^4(X)$ et $H^3_{nr}(X,\mu_{n}^{\otimes 2})$}

Le th\'{e}or\`{e}me \ref{theoctojan} montre que l'invariant $H^3_{nr}(X,\Z/2)$, ou encore la $2$-torsion
du groupe $Z^4(X)$, peuvent \^{e}tre non nuls pour des vari\'{e}t\'{e}s unirationnelles, donc
 rationnellement connexes,  de dimension $6$.
Cependant, on a le :
 \begin{theo} (Voisin \cite{voisinuniruled}) \label{theovoisinuniruled}
 Soit $X$ une vari\'{e}t\'{e} projective et lisse de dimension 3 sur $\C$.
 Le groupe   $Z^4(X)$ est nul  dans chacun des cas suivants :

 (i) $X$ est   unir\'{e}gl\'{e}e;

 (ii) $X$ est une vari\'{e}t\'{e} de Calabi-Yau.

 \end{theo}

 \begin{coro}\label{corodevoisinuniregle}
 Soit $k$ un corps alg\'{e}briquement clos de caract\'{e}ristique z\'{e}ro et $X$
un solide projectif et lisse sur $k$. Si  la vari\'{e}t\'{e} $X$ est unir\'{e}gl\'{e}e, alors
pour tout entier $n>0$, les groupes $H^3_{nr}(X,\mu_{n}^{\otimes 2})$
 sont nuls.
  \end{coro}
   {\bf D\'{e}monstration.}  Pour obtenir le
   r\'{e}sultat lorsque $k=\C$,
il suffit de combiner les th\'{e}or\`{e}mes
  \ref{principal2},
   \ref{blochsrinivasgriffnul} et   \ref{theovoisinuniruled}.
  Il est par ailleurs connu que les groupes de cohomologie non   ramifi\'{e}e
  \`{a} coefficients de torsion
  sont invariants par extension de corps de base alg\'{e}briquement clos
  (\cite[Thm. 4.4.1]{ctsantabarbara}).
  \cqfd

  \begin{prop}\label{cyclemodnratcon}
  Soit $X$
un solide projectif et  lisse sur $\C$, unir\'{e}gl\'{e} et sans torsion dans son groupe de Picard, et tel que
$H^2(X,\mathcal{O}_X)=0$.
Pour tout entier $n>0$ :

(i) L'application classe de cycle
$CH^2(X)/n \to H^4_{\et}(X,\mu_{n}^{\otimes 2})$
est un isomorphisme.

(ii) Les groupes $H^{i}(X,\H^3_{X}(\mu_{n}^{\otimes 3}))$ sont nuls
pour $i\neq 3$, et $H^3(X,\H^3_{X}(\mu_{n}^{\otimes 3}))\simeq \Z/n$.

  \end{prop}
  Notons que les hypoth\`{e}ses impliquent
  $H^{i}(X,\O_{X})=0$  pour $i>0$.

   {\bf D\'{e}monstration.}
D'apr\`{e}s le corollaire  \ref{corodevoisinuniregle}, on a
   $H^{0}(X,\H^3_{X}(\mu_{n}^{\otimes 2}))=0$, et  on a donc  aussi  $H^{0}(X,\H^3_{X}(\mu_{n}^{\otimes i}))=0$ pour tout $i$.

 Pour  une vari\'{e}t\'{e} $X$ de dimension 3 satisfaisant
 $H^2(X,\mathcal{O}_X)=0$, la structure de Hodge sur $H^4(X(\mathbb{C}),\mathbb{Q}(2))$ est triviale. On a donc
 $H^4(X(\mathbb{C}),\mathbb{Z}(2))=  Hdg^4(X,\mathbb{Z})$. Si de plus
 $X$ est unir\'{e}gl\'{e}e, le th\'{e}or\`{e}me \ref{theovoisinuniruled} dit alors que l'application classe de
 cycle $CH^2(X)\rightarrow H^4(X(\mathbb{C}),\mathbb{Z}(2))$ est surjective.

On a la suite exacte
$$ H^4(X(\C),\Z(2))/n \to H^4(X(\C),\mu_{n}^{\otimes 2}) \to H^5(X(\C),\Z(2))[n].$$
Par dualit\'{e} de Poincar\'{e}, la $n$-torsion de $H^5(X(\C),\Z)$ s'identifie \`{a} la torsion
de $H_{1}(X(\C),\Z)$. Par hypoth\`{e}se, $X(\C)$ ne poss\`{e}de pas de rev\^{e}tements finis
ab\'{e}liens non ramifi\'{e}s. On a donc $H^5(X(\C),\Z(2))[n]=0$, d'on l'on d\'{e}duit que
l'application $CH^2(X)/n \to H^4(X(\C), \mu_{n}^{\otimes 2}    )$ est  surjective, ce qui entra\^{\i}ne
le m\^{e}me \'{e}nonc\'{e} en cohomologie \'{e}tale. Enfin l'injectivit\'{e} de cette fl\`{e}che r\'{e}sulte
de l'annulation
$H^{0}(X,\H^3_{X}(\mu_{n}^{\otimes 2}))=0$. Ceci \'{e}tablit l'\'{e}nonc\'{e} (i).

   On consid\`{e}re, pour tout $i$ entier, la suite spectrale de Bloch--Ogus pour la cohomologie \'{e}tale :
$$E_{2}^{pq} = H^p(X,\H^q_{X}(\mu_{n}^{\otimes i})) \Longrightarrow H^n_{\et}(X,\mu_{n}^{\otimes i}).$$
Le th\'{e}or\`{e}me \ref{boresolutionfinie}  et le th\'{e}or\`{e}me de
Lefschetz faible pour la cohomologie \'{e}tale montrent que
cette suite spectrale est concentr\'{e}e dans le triangle $0 \leq p \leq q \leq 3$.
On en d\'{e}duit alors
 l'isomorphisme  $H^{2}(X,\H^3_{X}(\mu_{n}^{\otimes 2}))\oi H^5_{\et}(X,\mu_{n}^{\otimes 2})$.
Par la dualit\'{e} de Poincar\'{e} en cohomologie \'{e}tale, ce dernier groupe
est dual de $H^1_{\et}(X,\mu_{n})$. Mais ce groupe  n'est autre que
la $n$-torsion du groupe de Picard de $X$, par hypoth\`{e}se nulle.
De la suite spectrale on d\'{e}duit
 l'isomorphisme  $H^{3}(X,\H^3_{X}(\mu_{n}^{\otimes 3}))\oi H^6_{\et}(X,\mu_{n}^{\otimes 3}) \simeq \Z/n$.
En utilisant la nullit\'{e} de $H^{0}(X,\H^3_{X}(\mu_{n}^{\otimes 2}))$,
on d\'{e}duit de la suite spectrale la suite
exacte
$$0 \to  CH^2(X)/n \to H^4_{\et}(X,\mu_{n}^{\otimes 2}) \to H^1(X,\H^3_{X}(\mu_{n}^{\otimes 2})) \to 0,$$
o\`{u} la premi\`{e}re application est l'application classe de cycle.
L'\'{e}nonc\'{e}  (i) montre alors que l'on a  $H^1(X,\H^3_{X}(\mu_{n}^{\otimes 2}))=0$
et donc  $H^{1}(X,\H^3_{X}(\mu_{n}^{\otimes i}))=0$ pour tout $i$.
\cqfd

 \begin{rema} {\rm
Des propri\'{e}t\'{e}s de rigidit\'{e} des groupes $H^{i}(X,\H^{j}(\Z/n))$ en principe connues
(\cite[Rmk. 4.4.2]{ctsantabarbara};  U. Jannsen, \og Rigidity results on K-theory and other functors\ \fg \
(non publi\'{e}, 1995)) permettent de d\'{e}duire du th\'{e}or\`{e}me ci-dessus le m\^{e}me \'{e}nonc\'{e}
sur un corps alg\'{e}briquement clos de caract\'{e}ristique z\'{e}ro.}
  \end{rema}

  \begin{rema}
{\rm  La proposition ci-dessus s'applique en particulier \`{a} tout solide rationnellement connexe.
 On sait en effet que toute vari\'{e}t\'{e} projective lisse rationnel\-lement connexe est
 alg\'{e}bri\-quement simplement connexe et satisfait de plus $H^2(X,\mathcal{O}_X)=0$.}
 \end{rema}

 \medskip

On pourrait se poser la question :
Pour une vari\'{e}t\'{e} $X$ de dimension 3 dont le groupe de Chow  des z\'{e}ro-cycles est support\'{e} sur
  une surface, le groupe fini  $H^3_{nr}(X, \Q/\Z(2))  \simeq Z^4(X)$ (voir le th\'{e}or\`{e}me
  \ref{blochsrinivasgriffnul})
est-il nul ?
La r\'{e}ponse est tr\`{e}s probablement non,
  car si la conjecture de Bloch
 (voir \cite[lecture 1]{blochbook} pour les surfaces et
 \cite[Conjecture 23.23]{voisinlivre} en dimension arbitraire)
 est satisfaite par les vari\'{e}t\'{e}s $X$ construites
 au paragraphe \ref{nouveaunonnul},
 celles-ci fournissent alors un contre-exemple. De fait  on devrait avoir pour ces vari\'{e}t\'{e}s
 un isomorphisme $\deg : CH_{0}(X)\oi \Z$.

  \medskip

 Les r\'{e}sultats  et exemples ci-dessus laissent n\'{e}anmoins  ouvertes les  questions suivantes.

  \medskip

\begin{question}\label{question6.2}
 Pour une vari\'{e}t\'{e} $X$ rationnellement
 connexe de dimension 4 ou 5, le groupe fini \break  $H^3_{nr}(X, \Q/\Z(2))  \simeq Z^4(X)$
 est-il nul ?
 \end{question}

  \medskip

 \begin{question}\label{question6.3}
 (\cite[Question 16]{voisinjapjmath}) Pour une vari\'{e}t\'{e} $X$ rationnellement
 connexe de dimension $d$ quelconque, le groupe fini  $Z^{2d-2}(X)$ est-il nul ?
 \end{question}

 D'apr\`{e}s \cite{voisinuniruled}, la r\'{e}ponse  \`{a} cette derni\`{e}re question
 est affirmative en dimension
$d\leq 3$. En dimension $4$, mentionnons    le  r\'{e}sultat
(cf. \cite[Theorem 1.7]{hoeringvoisin}) : Pour une vari\'{e}t\'{e} de Fano $X$
lisse de dimension $4$, le groupe $Z^6(X)$ est nul.

 \medskip

Pour la question \ref{question6.2}
une r\'{e}ponse affirmative pour les hypersurfaces cubiques de dimension
$4$ est donn\'{e}e dans \cite{voisinjapjmath}.
On va \'{e}tablir ici un r\'{e}sultat plus g\'{e}n\'{e}ral. On s'int\'{e}resse
 ici au  groupe $Z^4(X)$, o\`{u} $f:X\rightarrow \Gamma$ est une fibration en vari\'{e}t\'{e}s projectives
 $X_t$ de dimension $3$ satisfaisant
 $H^3(X_t,\mathcal{O}_{X_t})=H^2(X_t,\mathcal{O}_{X_t})=0$.
  Pour $X_t,\,t\in \Gamma(\mathbb{C})$, comme ci-dessus, la jacobienne interm\'{e}diaire
  $J^2(X_t)$ est une vari\'{e}t\'{e} ab\'{e}lienne, du fait que $H^3(X_t,\mathcal{O}_{X_t})=0$ (cf. \cite[12.2.2]{voisinlivre})
  et pour toute classe
 $\alpha\in H^4(X_t(\mathbb{C}),\mathbb{Z})=Hdg^4(X_t,\mathbb{Z})$, on dispose d'un torseur
 $J^2(X_t)_\alpha$ sous $J^2(X_t)$ qui est une vari\'{e}t\'{e} alg\'{e}brique et qu'on peut d\'{e}finir
 par la formule:
 $$J^2(X_t)_\alpha:=d^{-1}(\alpha),$$ o\`{u} $d: H^4_D(X_t(\mathbb{C}),\mathbb{Z}(2))\rightarrow   Hdg^4(X_t,\mathbb{Z})$
 est l'application naturelle de la cohomologie de Deligne vers la cohomologie de Betti \`{a} coefficients entiers, d'image
  $Hdg^4(X,\mathbb{Z})$ (cf. \cite[12.3.1]{voisinlivre}).
 Pour toute famille
 de $1$-cycles $\mathcal{Z}\subset B_t\times X_t$ sur $X_t$ de classe $[Z_t]=\alpha$
 param\'{e}tr\'{e}e par une vari\'{e}t\'{e}
 alg\'{e}brique $B_t$, l'application d'Abel-Jacobi (ou plut\^{o}t classe de Deligne) de $X_t$ induit un morphisme
 de vari\'{e}t\'{e}s alg\'{e}briques complexes:
 $$AJ_{X_t}^2:B_t\rightarrow J^2(X_t)_\alpha.$$
 Cette situation se met en famille au-dessus de l'ouvert $\Gamma_0\subset \Gamma$
 de lissit\'{e} de $f$. On dispose donc d'une famille de vari\'{e}t\'{e}s
 ab\'{e}liennes $\mathcal{J}^2\rightarrow \Gamma_0$, et pour toute section   $\alpha $
 de $R^4f_*\mathbb{Z}$ sur $\Gamma_{0}$, de la famille tordue $\mathcal{J}^2_\alpha\rightarrow \Gamma_0$.
 Etant donn\'{e}s une vari\'{e}t\'{e} alg\'{e}brique $B$, un morphisme $g:B\rightarrow \Gamma$ et
 une famille de $1$-cycles relatifs $\mathcal{Z}\subset B\times_\Gamma X$ de classe
 $[Z_t]=\alpha_t\in H^4(X_t,\mathbb{Z})$, l'application d'Abel-Jacobi relative fournit un morphisme
 $\phi_\mathcal{Z}:B_0\rightarrow \mathcal{J}^2_\alpha$
 au-dessus de $\Gamma_0$, o\`{u} $B_0:=g^{-1}(\Gamma_0)$.
\begin{theo} \label{theocritereabeljacobi} Supposons que  $f:X\rightarrow \Gamma$ satisfait  les hypoth\`{e}ses suivantes:

(i)  Les fibres lisses $X_t$ de $f$  sont
  de dimension $3$ et satisfont
 $H^3(X_t,\mathcal{O}_{X_t})=H^2(X_t,\mathcal{O}_{X_t})=0$.

(ii) Pour $X_t$ lisse, le groupe $H^3(X_t(\mathbb{C}),\mathbb{Z})$ est sans torsion.

(iii) Les fibres singuli\`{e}res de $f$ sont r\'{e}duites
     et ont au pire des singularit\'{e}s quadratiques ordinaires.

(iv) Pour tout   $\alpha \in H^0(\Gamma,R^4f_*\mathbb{Z})$,
     il existe une famille $g_\alpha:B_\alpha\rightarrow \Gamma$  et une famille de $1$-cycles relatifs $\mathcal{Z}_\alpha\subset B_\alpha\times_\Gamma X$ de classe
 $\alpha$, telle que le morphisme
 $\phi_{\mathcal{Z}_\alpha}:B\rightarrow \mathcal{J}^2_\alpha$
     soit surjectif     \`{a} fibre g\'{e}n\'{e}rale rationnellement connexe.

     Alors $Z^4(X)=0$.
\end{theo}
Ce th\'{e}or\`{e}me s'applique concr\`{e}tement aux fibrations en cubiques de dimension $3$ sur une courbe, \`{a} fibres au pire  nodales (cf. \cite{voisinabeljacobi}, o\`u l'on \'{e}tudie d'autres cons\'{e}quences de la condition
(iv)).

Le  th\'{e}or\`{e}me  \ref{theocritereabeljacobi} est une cons\'{e}quence imm\'{e}diate de  l'\'{e}nonc\'{e} suivant:
\begin{prop} \label{critereabeljacobi} Sous les hypoth\`{e}ses
(i), (ii) et (iii)
 ci-dessus, soit
$\alpha$ une section de $R^4f_*\mathbb{Z}$ pour laquelle il existe une famille $B_\alpha$ et un cycle
$\mathcal{Z}_\alpha$ satisfaisant la conclusion de  l'hypoth\`{e}se (iv).
Alors pour toute classe de Hodge $\tilde{\alpha}$ sur $X$, induisant par restriction la section $\alpha$ de $R^4f_*\mathbb{Z}$,
$\tilde{\alpha}$ est la classe d'un cycle alg\'{e}brique sur $X$.
\end{prop}
{\bf D\'{e}monstration.} D'apr\`{e}s Zucker \cite{zucker}, une classe de Hodge enti\`{e}re $\tilde{\alpha}\in Hdg^4(X,\mathbb{Z})\subset H^4(X (\mathbb{C}),\mathbb{Z})$ induit une section $\alpha$ de $R^4f_*\mathbb{Z}$ qui a un rel\`{e}vement canonique
dans la famille tordue de jacobiennes
 $\mathcal{J}^2_\alpha$. Ce rel\`{e}vement est une section alg\'{e}brique  $\sigma:\Gamma\rightarrow\mathcal{J}^2_\alpha$
du morphisme structurel $\mathcal{J}^2_\alpha \rightarrow \Gamma$.
Nous disposons maintenant par hypoth\`{e}se du morphisme $\phi_{\mathcal{Z}_\alpha}:B_\alpha\rightarrow \mathcal{J}^2_\alpha$ qui est alg\'{e}brique,
      surjectif     \`{a} fibre g\'{e}n\'{e}rale rationnellement connexe. Nous appliquons le r\'{e}sultat suivant
       (cf. \cite[Proposition 2.7]{GHMS}):

      \begin{theo} \label{hoxu} Soit $f:W\rightarrow B$ un morphisme surjectif projectif entre vari\'{e}t\'{e}s lisses sur $\mathbb{C}$ dont la fibre g\'{e}n\'{e}rale est rationnellement connexe. Alors pour tout morphisme
$g:C\rightarrow B$, o\`{u} $C$ est une courbe lisse, il existe un rel\`{e}vement $\tilde{g}:C\rightarrow W$
de $g$ dans $W$.
\end{theo}

\medskip
      D'apr\`{e}s ce th\'{e}or\`{e}me, la  section
      $\sigma$ se rel\`{e}ve  en une section $\tilde{\sigma}:\Gamma\rightarrow B_\alpha$.
      Rappelons que $B_\alpha$ param\`{e}tre une famille de $1$-cycles relatifs $\mathcal{Z}_\alpha\subset B_\alpha\times_\Gamma X$ et restreignons cette famille \`{a} $\tilde{\sigma}(\Gamma)$. Ceci fournit
       un cycle $Z\subset X\cong\tilde{\sigma}(\Gamma)\times_\Gamma X$ qui a par construction
      la propri\'{e}t\'{e} que la \og fonction normale \fg \ (cf. \cite{voisinlivre}) $\nu_Z$ associ\'{e}e \`{a} $Z$, d\'{e}finie par
      $$\nu_Z(t)=AJ_{X_t}(Z_{\mid X_t})$$
       est \'{e}gale \`{a} $\sigma$.
      On en d\'{e}duit d'apr\`{e}s \cite{griffiths} (voir aussi \cite[20.2.2]{voisinlivre}),
      gr\^{a}ce \`{a} l'hypoth\`{e}se (ii),
      que les classes de cohomologie $[Z]\in H^4(X(\mathbb{C}),\mathbb{Z})$ de $Z$ et
      $\tilde{\alpha}$ co\"\i ncident apr\`{e}s restriction  \`{a} $X_U$, o\`{u} $U\subset \Gamma_0$ est un ouvert
      affine  de lissit\'{e} de
      $f$.

      Mais le noyau de la fl\`{e}che de restriction $H^4(X(\mathbb{C}),\mathbb{Z})\rightarrow H^4(X_U(\mathbb{C}),\mathbb{Z})$ est engendr\'{e}
      par les $i_{t*}H_4(X_t(\mathbb{C}),\mathbb{Z})$ o\`{u} $t\in \Gamma\setminus U$, et
      $i_t:X_t\rightarrow X$ est l'inclusion. On conclut en notant que gr\^{a}ce \`{a} l'hypoth\`{e}se (iii)
      et au fait que la fibre g\'{e}n\'{e}rale de $f$ satisfait la condition
      $H^2(X_t,\mathcal{O}_{X_t})=0$, les fibres (singuli\`{e}res ou non) $X_t$ ont la propri\'{e}t\'{e} que leur homologie enti\`{e}re de degr\'{e} $4$ est engendr\'{e}e
      par des classes d'homologie de cycles alg\'{e}briques; il en r\'{e}sulte que  $[Z]-\tilde{\alpha}$ est alg\'{e}brique,
      et donc que $\tilde{\alpha}$ est alg\'{e}brique.

\cqfd

\section{Liens entre le groupe $Z_{2}(X)$ et l'indice des vari\'{e}t\'{e}s sur le corps de fonctions d'une courbe sur $\C$}

Soient $\Gamma$ une courbe projective  lisse connexe et $X$  une vari\'{e}t\'{e}  projective lisse de dimension $n$ \'{e}quip\'{e}e d'un morphisme surjectif  $f :X\rightarrow \Gamma$ de fibre g\'{e}n\'{e}rique
$X_\eta$  g\'{e}n\'{e}riquement int\`{e}gre sur $F=\C(\Gamma)$.
 L'indice $I(X_\eta)=I(X_{\eta}/F)$ de $X_\eta$ est d\'{e}fini comme le pgcd des degr\'{e}s
$\deg_{F}(z)$, o\`{u} $z\in CH_0(X_\eta)$. C'est aussi le pgcd des   degr\'{e}s sur $F$ des points ferm\'{e}s
de  la $F$-vari\'{e}t\'{e}  $X_{\eta}$.
Notons $NS(Y)$ le groupe de N\'{e}ron-Severi d'une vari\'{e}t\'{e} projective et lisse $Y$.

\begin{lemm}\label{facteurdirect}
Soit $f:X\rightarrow \Gamma$ un morphisme dominant entre vari\'{e}t\'{e}s connexes, projectives et lisses sur $\C$,
o\`{u} $\Gamma $ est une courbe. On suppose  la fibre g\'{e}n\'{e}rique  ${X}_\eta$ g\'{e}om\'{e}triquement int\`{e}gre
sur son corps de base
  $F=\C(\Gamma)$.
Les conditions suivantes sont \'{e}quivalentes :

(i) On a $I(X_{\eta})=n$.

(ii) L'application compos\'{e}e
$$CH_{1}(X) \to H_2(X(\C),\Z)   \stackrel{f_{*}}{\rightarrow}   H_2(\Gamma(\C),\Z)=\Z$$
a pour image $n\mathbb{Z}$.

 Ceci implique pour un certain entier $m$ divisant $n$ chacune des conditions \'{e}quivalentes suivantes  :

(iii) La projection $p_{*} : H_{2}(X(\C),\Z) \to H_{2}(\Gamma(\C),\Z)$ a pour image $m\mathbb{Z}$.

(iv)  Le  conoyau   de la fl\`{e}che
$p^* : \Z=H^2(\Gamma(\C),\Z(1)) \to H^2(X(\C),\Z(1))/{\rm tors}$  a sa torsion  isomorphe \`{a} $\mathbb{Z}/m\mathbb{Z}$.

(v) Le  conoyau de la fl\`{e}che  $p^* : \Z=NS(\Gamma) \to NS(X)/{\rm tors}$ a sa torsion  isomorphe \`{a} $\mathbb{Z}/m\mathbb{Z}$.

Pour $m=1$, ces \'{e}nonc\'{e}s sont \'{e}quivalents \`{a} chacun des \'{e}nonc\'{e}s suivants :

(vi) La fl\`{e}che $p^* :  \Z=NS(\Gamma) \to NS(X)$ admet une r\'{e}traction.

(vii) Pour tout entier $n>0$, la fl\`{e}che $p^* :  \Z/n=\Pic(\Gamma)/n \to \Pic(X)/n$
est injective.

(viii) Pour tout entier $n>0$, la fl\`{e}che $p^* : \Z/n=H^2_{\et}(\Gamma,\mu_{n}) \to H^2_{\et}(X,\mu_{n})$
est injective.

Ces conditions \`{a} leur tour impliquent que pour tout point $P\in \Gamma(\C)$,
le diviseur $f^{-1}(P)$ n'est pas multiple : le pgcd des degr\'{e}s de ses composantes
est 1.
\end{lemm}

{\bf D\'{e}monstration.}
L'\'{e}quivalence de (i) et (ii) est bien connue : tout point ferm\'{e} $x$ sur la fibre
g\'{e}n\'{e}rique $X_{\eta}$ s'\'{e}tend en un 1-cycle $Z$ sur $X$. De plus, le degr\'{e} de l'extension
$F(x)/F$ est aussi le degr\'{e} de $Z$ sur $\Gamma$.

De (ii) on tire trivialement (iii) pour un certain entier $m$ divisant $n$.

On a une dualit\'{e} parfaite $H^2(X(\C),\Z(1))/{\rm tors} \times H_{2}(X(\C),\Z)/{\rm tors} \to \Z$,
compatible avec les fl\`{e}ches $p^*$ et $p_{*}$. Ceci \'{e}tablit l'\'{e}quivalence de
(iii) et (iv). L'application cycle
$NS(X) \to H^2(X(\C),\Z(1))$ a un conoyau sans torsion (cons\'{e}quence du
th\'{e}or\`{e}me de Lefschetz sur les classes $(1,1)$). Ceci \'{e}tablit l'\'{e}qui\-valence de (iv) et (v).
Enfin  (v) pour $m=1$ est clairement \'{e}quivalent \`{a} (vi).

Pour tout $X/\C$ connexe, projective et lisse, la suite exacte
$$ 0 \to \Pic^0(X) \to \Pic(X) \to NS(X) \to 0,$$
o\`u $\Pic^0(X)$ est divisible donne des isomorphismes
$\Pic(X)/n \simeq NS(X)/n$. Ainsi (vi) implique (vii).
L'\'{e}qui\-valence de (vii) et (viii) r\'{e}sulte de la suite de Kummer
en cohomologie \'{e}tale et de la nullit\'{e} de ${\rm Br}(\Gamma)$.

Sous l'hypoth\`{e}se (viii), les applications $H^2(\Gamma(\C),\Z(1))/n \to H^2(X(\C),\Z(1))/n$ sont injectives.
On en conclut que dans la suite exacte d\'{e}finissant le groupe $K$
$$ 0 \to H^2(\Gamma(\C),\Z(1)) \to H^2(X(\C),\Z(1)) \to K \to 0,$$
o\`u $H^2(\Gamma(\C),\Z(1))=\Z$, le sous-groupe de torsion de $H^2(X(\C),\Z(1))$ est envoy\'{e} isomorphiquement
sur le sous-groupe de torsion de $K$. On a donc la suite exacte
$$0 \to H^2(\Gamma(\C),\Z(1)) \to H^2(X(\C),\Z(1))/{\rm tors} \to K/{\rm tors} \to 0,$$
et donc la condition
 (iv) avec $m=1$
  est satisfaite.

\medskip

Si la fibre en  $P$ est multiple, on a $f^{-1}(P)=nD$ pour un entier $n>1$
et un diviseur $D$. Ainsi $\Z/n=\Pic(\Gamma)/n \to \Pic(X)/n$ n'est pas injectif.
 \cqfd

\begin{prop}\label{propindice0}
Soit $f:X\rightarrow \Gamma$ un morphisme dominant entre vari\'{e}t\'{e}s connexes, projectives et lisses sur $\C$,
o\`{u} $\Gamma $ est une courbe. On suppose  la fibre g\'{e}n\'{e}rique  ${X}_\eta$ g\'{e}om\'{e}triquement int\`{e}gre
sur son corps de base
  $F=\C(\Gamma)$.

(i) Supposons $H^2(X,\O_{X})=0$. Alors on a
$Hdg_{2}(X,\Z)=H_{2}(X(\mathbb{C}),\Z)$, et le groupe
$Z_{2}(X)=H_{2}(X(\C),\Z) /\Im[CH_{1}(X)]$ est fini.

(ii)
Si  le  conoyau de $\Z=NS(\Gamma) \to NS(X)/{\rm tors}$ a sa torsion isomorphe \`{a} $\mathbb{Z}/n\mathbb{Z}$, et
si    de plus  $Z_{2}(X)=0,$
  alors $I(X_{\eta})=n$.

  (iii) Si  le  conoyau de $\Z=NS(\Gamma) \to NS(X)/{\rm tors}$ est sans torsion, et $Z_{2}(X)$
  est d'ordre $n$, alors $I(X_{\eta})$ divise $n$.
\end{prop}

{\bf D\'{e}monstration.}
 De fa\c con g\'{e}n\'{e}rale, le quotient $Hdg_{2}(X,\Z)/\Im[CH_{1}(X)]$ est fini.
Soit $d=\dim(X)$. Sous l'hypoth\`{e}se $H^2(X,\O_{X})=0$, la structure de Hodge sur $H^2(X(\mathbb{C}),\Q(1))$
est triviale. Par dualit\'{e} de Poincar\'{e}, il en est donc de m\^{e}me de la structure
de Hodge sur l'espace vectoriel $H^{2d-2}(X (\mathbb{C}),\Q(d-1)) \simeq H_{2}(X (\mathbb{C}),\Q)$.   On a donc bien $Hdg_{2}(X,\Z)=H_{2}(X(\C),\Z)$.
Ceci \'{e}tablit le point (i).

Le point (ii)  r\'{e}sulte clairement de  (i) et des \'{e}quivalences entre les points (i) et (ii) d'une part, et
 (iii) et (v) d'autre part, du lemme \ref{facteurdirect}.

(iii)  D'apr\`{e}s l'\'{e}quivalence entre les points
(iii) et (v) du  lemme \ref{facteurdirect} et l'\'{e}nonc\'{e} (i) ci-dessus, la premi\`{e}re condition dans (iii) entra\^{\i}ne
qu'il existe une classe de Hodge $\alpha\in Hdg_2(X,\mathbb{Z})$ telle que $p_*\alpha=[\Gamma]$.
Cette classe \'{e}tant annul\'{e}e par $n$ dans le groupe $Z_{2}(X)$, $n\alpha$ est alg\'{e}brique et on a
$p_*(n\alpha)=n[\Gamma]$. Donc $I(X_\eta)$ divise $n$ par l'\'{e}quivalence entre les points (i) et (ii) du lemme
\ref{facteurdirect}.

 \cqfd

 \begin{prop}\label{divisibilite} Soit $f:X\rightarrow \Gamma$ un morphisme dominant entre vari\'{e}t\'{e}s connexes, projectives et lisses,
sur $\C$, o\`{u} $\Gamma $ est une courbe. On suppose que la fibre g\'{e}n\'{e}rique  ${X}_\eta$ est g\'{e}om\'{e}triquement int\`{e}gre
et satisfait $H^i({X}_\eta,\mathcal{O}_{{X}_\eta})=0$  pour $i>0$.
Alors :

(i) Pour tout $i\geq 2$, on a $H^{i}(X,\O_{X})=0$.

(ii)  La fl\`{e}che $f^* : \Z=NS(\Gamma) \to NS(X)/{\rm tors}$ a son conoyau sans torsion.

(iii) Pour tout point $P\in \Gamma(\C)$, le pgcd des multiplicit\'{e}s des composantes
de la fibre $f^{-1}(P)$ est \'{e}gal \`{a} 1.
\end{prop}

{\bf D\'{e}monstration.}
 Les hypoth\`{e}ses sur ${X}_\eta$ impliquent  $R^0f_*\mathcal{O}_{{X}}=\mathcal{O}_\Gamma$ et  $R^if_*\mathcal{O}_{{X}}=0$ pour
$i>0$. L'\'{e}nonc\'{e} (i)  r\'{e}sulte par exemple du fait que
$R^if_*\mathcal{O}_{{X}}$ est g\'{e}n\'{e}riquement nul par hypoth\`{e}se, tandis que Koll\'ar (\cite[Theorem 2.6]{kollarvanishing})
montre que $R^if_*\mathcal{O}_{{X}}$ est sans torsion.
 Par un argument de suite spectrale de Leray, il en r\'{e}sulte que pour tout $L\in \Pic(\Gamma)$, et pour tout $i$, on a:
$$H^i({X},f^*L)=H^i(\Gamma,L).$$
Pour tout $i\geq 2$,  cette \'{e}galit\'{e} implique
$H^i({X},f^*L)=0$
ce qui avec $L=\O_{X}$ \'{e}tablit   la premi\`{e}re assertion.
Cette \'{e}galit\'{e} implique par ailleurs
\begin{eqnarray}
\label{tttruc}\chi({X},f^*L)=\chi(\Gamma,L).
\end{eqnarray}
Pour $L=\mathcal{O}_\Gamma$,
on trouve
 \begin{eqnarray}
\label{truc}
\chi({X},\mathcal{O}_{X})=\chi(\Gamma,\mathcal{O}_\Gamma).
\end{eqnarray}
On applique maintenant le th\'{e}or\`{e}me de Riemann-Roch sur ${X}$.
   On conclut qu'il existe
$Z\in CH_1({X})\otimes \mathbb{Q}$ tel que pour tout $H\in \Pic(X)$ avec $c_1(H)^2=0$ dans $CH^2({X})\otimes\mathbb{Q}$, on ait
\begin{eqnarray}
\label{tructruc}\chi({X},H)=\chi({X},\mathcal{O}_{{X}})+c_1(H)\cdot Z=\chi(\Gamma,\mathcal{O}_\Gamma)+c_1(H)\cdot Z,
\end{eqnarray}
o\`{u} la seconde \'{e}galit\'{e} vient de (\ref{truc}).
Appliquant cette formule \`{a} $H=f^*L,\,L\in \Pic(\Gamma)$, et combinant ceci avec (\ref{tttruc}), on obtient
$$\chi(\Gamma,\mathcal{O}_\Gamma)+c_1(f^*L)\cdot Z=\chi(L)=\chi(\Gamma,\mathcal{O}_\Gamma)+\deg\,L.$$
On en conclut que si $\deg\,L=1$ alors   $c_1(f^*L)\cdot Z=1$.

Supposons maintenant  $f^*L=H^{\otimes k}\otimes M$ dans $\Pic(X)$, avec $k>0$
et $M$ un fibr\'{e} en droites  num\'{e}riquement trivial sur $X$.
La formule de Riemann-Roch (\ref{tructruc}) reste valable pour $H$ et on trouve
$$\chi({X},H)=\chi(\Gamma,\mathcal{O}_\Gamma)+c_1(H)\cdot Z=\chi(\Gamma,\mathcal{O}_\Gamma)+\frac{c_1(f^*L)\cdot Z}{k},$$
et comme $c_1(f^*L)\cdot Z=1$, on trouve que $\chi({X},H)$ n'est pas un nombre entier, ce qui est une contradiction et \'{e}tablit le point (ii).
Ceci implique (iii) (voir le lemme \ref{facteurdirect}).
\cqfd

\begin{rema}

{\rm L'argument ci-dessus montre en particulier
que, sous l'hypoth\`{e}se $H^i({X}_\eta,\mathcal{O}_{{X}_\eta})=0$  pour  tout $i>0$,
aucune fibre de $X \to \Gamma$
n'est multiple. Cela signifie que partout  localement pour la topologie \'{e}tale
sur $\Gamma$, l'indice relatif est 1.
Ce r\'{e}sultat a d\'{e}j\`{a} \'{e}t\'{e} obtenu pr\'{e}c\'{e}demment par
H.~Esnault et O.~Wittenberg
et aussi par J.~Nicaise. }
\end{rema}

\begin{rema}
{\rm Les d\'{e}nominateurs de $Z \in CH_{1}(X)\otimes \mathbb{Q}$ sont connus.
Dans le cas particulier $\dim(X)=3$, on a $Z=(\gamma_{1}^2+\gamma_{2})/12,$
o\`u les $\gamma_{i} \in CH^{i}(X)$ sont les classes de Chern du fibr\'{e} tangent.
On voit donc qu'il existe un 1-cycle $Z$ tel que $\deg(f^*L.Z)=12$.
Ainsi $I(X_{\eta})$ divise 12.  De fa\c con plus g\'{e}n\'{e}rale, si $W/F$
est une surface projective, lisse et g\'{e}om\'{e}triquement connexe sur un corps quelconque $F$,
et si $\chi(W,\O_{W})=1$, alors la formule de Riemann-Roch pour les surfaces
donne $12 = \deg(\gamma_{1}^2+\gamma_{2})$, o\`u  les $\gamma_{i} \in CH^{i}(X)$ sont les classes de Chern du fibr\'{e} tangent.  Ainsi $I(W)$ divise 12.
Le m\^{e}me argument montre que si $W$ est une surface $K3$ sur un corps $F$
alors $I(W)$ divise 24.}

\end{rema}

\begin{theo} \label{propindice1}
 Soit $f:X\rightarrow \Gamma$ un morphisme dominant entre vari\'{e}t\'{e}s connexes, projectives et lisses sur $\C$,
o\`{u} $\Gamma $ est une courbe. On suppose que la fibre g\'{e}n\'{e}rique  ${X}_\eta$ est g\'{e}om\'{e}triquement int\`{e}gre
sur son corps de base
  $F=\C(\Gamma)$.
Si  l'on a $H^{i}(X_{\eta},\O_{X_{\eta}})=0$ pour tout  $i>0$,  alors $I(X_{\eta})$ divise
l'ordre du groupe de torsion $Z_{2}(X)$. En particulier, si $Z_{2}(X)=0$, alors $I(X_{\eta})=1$.
\end{theo}

{\bf D\'{e}monstration.}
Sous l'hypoth\`{e}se  d'annulation des $H^i(X_\eta,\mathcal{O} _{X_\eta})$ pour tout $ i>0$,
d'apr\`{e}s le point (i) de la proposition \ref{divisibilite},
on
a $H^2(X,\mathcal{O}_X)=0$.
L'hypoth\`{e}se (i) de la proposition \ref{propindice0} est donc satisfaite.
La proposition \ref{divisibilite}  assure aussi que sous les hypoth\`{e}ses
 $H^i(X_\eta,\mathcal{O} _{X_\eta})=0$ pour tout $i>0$, la fl\`{e}che
 $\Z=NS(\Gamma) \to NS(X)/{\rm tors}$ a un conoyau sans torsion.
 L'hypoth\`{e}se (iii) de la proposition \ref{propindice0}  est donc bien satisfaite.
 Il ne reste plus qu'\`{a} appliquer cette proposition.
 \cqfd

En dimension 3, on a une r\'{e}ciproque partielle du th\'{e}or\`{e}me \ref{propindice1}.

\begin{theo}\label{propindice2}
 Soit $f:X\rightarrow \Gamma$ un morphisme dominant entre vari\'{e}t\'{e}s connexes, projectives et lisses sur $\C$,
o\`{u} $\Gamma $ est une courbe et $X$ est de dimension 3.
On suppose que  la  fibre g\'{e}n\'{e}rique  ${X}_\eta$, qui est une surface,  est   g\'{e}om\'{e}triquement int\`{e}gre
et satisfait
 $H^i(X_\eta,\mathcal{O} _{X_\eta})=0,\,i=1,\,2$. On suppose en outre que les fibres lisses de $f$ n'ont pas de torsion
dans leur cohomologie  de Betti enti\`{e}re de degr\'{e} $3$, et que les fibres singuli\`{e}res ont au plus
des singularit\'{e}s quadratiques ordinaires. Alors  $Z^4(X)\cong \mathbb{Z}/n\mathbb{Z}$, avec $I(X_\eta)=n$. En particulier  si $I( X_\eta)=1$,  alors $Z^4(X)=0$.
\end{theo}

{\bf D\'{e}monstration.}
Les fibres $X_t$ de $f$ sont de dimension $2$, et l'hypoth\`{e}se
$H^1(X_\eta,\mathcal{O} _{X_\eta})=0$ entra\^{\i}ne
$H^1(X_t,\mathcal{O}_{X_t})=0$  pour toute fibre lisse $X_{t}$,
et donc
 que $H^3(X_t(\mathbb{C}),\Z)$ est de torsion. Ce groupe est donc nul
car sa torsion est nulle par hypoth\`{e}se. Soit $\Gamma^0 \subset \Gamma$
un ouvert affine maximal au-dessus duquel $f$ est lisse, et $X^0=f^{-1}(\Gamma^0)$.
La suite spectrale de Leray
de la restriction \`{a} $f$ \`{a} l'ouvert $X^0$ montre
alors
(du fait que $\Gamma^0$ est affine, et donc a le type d'homotopie d'un CW-complexe de
dimension $\leq1$)
 que $H^4(X^0(\mathbb{C}),\Z)=H^0(\Gamma^0,R^4f_*\mathbb{Z})=\Z$.
 Or on sait par la proposition
\ref{divisibilite}, ii) qu'il existe une classe enti\`{e}re
$\alpha \in H^4(X(\mathbb{C}),\mathbb{Z})$ dont la restriction  aux fibres lisses $X_t$
de $f$ engendre $H^4(X_t(\mathbb{C}),\mathbb{Z})$.  La classe $\alpha$ est une classe de Hodge car
la structure de Hodge sur $H^4(X(\mathbb{C}),\mathbb{Q})$ est triviale du fait que $H^2(X,\mathcal{O}_X)=0$.
 On conclut
que
$H^4(X(\mathbb{C}),\Z)$ est engendr\'{e} par  $\alpha$  et par
$\Ker\, [j_0^*:H^4(X(\mathbb{C}),\mathbb{Z})\rightarrow H^4(X^0(\mathbb{C}),\mathbb{Z})]$,
ce dernier groupe  \'{e}tant  la cohomologie de degr\'{e} $4$ de $X$
support\'{e}e sur les fibres au-dessus de $\Gamma\setminus\Gamma^0$,
c'est-\`{a}-dire le sous-groupe
 $$\oplus_{t\in\Gamma\setminus\Gamma^0}i_{X_t*}H_2(X_t(\mathbb{C}),\Z)\subset H_2(X(\mathbb{C}),\Z)=H^4(X(\mathbb{C}),\Z(2)).$$
Mais sous nos hypoth\`{e}ses sur les fibres singuli\`{e}res on conclut facilement
que leur homologie
$H_2(X_t(\mathbb{C}),\Z)$ est engendr\'{e}e par des classes d'homologie de cycles alg\'{e}briques,  comme c'est le cas
pour les fibres lisses $X_{t}$.
 En effet, les singularit\'{e}s des fibres
$X_t$ sont au pire des singularit\'{e}s
quadratiques ordinaires. Par r\'{e}solution simultan\'{e}e, on conclut
que la d\'{e}singularisation $\widetilde{X}_t$ satisfait encore
$H^2(\widetilde{X}_{t},\mathcal{O}_{X_{t}})=0$. Donc la cohomologie enti\`{e}re (ou encore l'homologie
enti\`{e}re) de degr\'{e}
$2$ de $\widetilde{X}_{t}(\mathbb{C})$ est engendr\'{e}e par des classes de cycles alg\'{e}briques par le th\'{e}or\`{e}me de
Lefschetz sur les classes $(1,1)$. Or la r\'{e}solution par \'{e}clatement d'un point double de surface
 a un diviseur exceptionnel isomorphe \`{a} $\mathbb{P}^1$, et il en r\'{e}sulte que
 la fl\`{e}che
 naturelle  $$H_2(\widetilde{X}_{t}(\mathbb{C}),\mathbb{Z})\rightarrow
 H_2({X}_{t}(\mathbb{C}),\mathbb{Z})$$
 est surjective. Ceci entra\^{\i}ne que l'homologie
enti\`{e}re de degr\'{e}
$2$ de toutes les fibres ${X}_{t}$ est engendr\'{e}e par des classes de cycles alg\'{e}briques.

 On en d\'{e}duit que le groupe $Z^4(X)$ est engendr\'{e} par la classe $\alpha$ ci-dessus et donc est cyclique.
 Que l'ordre de $\alpha$ dans ce groupe soit \'{e}gal \`{a} $I(X_\eta)$ r\'{e}sulte de l'\'{e}quivalence
 des points (i) et (ii) du  lemme \ref{facteurdirect} et du fait que
 $f_*\alpha=[\Gamma]$.

\cqfd

\begin{rema} {\rm Si la surface
 $X_{\eta}$ satisfait
  $H^i(X_\eta,\mathcal{O} _{X_\eta})=0$ pour tout $i>0$, alors la conjecture de Bloch
\cite[lecture 1]{blochbook}
 implique $CH_0(X_{\overline{\eta}})=\Z$.
 Lorsque ceci est satisfait,  sous des hypoth\`{e}ses
moins restrictives, on peut par des m\'{e}thodes de $K$-th\'{e}orie
\'{e}tablir le th\'{e}or\`{e}me \ref{propindice2} sans faire d'hypoth\`{e}ses
  sur les fibres singuli\`{e}res : voir la proposition \ref{indicesurfacesviaKtheorie}
  et le corollaire \ref{famillesurfacesrationnelles}.

Cependant, il  existe
des surfaces simplement connexes (et donc sans torsion dans leur cohomologie enti\`{e}re)
telles que $q=p_g=0$, et pour lesquelles la conjecture de Bloch n'est pas connue. De telles surfaces
(par exemple les surfaces de Barlow g\'{e}n\'{e}riques)
ne sont de plus pas n\'{e}cessairement rigides, pouvant donner lieu \`{a} des familles comme ci-dessus non isotriviales.}
\end{rema}

Au paragraphe \ref{nouveaunonnul}, on a construit
 des  vari\'{e}t\'{e}s $X$ de dimension $3$ avec $H^i(X,\mathcal{O}_X)=0$ pour  $i>0,$
mais $Z^4(X)\not=0$. Ces vari\'{e}t\'{e}s sont fibr\'{e}es en surfaces $K3$ sur une courbe.
Il est naturel de se demander si des exemples analogues existent avec
 des vari\'{e}t\'{e}s  fibr\'{e}es en surfaces d'Enriques ou plus g\'{e}n\'{e}ralement
en surfaces $Y$ avec $H^i(Y,\mathcal{O}_{Y})=0,\,i>0$. D'apr\`{e}s les th\'{e}or\`{e}mes \ref{propindice1}
et \ref{propindice2}, ceci est intimement li\'{e}
aux questions suivantes :
\begin{question}\label{question1}
 Existe-t-il des surfaces  $Y$ d\'{e}finies sur le corps de fonctions d'une courbe
complexe $\Gamma$, telles que
$H^i(Y,\mathcal{O}_Y)=0$ pour $i>0$,  mais $I(Y)\not=1$?
\end{question}

\begin{question}\label{question2}
Existe-t-il de telles surfaces avec de plus $H^3(X_t(\C),\Z)$  sans torsion,  pour une fibre complexe lisse $X_t$ de la fibration
correspondante $X \rightarrow\Gamma$?
\end{question}

Notons que si on remplace la condition $I(Y)\not=1$ par la condition: \og il n'existe pas de point rationnel \fg, la premi\`{e}re de ces questions a une r\'{e}ponse positive, d'apr\`{e}s  \cite{lafon}
et  \cite{GHMS}.
 Dans  \cite{starr}, Starr \'{e}crit que, vraisemblablement, les exemples
de familles de surfaces d'Enriques sans sections construits dans \cite{GHMS} ont un indice diff\'{e}rent de $1$,
mais cela n'est pas montr\'{e} \`{a} notre connaissance.

\section{Quelques r\'{e}sultats obtenus par la $K$-th\'{e}orie alg\'{e}brique}\label{methodesKtheoriques}

Soient $X$ et $Y$ deux vari\'{e}t\'{e}s connexes, projectives et lisses sur $\C$,
et soit $f : X \to Y$ un morphisme dominant \`{a} fibre g\'{e}n\'{e}rale lisse
et connexe. Soit $F=\C(Y)$ le corps des fonctions de $Y$, et soit
$V=X \times_{Y} \Spec(F)$ la fibre g\'{e}n\'{e}rique de $p$.

La comparaison des th\'{e}or\`{e}mes \ref{boresolution} et \ref{boresolutionfinie},
rappel\'{e}e dans   la section \ref{secrappelblochkato},
 donne en particulier
des inclusions
$H^3_{nr}(X,\mu_{n}^{\otimes 2}) \subset H^3_{nr}(V,\mu_{n}^{\otimes 2})$.
Sous certaines hypoth\`{e}ses sur le corps $F$
et sur la $F$-vari\'{e}t\'{e} $V$, des techniques de $K$-th\'{e}orie alg\'{e}brique
permettent de montrer $H^3_{nr}(V,\mu_{n}^{\otimes 2})=0$.
Dans ces cas-l\`{a}, le th\'{e}or\`{e}me \ref{principal1}
donne $Z^4(X)=0$.

\medskip

Soit  d\'{e}sormais $F$ un corps de caract\'{e}ristique z\'{e}ro, $\ovF$ une cl\^{o}ture alg\'{e}brique de $F$,
et $G={\rm Gal}(\ovF/F)$ le groupe de Galois absolu. Soit $V$ une $F$-vari\'{e}t\'{e}
projective, lisse, g\'{e}om\'{e}triquement connexe. On note $\ovV=V\times_{F}\ovF$.
La $K$-th\'{e}orie alg\'{e}brique permet de faire un lien entre les groupes
$H^3_{nr}(V,\Q/\Z(2))$ et le conoyau de l'application
$ CH^2(V) \to CH^2(\ovV)^G.$ Lorsque $V$ est une surface, on trouve
des liens entre $H^3_{nr}(V,\Q/\Z(2))$ et l'indice $I(V)$.

\medskip

Kahn, Rost et Sujatha  \cite[Thm. 5 et Cor. 10 (2)]{KRS})
 ont   \'{e}tabli le  r\'{e}sultat g\'{e}n\'{e}ral suivant.

 \begin{theo}\label{krs}
Soit $Q$   une quadrique  lisse  de dimension au moins 1 sur un corps $F$
de caract\'{e}ristique diff\'{e}rente de 2. Supposons que
cette quadrique n'est pas d\'{e}finie par une  forme  d'Albert
  (forme  quadratique  de rang 6, de la forme
$<a,b,ab, -c,-d,-cd>$).
Alors l'application de restriction
$H^3(F,\Q/\Z(2)) \to H^3_{nr}(Q,\Q/\Z(2))$ (o\`u l'on se limite aux coefficients de
torsion premi\`{e}re \`{a} la caract\'{e}ristique)
est surjective.
 \end{theo}

Ainsi, si $cd(F)\leq 2$, le groupe $ H^3_{nr}(Q,\Q/\Z(2))$ est nul, sauf
peut-\^{e}tre si $Q \subset \mathbb{P}^5_{F}$ est une quadrique d'Albert.

\begin{coro}\label{consequencekrs}
Soit $X \to Y$ un morphisme dominant de vari\'{e}t\'{e}s connexes, projectives et
lisses sur $\C$, de fibre g\'{e}n\'{e}rique une quadrique $Q$ de dimension au moins 1,
de base une surface $Y$. Alors $H^3_{nr}(X,\Q/\Z(2))=0$ et $Z^4(X)=0$.
\end{coro}
{\bf D\'{e}monstration.} Si la dimension de la fibre est au moins 3,
la quadrique $Q$  a un point rationnel sur le corps
 $F=\C(Y)$, qui est un corps
$C_{2}$.  Ceci implique que $Q$ est $F$-birationnelle \`{a} un espace projectif,
et donc  (cf. \cite[Thm. 4.15]{ctsantabarbara})
 que la fl\`{e}che
de restriction $H^3(F,\Q/\Z(2)) \to H^3_{nr}(Q,\Q/\Z(2))$ est un isomorphisme.
Comme la dimension cohomologique de $F$ est 2, on a $H^3(F,\Q/\Z(2))=0$.
Si la dimension de la fibre est 1 ou 2, le th\'{e}or\`{e}me \ref{krs}
et la nullit\'{e} de $H^3(F,\Q/\Z(2))=0$ impliquent $H^3_{nr}(Q,\Q/\Z(2))=0$.
L'inclusion $H^3_{nr}(X,\mu_{n}^{\otimes 2}) \subset H^3_{nr}(V,\mu_{n}^{\otimes 2})$
et le th\'{e}or\`{e}me \ref{principal1} permettent de conclure. \cqfd

\begin{rema}\label{casparticulier}
{\rm
On obtient donc par la $K$-th\'{e}orie
 une d\'{e}monstration  d'un cas particulier du r\'{e}sultat
de C. Voisin sur les solides unir\'{e}gl\'{e}s
(th\'{e}or\`{e}me \ref{theovoisinuniruled} ci-dessus).
La nullit\'{e} de $H^3_{nr}(X,\Z/2)$
pour les solides fibr\'{e}s en coniques
avait \'{e}t\'{e} en substance \'{e}tablie par S. Bloch en 1977 (non publi\'{e}).
Une d\'{e}monstration fond\'{e}e sur les r\'{e}sultats de Merkur'ev et Suslin sur
 la
 $K$-th\'{e}orie des coniques
fut donn\'{e}e en 1989 (voir l'appendice de \cite{parimala}).
}
\end{rema}

 \bigskip

 \subsection{Vari\'{e}t\'{e}s sur un corps de dimension cohomologique 1} \label{varsurcd1}

 Dans cette sous-section, on s'int\'{e}resse de fa\c con g\'{e}n\'{e}rale
 aux vari\'{e}t\'{e}s projectives et lisses
 sur un corps $F$ de caract\'{e}ristique z\'{e}ro et de dimension cohomologique 1,
  le cas principal que nous avons en vue \'{e}tant le corps des fonctions d'une courbe complexe,
  discut\'{e} \`{a} la sous-section \ref{fibrationcourbe}.
 La plupart des \'{e}nonc\'{e}s \'{e}tablis ici   valent encore pour les vari\'{e}t\'{e}s
 projectives et lisses sur un corps fini. Les arguments et r\'{e}sultats sp\'{e}cifiques
 \`{a} cette situation seront trait\'{e}s dans une autre publication.

 Les  consid\'{e}rations qui suivent remontent \`{a} S.~Bloch, elles furent d\'{e}velopp\'{e}es
 par  W.~Raskind et l'un des auteurs dans \cite{ctraskind}, puis dans l'article \cite{kahn} de B. Kahn. Comme on va le voir,
  lorsque le corps $F$ est de dimension cohomologique 1, la m\'{e}thode
  de \cite{ctraskind} suffit pour recouvrer \`{a} moindres frais, et  m\^{e}me g\'{e}n\'{e}raliser,   les r\'{e}sultats de \cite{kahn}.
  Il devrait cependant \^{e}tre clair au lecteur que ce qui suit est inspir\'{e} des r\'{e}sultats de  \cite{kahn},
  et en particulier du Th\'{e}or\`{e}me~1,  de son corollaire, et  de la suite (6)   annonc\'{e}e
  \`{a} la page 397.

\medskip

 \begin{prop}\label{sansbk}
 Soit $F$ un corps de caract\'{e}ristique z\'{e}ro et de dimension cohomologique 1.
Soit $V$ une $F$-vari\'{e}t\'{e} g\'{e}om\'{e}triquement int\`{e}gre.
Soient ${\overline F}$ une cl\^{o}ture alg\'{e}brique de $F$
et $G={\rm Gal}({\overline F}/F)$.

(i) Il existe alors un isomorphisme naturel
$$ \Ker [H^3(F(V),\Q/\Z(2)) \to H^3({\overline F}(V),\Q/\Z(2))]  \simeq H^2(G, K_{2}({\overline F}(V)).$$

(ii) Si $V$ est lisse, cet isomorphisme induit un isomorphisme
$$ \Ker [H^0(V,\H^3(\Q/\Z(2))) \to H^0(\ovV,\H^3(\Q/\Z(2)))]  \simeq
 \Ker[H^2(G, K_{2}({\overline F}(V))) \to H^2(G,\oplus_{x \in \ovV^{(1)} } \ovF(x)^{\times})].$$
 \end{prop}

{\bf D\'{e}monstration.}
Le th\'{e}or\`{e}me de Merkur'ev et Suslin dit que l'application  $G$-\'{e}quivariante
$$K_{2}({\overline F}(V))\otimes \Q/\Z  \to H^2({\overline F}(V),\Q/\Z(2))$$
est un isomorphisme.

On a une suite exacte longue
$$0 \to K_{2}({\overline F}(V))_{\rm tors}  \to K_{2}({\overline F}(V)) \to K_{2}({\overline F}(V))\otimes \Q \to K_{2}({\overline F}(V))\otimes \Q/\Z \to 0$$
que l'on peut d\'{e}couper en deux suites exactes
$$0 \to K_{2}({\overline F}(V))_{\rm tors}  \to K_{2}({\overline F}(V)) \to L \to 0$$
et
$$0 \to L \to K_{2}({\overline F}(V)) \otimes \Q \to K_{2}({\overline F}(V))\otimes \Q/\Z \to 0.$$
La seconde suite donne
$H^1(G, K_{2}({\overline F}(V))\otimes \Q/\Z )  \oi  H^2(G,L).$
La premi\`{e}re suite donne un homomorphisme
$H^2(G,K_{2}({\overline F}(V))) \to H^2(G,L)$
 qui est
bijectif car  $cd(F) \leq 1$. On a donc un isomorphisme
$$H^2(G,K_{2}({\overline F}(V))) \simeq H^1(G, K_{2}({\overline F}(V))\otimes \Q/\Z ).$$
La
 suite spectrale de Hochschild-Serre
 $$E_{2}^{pq} = H^p(G,H^q({\overline F}(V), \Q/\Z(2)) \Longrightarrow H^n(F(V),\Q/\Z(2)).$$
  donne naissance \`{a} un homomorphisme
$$  \Ker [H^3(F(V),\Q/\Z(2)) \to H^3({\overline F}(V),\Q/\Z(2))]  \to H^1(G, H^2({\overline F}(V),\Q/\Z(2)) )$$
qui est
un isomorphisme car $cd(F) \leq 1$. Ceci \'{e}tablit la premi\`{e}re partie de la proposition.

\medskip

Supposons $V$ lisse sur $F$. Soit $x$ un point de $V^{(1)}$.
Comme on l'a rappel\'{e} au paragraphe \ref{secrappelblochkato},
on dispose aussi pour l'anneau semilocal ${\mathcal O}_{\ovV,x}$
(semilocalis\'{e} de $\ovV$ aux points de codimension 1 de $\ovV$ d'image $x \in V$)
 d'un isomorphisme $G$-\'{e}quivariant
 $$K_{2}({\mathcal O}_{\ovV,x})\otimes \Q/\Z  \oi  H^2_{\et}({\mathcal O}_{\ovV,x},\Q/\Z(2)).$$
 Proc\'{e}dant comme ci-dessus, on obtient un isomorphisme
 $$ \Ker [H^3_{\et}({\mathcal O}_{V,x},\Q/\Z(2)) \to H^3_{\et}({\mathcal O}_{\ovV,x},\Q/\Z(2))]  \simeq H^2(G, K_{2}({\mathcal O}_{\ovV,x})).$$

 Pour tout point $x \in V^{(1)}$, la conjecture de Gersten pour l'anneau semilocal ${\mathcal O}_{\ovV,x}$
(cf. \cite{ctkh}) donne une suite exacte courte de $G$-modules
 $$0 \to K_{2}({\mathcal O}_{\ovV,x}) \to K_{2} ({\overline F}(V)) \to {\overline F}(x)^{\times} \to 0,$$
 o\`u l'on a not\'{e} $\ovF(x)=\ovF \otimes_{F}F(x)$.
 Ainsi le noyau de
 $$H^2(G,K_{2} ({\overline F}(V))) \to   H^2(G,\oplus_{x \in V^{(1)} } \ovF(x)^{\times})$$
 est le groupe des \'{e}l\'{e}ments de $H^2(G,K_{2} ({\overline F}(V))) $
 qui sont dans l'image de $H^2(G,K_{2}({\mathcal O}_{\ovV,x}))$ pour tout $x \in V^{(1)}$.

 La conjecture de Gersten pour la cohomologie \'{e}tale (Bloch--Ogus, voir le th\'{e}or\`{e}me \ref{boresolution})
 implique que le groupe $H^0(V,\H^3(\Q/\Z(2))) \subset H^3(F(V),\Q/\Z(2))$ est le sous-groupe
 des \'{e}l\'{e}ments de $H^3(F(V),\Q/\Z(2))$ qui en chaque point $x \in V^{(1)}$ sont dans
 l'image de $H^3_{\et}({\mathcal O}_{V,x}, \Q/\Z(2))$.

 Ceci \'{e}tablit la deuxi\`{e}me partie de la proposition.
 \cqfd

\begin{theo}\label{cdFleq1}
Soit  $F$ un corps de caract\'{e}ristique z\'{e}ro et de dimension cohomologique $\leq 1$.
Soit $V$ une $F$-vari\'{e}t\'{e} projective, lisse, g\'{e}om\'{e}triquement int\`{e}gre.
Soit $\ovF$ une cl\^{o}ture alg\'{e}brique de $F$ et  $\ovV=V\times_{F}\ovF$.
On a
alors une suite exacte
 $$0 \to \Ker [CH^2(V) \to CH^2(\ovV)]  \to  H^1(G, H^1(\ovV, \K_{2})     ) \to $$ $$
\Ker[ H^3_{nr}(V,\Q/\Z(2)) \to H^3_{nr}(\ovV,\Q/\Z(2))]
\to
  \Coker [CH^2(V) \to CH^2(\ovV)^G] \to 0.$$
\end{theo}

{\bf D\'{e}monstration.}
Soit $F$ un corps de caract\'{e}ristique z\'{e}ro, $E/F$ une extension galoisienne,
et $G={\rm Gal}(E/F)$ le groupe de Galois.
  Soit $V$ une $F$-vari\'{e}t\'{e} lisse g\'{e}om\'{e}triquement int\`{e}gre.

On consid\`{e}re le complexe, gradu\'{e} en degr\'{e}s  (0, 1, 2)
$$ 0 \to K_{2}E(V) \to
 \oplus_{      x  \in  V_{E}^{(1)}      }  E(x)^{\times} \to
\oplus_{x \in V_{E}^{(2)  }     } \Z  \to 0.$$

Par la conjecture de Gersten pour la $K$-th\'{e}orie, \'{e}tablie par Quillen
(Th\'{e}or\`{e}me \ref{quillenresolution})
l'homologie de ce complexe en degr\'{e} $i$ est $H^{i}(V_{E}, \K_{2})$.
Notons $\mathcal Z$ le noyau de la deuxi\`{e}me fl\`{e}che,
et $\cal I$ son image.

On a donc des suites exactes de modules galoisiens
$$ 0 \to {\cal Z} \to \oplus_{x \in V_{E}^{(1)} }E(x)^{\times} \to {\cal I} \to 0,$$
$$0 \to K_{2}E(V)  /H^0(V_{E},\K_{2}) \to {\cal Z} \to H^1(  V_{E},\K_{2}) \to 0$$
 $$0 \to {\cal I}  \to \oplus_{x \in   V_{E}^{(2)} } \Z \to CH^2(  V_{E}) \to 0.$$
et les suites analogues avec $F$ en place de $E$.

On sait que l'on a
$\oplus_{x \in V^{(2)} } \Z = (\oplus_{x \in  V_{E}^{(2)} } \Z)^G$, que l'on a
$H^1(G,\oplus_{x \in  V_{E}^{(2)} } \Z)=0$ (lemme de Shapiro et nullit\'{e} de $H^1(G,\Z)$ pour
$G$ profini)
et qu'enfin  on a $H^1(G,  \oplus_{x \in  V_{E}^{(1)} } E(x)^{\times})=0$ (lemme de Shapiro et th\'{e}or\`{e}me 90 de Hilbert).

\medskip

De ceci on d\'{e}duit les suites exactes
$$0 \to H^1(G,{\cal Z}) \to CH^2(V) \to CH^2(V_{E})^G \to H^1(G,{\cal I}) \to 0$$
$$0 \to H^1(G,{\cal I}) \to H^2(G,{\cal Z}) \to
H^2(G,
 \oplus_{    x \in  {V_E}^{(1)}   }
   E(x)^{\times}).$$

Dans la suite on prend $E=\ovF$.

\medskip

Sur tout corps $F$ de caract\'{e}ristique z\'{e}ro, on sait
  (\cite[Thm. 1.8]{ctraskind}) que le groupe $H^0(\ovV,\K_{2})$ est extension d'un groupe fini
par un groupe divisible. On a le r\'{e}sultat analogue pour  $H^1(\ovV,\K_{2})$ (\cite[Thm. 2.2]{ctraskind}).
La d\'{e}monstration de ces r\'{e}sultats repose sur des travaux ant\'{e}rieurs de Bloch, Merkur'ev--Suslin, et Suslin,
et utilise  les r\'{e}sultats de Deligne sur les conjectures de Weil.

Sous l'hypoth\`{e}se $cd(F) \leq 1$, cette  structure des groupes  $H^{i}(\ovV,\K_{2})$  pour $i=0,1$ implique
$H^r(G,H^{i}(\ovV,\K_{2}))=0$ pour $i=0,1$ et
  $r\geq 2$.

On a $H^1(G,K_{2}\ovF(V)/H^0(\ovV,\K_{2})) =0$
car $H^1(G,K_{2}\ovF(V))=0$ si $cd(F) \leq 1$
(\cite[Thm. B]{ctHilbert90}),
et la fl\`{e}che
$H^2(G,K_{2}\ovF(V)) \to H^2(G,K_{2}\ovF(V)/H^0(\ovV,\K_{2}))$
est un isomorphisme.

La cohomologie galoisienne des suites exactes ci-dessus
donne alors   la suite exacte
 $$0 \to H^1(G,{\cal Z}) \to
 H^1(G,H^1(\ovV, \K_{2})) \to H^2(G, K_{2} \ovF(V)) ) \to H^2(G, {\cal Z} ) \to  0.$$
et donc la suite exacte
$$0 \to \Ker [CH^2(V) \to CH^2(\ovV)]  \to  H^1(G,H^1(\ovV, \K_{2})) \to H^2(G, K_{2}\ovF(V)) \to H^2(G, {\cal Z}) \to 0.$$
La fl\`{e}che
$H^2(G, K_{2}\ovF(V)) \to H^2(G,\oplus_{x \in \ovV^{(1)} } \ovF(x)^{\times})$
est induite par les r\'{e}sidus. Elle induit    la fl\`{e}che
$H^2(G, {\cal Z}) \to H^2(G,\oplus_{x \in \ovV^{(1)} } \ovF(x)^{\times}).$
On a donc   la suite exacte
$$0 \to \Ker [CH^2(V) \to CH^2(\ovV)]  \to  H^1(G,H^1(\ovV, \K_{2})) \to  \hskip7cm
$$ $$  \Ker[H^2(G, K_{2}\ovF(V))  \to H^2(G,\oplus_{x \in  \ovV^{(1)} } \ovF(x)^{\times})]
\to \Ker[H^2(G, {\cal Z}) \to H^2(G,\oplus_{x \in \ovV^{(1)} } \ovF(x)^{\times})] \to 0,$$
soit encore
$$0 \to \Ker [CH^2(V) \to CH^2(\ovV)]  \to  H^1(G,H^1(\ovV, \K_{2})) \to \hskip6cm$$ $$
\hskip2cm  \Ker[H^2(G, K_{2}\ovF(V))   \to H^2(G,\oplus_{x \in  \ovV^{(1)} } \ovF(x)^{\times})]
\to \Coker [CH^2(V) \to CH^2(\ovV)^G]  \to 0.$$
On conclut avec la proposition \ref{sansbk}. \cqfd

\bigskip

Pour tirer du th\'{e}or\`{e}me des cons\'{e}quences pratiques, il faut contr\^{o}ler le module galoisien
$H^1(\ovV, \K_{2}) $. Le th\'{e}or\`{e}me suivant regroupe des r\'{e}sultats de \cite{ctraskind}
 (Thm. 2.1, Thm. 2.2, Thm. 2.12).

\begin{theo} (Colliot-Th\'{e}l\`{e}ne et Raskind)
Soit  $F$ un corps de  caract\'{e}ristique z\'{e}ro.
Soit $V$ une $F$-vari\'{e}t\'{e} projective, lisse, g\'{e}om\'{e}triquement int\`{e}gre.
Soit $M=M(\ovV)$ le module galoisien fini d\'{e}fini par
$M =\oplus_{l} H^3_{\et}(\ovV,\Z_{l}(2))\{l\}$.

(i) Il existe une suite exacte naturelle
$$0 \to D \to H^1(\ovV,\K_{2}) \to M \to 0,$$
o\`u le groupe de droite est fini et le groupe $D$ est divisible.

(ii) Pour tout premier $l$ il existe un isomorphisme naturel
$$H^2_{\et}(\ovV, \Q_{l}/\Z_{l}(2)) \oi H^1(\ovV,\K_{2})\{l\}.$$

(iii) Soit $K$, resp. $C$, le noyau, resp. le conoyau de
la fl\`{e}che naturelle $ \Pic(\ovV) \otimes {\overline F}^{\times} \to  H^1(\ovV,\K_{2}).$
 Si l'on a $H^2(V,{\mathcal O}_{V})=0$, alors $K$ est uniquement divisible et
 $C$ est  la somme directe d'un groupe uniquement divisible et du groupe
$M$.
\end{theo}

 Sous l'hypoth\`{e}se  $H^2(V,{\mathcal O}_{V})=0$, on a donc une suite exacte de $G$-modules
$$0  \to K \to \Pic({\overline V}) \otimes\overline{k}^* \to H^1({\overline V},K_{2}) \to C \to 0,$$
o\`u $K$ est un groupe uniquement divisible et $C$
est somme directe de $M$ et d'un groupe uniquement divisible.

Ceci donne naissance \`{a} deux suites exactes
$$0  \to K \to \Pic({\overline V}) \otimes\overline{k}^*  \to R \to 0$$
$$0 \to R  \to H^1({\overline V},K_{2}) \to C \to 0,$$
o\`u $R$ est divisible.

On a aussi la suite exacte
$$0 \to \Pic^0({\overline V}) \to \Pic({\overline V})  \to NS({\overline V}) \to 0 $$
et donc la suite exacte
$$ Tor^1(NS({\overline V}), \overline{k}^*) \to \Pic^0(\overline{V}) \otimes  \overline{k}^*
 \to \Pic({\overline V}) \otimes \overline{k}^*  \to NS({\overline V}) \otimes \overline{k}^*  \to 0.$$

On en d\'{e}duit les suites exactes
$$ 0 \to P \to  \Pic^0({\overline V} ) \otimes \overline{k}^*  \to Q \to 0$$
$$0 \to Q \to \Pic({\overline V} ) \otimes \overline{k}^*  \to NS({\overline V}) \otimes \overline{k}^*  \to 0$$
o\`u $P$ est un groupe fini et $Q$  est divisible.

Si l'on tient maintenant compte des faits suivants :

\smallskip

(1) $\Pic^0({\overline V} )\otimes \overline{k}^*$ est uniquement divisible car c'est un produit tensoriel de deux groupes divisibles;

(2) pour un corps $k$ avec $cd(k) \leq 1$ et  $W$  un module galoisien sur $k$,
 on a $H^2(G,W)=0$ si  $W$ est fini ou divisible, et $H^{r}(G,W)=0$ pour
$r\geq 1$ si $W$ est uniquement divisible,

(3) $H^1(G, NS({\overline V}) \otimes\overline{k}^*)=H^1(G, NS({\overline V})/tors \otimes\overline{k}^*) = 0$
car pour  $k$ avec $cd(k) \leq 1$ et $T$ un $k$-tore, $H^1(G,T({\overline k}))=0$,

\smallskip

\noindent alors on obtient  $H^1(G,Q)=0$, puis $H^1(G,  \Pic({\overline V} ) \otimes\overline{k}^*)=0$,
puis  $H^1(G,R)=0$ et $H^2(G,R)=0$, et
en fin de compte $$H^1(G, H^1({\overline V},K_{2}))  \simeq H^1(G,M).$$

En combinant ce r\'{e}sultat avec le th\'{e}or\`{e}me \ref{cdFleq1} , on obtient le th\'{e}or\`{e}me suivant.

\begin{theo}\label{H2OXtrivial}
Soit $F$ un corps de caract\'{e}ristique z\'{e}ro et de dimension cohomologique au plus~1.
Soit $\ovF$ une cl\^{o}ture alg\'{e}brique de $F$, et $G$ le groupe de Galois de $\ovF$ sur $F$.
Soit $V$ une $F$-vari\'{e}t\'{e} projective, lisse, g\'{e}om\'{e}triquement int\`{e}gre,
satisfaisant $H^2(V,\O_{V})=0$.
Soit $M$ le module galoisien fini
$\oplus_{l} H^3_{\et}(\ovV,\Z_{l}(2))\{l\}$.
On a alors une suite exacte
$$0 \to \Ker [CH^2(V) \to CH^2(\ovV)]  \to H^1(G,M)   \to \hskip7cm
 $$ $$\hskip2cm \Ker[ H^3_{nr}(V,\Q/\Z(2)) \to H^3_{nr}(\ovV,\Q/\Z(2))]
\to   \Coker[CH^2(V) \to CH^2(\ovV)^G] \to 0.$$
\end{theo}

Le groupe de Brauer $Br(\ovV)$ est
une extension du groupe fini $\oplus_{l} H^3_{\et}(\ovV,\Z_{l}(1))\{l\} $
par le groupe $(\Q/\Z)^{b_{2}-\rho}$ (\cite{grothendieck}, voir le paragraphe \ref{groupebrauer} ci-dessus).
En caract\'{e}ristique z\'{e}ro, $b_{2}-\rho=0$
 si et seulement si $H^2(V,\O_{V})=0$.

\begin{theo}\label{brauertrivial}
Soit $F$ un corps de caract\'{e}ristique z\'{e}ro et de dimension cohomologique au plus~1.
Soit $\ovF$ une cl\^{o}ture alg\'{e}brique de $F$, et $G=Gal(\ovF/F)$.
Soit $V$ une $F$-vari\'{e}t\'{e} projective et lisse, g\'{e}om\'{e}triquement int\`{e}gre.

(i) Si le groupe de Brauer de $\ovV$ est nul,
alors on a une suite exacte
$$0 \to CH^2(V) \to CH^2(\ovV)^G \to H^3_{nr}(V,\Q/\Z(2)) \to H^3_{nr}(\ovV,\Q/\Z(2)).$$

(ii) Si $\ovV$ est une vari\'{e}t\'{e} rationnelle,  on a une suite exacte
$$0 \to CH^2(V) \to CH^2(\ovV)^G \to H^3_{nr}(V,\Q/\Z(2)) \to 0.$$

(iii) Si $\ovV$ est une vari\'{e}t\'{e} rationnellement connexe de dimension 3,
sans torsion dans sa cohomologie enti\`{e}re de degr\'{e} 3,
 alors on a une suite exacte
$$0 \to CH^2(V) \to CH^2(\ovV)^G \to H^3_{nr}(V,\Q/\Z(2)) \to 0.$$
\end{theo}

{\bf D\'{e}monstration.}   L'\'{e}nonc\'{e} (i) est une cons\'{e}quence imm\'{e}diate
du th\'{e}or\`{e}me \ref{H2OXtrivial} et du rappel ci-dessus sur le groupe de Brauer.
Si  $\ovV$ est rationnelle, son groupe de Brauer
est nul et  l'invariant birationnel  $H^3_{nr}(\ovV,\Q/\Z(2))$
 est aussi  nul.  Ceci \'{e}tablit (ii).
 Si $\ovV$ est rationnellement connexe de dimension 3, alors
 $H^2(V,\O_{V})=0$, donc le groupe de Brauer de $\ovV$ est, sous l'hypoth\`{e}se de (iii), nul.
 Par ailleurs pour $\ovV$ rationnellement connexe de dimension 3, on a $H^3_{nr}(\ovV,\Q/\Z(2))=0$
 (\cite{voisinuniruled}, voir le corollaire \ref{corodevoisinuniregle} ci-dessus). Ceci \'{e}tablit (iii). \cqfd

 \begin{prop}\label{indicesurface}
Soit  $F$ un corps de caract\'{e}ristique z\'{e}ro et de dimension cohomologique $\leq 1$.
Soit $V$ une $F$-surface projective, lisse, g\'{e}om\'{e}triquement int\`{e}gre.
Supposons $H^1(V,\O_{V})=0.$ Si $H^3_{nr}(V,\Q/\Z(2))=0$,
 alors
l'indice $I(V)$ est \'{e}gal \`{a} 1.
\end{prop}

{\bf D\'{e}monstration.} Le degr\'{e} d\'{e}finit la suite exacte de modules galoisiens
$$ 0 \to A_{0}(\ovV) \to CH_{0}(\ovV) \to \Z \to 0.$$
 L'hypoth\`{e}se $H^1(V,\O_{V})=0$ implique
que la vari\'{e}t\'{e} d'Albanese de $V$ est nulle. Il r\'{e}sulte alors
du th\'{e}or\`{e}me de Roitman (cf. \cite[\S 4]{ctcime})
que le groupe $A_{0}(\ovV) $ est uniquement divisible.
Ainsi $H^1(G,A_{0}(\ovV))=0$ et
l'application $CH_{0}(\ovV)^G \to \Z$ est  surjective.
Le th\'{e}or\`{e}me \ref{cdFleq1}
 et l'hypoth\`{e}se $H^3_{nr}(V,\Q/\Z(2))=0$
impliquent que la fl\`{e}che $CH_{0}(V) \to CH_{0}(\ovV)^G$ est surjective.
On a donc $I(V)=1$. \cqfd

Voici un r\'{e}sultat un peu plus pr\'{e}cis que
la proposition \ref{indicesurface}.

\begin{prop}\label{H1etH3nrmodn}
Soit $F$ un corps de caract\'{e}ristique z\'{e}ro, de dimension cohomologique   1.
Soit $V$ une surface sur $F$,
projective, lisse,
g\'{e}om\'{e}triquement connexe.
  Soit $\ovF$ une cl\^{o}ture alg\'{e}brique de $F$ et $\ovV=V \times_{F}\ovF$.
  Soit $n>1$ un entier.
  Supposons

  (i) $H^1(V,\O_{V})=0$.

  (ii)  $NS(\ovV)[n]=0$.

  (iii) Le groupe de cohomologie \'{e}tale non ramifi\'{e}e $H^3_{nr}(V,\mu_{n}^{\otimes 2})$
  est nul.

  Alors l'indice $I(V)$ est premier \`{a} $n$.
  \end{prop}

  {\bf D\'{e}monstration.}
  La th\'{e}orie de Bloch--Ogus donne naissance \`{a} une suite exacte (\S \ref{secrappelblochkato})
  $$H^3_{\et}(V,\mu_{n}^{\otimes 2}) \to     H^3_{nr}(V,\mu_{n}^{\otimes 2}) \to CH^2(V)/n \to H^4_{\et}(V,\mu_{n}^{\otimes 2}).$$
  L'hypoth\`{e}se (iii) implique donc que l'application cycle
  $ CH^2(V)/n \to H^4_{\et}(V,\mu_{n}^{\otimes 2})$ est injective.
  Pour toute surface $V/F$, l'application
  $CH^2(\ovF)/n \to H^4(\ovV,\mu_{n}^{\otimes 2})=\Z/n$ est
  un isomorphisme.   La suite spectrale de Leray pour $V/F$ et la cohomologie
  \'{e}tale, sous l'hypoth\`{e}se $cd(F) \leq 1$, donne une suite
  exacte
  $$0 \to H^1(F,H^3(\ovV,\mu_{n}^{\otimes 2})) \to H^4(V,\mu_{n}^{\otimes 2})
  \to H^4(\ovV,\mu_{n}^{\otimes 2})^G \to 0.$$
  On a $H^4(\ovV,\mu_{n}^{\otimes 2})^G= H^4(\ovV,\mu_{n}^{\otimes 2})=\Z/n$.
 Les hypoth\`{e}ses (i) et (ii) donnent
  $H^1(\ovV,\mu_{n})=\Pic(\ovV)[n]=0$. Par dualit\'{e} de Poincar\'{e}
  en cohomologie \'{e}tale
  on en d\'{e}duit $H^3(\ovV,\mu_{n}^{\otimes 2})=0$.
  On a donc $H^4(V,\mu_{n}^{\otimes 2})
  \oi H^4(\ovV,\mu_{n}^{\otimes 2})$. Du diagramme commutatif
  $$\begin{array}{ccccccccccccccccccccccccccccccccccccccccccccccccccc}
CH^2(V)/n & \hookrightarrow H^4(V,\mu_{n}^{\otimes 2}) \cr
\downarrow{}{} &  \downarrow{}{\simeq} \cr
CH^2(\ovV)/n & \hookrightarrow H^4(\ovV,\mu_{n}^{\otimes 2})
  \end{array}$$
  on d\'{e}duit que l'application degr\'{e} $CH_{0}(V) \to \Z$
  dont l'image est $\Z.I(V)$
  induit un plongement $CH_{0}(V)/n \hookrightarrow \Z/n$
  et donc un plongement $\Z.I/\Z.(nI) \hookrightarrow \Z/n$. \cqfd

 \begin{rema} {\rm
 Dans \cite{ctmadore}, on donne un exemple
 surface cubique lisse $V$,
 sur un (grand) corps $F$ de dimension cohomologique 1,
 avec $I(V)=3$. Dans ce cas on a donc $H^3_{nr}(V,\Z/3) \neq 0$.}
  \end{rema}

   \begin{rema}
 {\rm
Soit  $S\subset \mathbb{P}^3_{\C}$ une surface quintique invariante sous l'action de Godeaux
de $G=\mathbb{Z}/5\mathbb{Z}$, $X=\mathbb{P}^1\times S/G$, o\`u l'action de $G$ est diagonale.
Alors $\pi : X\rightarrow \Gamma=\mathbb{P}^1/G$ a ses fibres lisses isomorphes \`{a} $S$
et donc la fibre g\'{e}n\'{e}rique $V=X_{\eta}$ satisfait $H^1(V,\O_{V})=0$
  et  $NS(\ovV)_{tors}=0$.
La fibration $\pi$  a des fibres de multiplicit\'{e} 5;  on en d\'{e}duit $I(X_{\eta})=5$.
D'apr\`{e}s la proposition  \ref{H1etH3nrmodn}, on  a $H^3_{nr}(V,\Z/5) \neq 0$.
Cependant la vari\'{e}t\'{e} $X$ est fibr\'{e}e en coniques au-dessus de la surface $S/G$.
Par  le corollaire \ref{consequencekrs}, on a donc   $H^3_{nr}(X,\Q/\Z(2))=0$ et  $Z^4(X)=0$.
Dans ce cas on voit donc que l'inclusion $H^3_{nr}(X,\Z/5) \subset H^3_{nr}(V,\Z/5)$
n'est pas un isomorphisme.}
 \end{rema}

\begin{prop}\label{indicesurfacesviaKtheorie}
Soit $F$ un corps de caract\'{e}ristique z\'{e}ro et de dimension cohomologique au plus~1.
Soit $\ovF$ une cl\^{o}ture alg\'{e}brique de $F$, et $G=Gal(\ovF/F)$.
Soit $V$ une $F$-surface projective et lisse, g\'{e}om\'{e}triquement int\`{e}gre.
Soit $I(V)$ l'indice de $V$ sur $F$.
Supposons   $H^2(V,\O_{V})=0$ et $NS(\ovV)_{\rm tors}=0$.
Alors on  a une suite exacte
$$0 \to CH_{0}(V) \to CH_{0}(\ovV)^G \to H^3_{nr}(V,\Q/\Z(2)) \to 0.$$

Supposons de plus $H^1(V,\O_{V})=0$. Alors
le quotient du groupe $H^3_{nr}(V,\Q/\Z(2)) $
par son sous-groupe divisible maximal est isomorphe \`{a} $\Z/I(V)$.

Si de plus $\deg : CH_{0}(\ovV) \to \Z$ est un isomorphisme
 (ce qui sous les hypoth\`{e}ses pr\'{e}c\'{e}dentes est une conjecture de Bloch), alors
$H^3_{nr}(V,\Q/\Z(2)) \oi \Z/I(V)$.
\end{prop}
{\bf D\'{e}monstration.}  On applique le th\'{e}or\`{e}me \ref{brauertrivial}(i).
Les hypoth\`{e}ses en sont satisfaites, car pour une surface projective et lisse $\ovV$,
le groupe $NS(\ovV)_{\rm tors}$ est nul si et seulement si  on a $\oplus_{l} H^3_{\et}(\ovV,\Z_{l}(2))\{l\}=0$.

Comme $V$ est une surface,  le corps $\ovF(V)$ est de dimension cohomologique~2,
on a donc $H^3_{nr}(\ovV,\Q/\Z(2)) \subset H^3(\ovF(V),\Q/\Z(2))=0$.
Sous l'hypoth\`{e}se $H^1(V,\O_{V})=0$, la vari\'{e}t\'{e} d'Albanese de $V$
est triviale. Il r\'{e}sulte alors du th\'{e}or\`{e}me de Roitman
(cf. \cite[ \S 4]{ctcime})
que le  groupe $A_{0}(\ovV)$ noyau de la fl\`{e}che degr\'{e} : $CH_{0}(\ovV) \to \Z$
est un groupe
uniquement divisible, i.e. est un $\Q$-vectoriel. Ceci implique que le groupe $A_{0}(\ovV)^G $
est  un $\Q$-vectoriel et que l'on a
 $H^1(G,A_{0}(\ovV))=0$.
Le lemme du serpent donne  alors la suite exacte
$$ A_{0}(\ovV)^G \to H^3_{nr}(V,\Q/\Z(2)) \to \Z/I(V) \to 0.$$
  \cqfd

\subsection{Applications aux fibrations au-dessus d'une courbe}\label{fibrationcourbe}

\begin{theo}\label{famillesurfacesrationnelles}
Soit $\Gamma$ une courbe projective et lisse sur $\C$ et
$X$ un solide projectif et lisse muni d'une fibration $X \to \Gamma$
dont la fibre g\'{e}n\'{e}rique g\'{e}om\'{e}trique est une surface rationnelle.
Alors les groupes $H^3_{nr}(X,\Q/\Z(2))$ et $Z^4(X)$ sont nuls.
\end{theo}
{\bf D\'{e}monstration.}
Soient $F=\C(\Gamma)$ et $V/F$ la fibre g\'{e}n\'{e}rique de  $X \to \Gamma$.
Le corps $F$, corps de fonctions d'une variable sur $\C$, est un corps $C_{1}$.
Rappelons le fait bien connu suivant.
De la  classification $F$-birationnelle des
$F$-surfaces g\'{e}om\'{e}triquement rationnelles (Enriques, Manin, Iskovskikh; Mori)
on d\'{e}duit par inspection que toute telle surface $V$ sur un corps $C_{1}$ poss\`{e}de un point
rationnel. Sur $F=\C(\Gamma)$, c'est aussi un cas particulier d'un th\'{e}or\`{e}me de
Graber, Harris et Starr \cite{GHS}.
En particulier l'application $CH^2(V) \to CH^2(\ovV)^G=\Z$
est surjective. De la proposition \ref{indicesurfacesviaKtheorie}
 on d\'{e}duit alors
$H^3_{nr}(V,\Q/\Z(2))=0$. Comme rappel\'{e} au d\'{e}but du paragraphe  \ref{methodesKtheoriques}, le groupe $H^3_{nr}(X,\Q/\Z(2))$ est un sous-groupe
de $H^3_{nr}(V,\Q/\Z(2))$, il est donc nul. L'\'{e}nonc\'{e} sur $Z^4(X)$ est alors
une cons\'{e}quence du th\'{e}or\`{e}me \ref{principal1}.
\cqfd

\begin{rema}{\rm
La d\'{e}monstration de ce th\'{e}or\`{e}me  ne repose pas sur l'\'{e}nonc\'{e}  \ref{corodevoisinuniregle}.
On obtient ainsi par la $K$-th\'{e}orie alg\'{e}brique un nouveau  cas particulier (voir la remarque \ref{casparticulier})  du th\'{e}or\`{e}me
de C. Voisin sur les solides unir\'{e}gl\'{e}s (th\'{e}or\`{e}me \ref{theovoisinuniruled} ci-dessus). }
\end{rema}

\begin{lemm} \label{9mai}
Soit   $\pi:X\rightarrow \Gamma$ un morphisme dominant de vari\'{e}t\'{e}s projectives et lisses
sur $\C$, avec $\Gamma$ une courbe. Soit $n$ la dimension de la fibre g\'{e}n\'{e}rique $V$
suppos\'{e}e g\'{e}om\'{e}triquement int\`{e}gre sur le corps $\C(\Gamma)$. Supposons
$H^2(V,\O_{V})=0$.

(i)  Si une fibre singuli\`{e}re $Y$ de $\pi$ a des singularit\'{e}s quadratiques ordinaires, on a $H^2(Z,\O_{Z})=0$, o\`{u} $Z$ est la d\'{e}singularis\'{e}e de $Y$ obtenue par \'{e}clatement de ses points doubles.

(ii) Si de plus une fibre lisse  de $\pi$ a sa cohomologie de Betti
enti\`{e}re de degr\'{e} 3  sans torsion, et $n\not=3$, il en va de m\^{e}me pour $Z$.
\end{lemm}
{\bf D\'{e}monstration.} (i) D'apr\`{e}s \cite{kollarvanishing}, les faisceaux $R^i\pi_*\mathcal{O}_X$ sont sans torsion donc localement libres sur $\Gamma$. Pour toute fibre $Y$ de $\pi$, on a donc
$H^2(Y,\mathcal{O}_Y)=0$. Ceci entra\^{\i}ne $H^2(Z,\mathcal{O}_Z)=0$ car les singularit\'{e}s quadratiques ordinaires
sont rationnelles.

(ii) Soit $0\in \Gamma(\mathbb{C})$ un point tel que la fibre $X_0=Y$ ait des points doubles ordinaires.
Soit $\Delta\subset \Gamma $ un petit disque analytique centr\'{e} en $0$, et soit $0\not=t\in \Delta$.
On dispose (cf. \cite[14.2.1]{voisinlivre}) d'une application continue
$f:X_t(\mathbb{C})\rightarrow Y(\mathbb{C})$ qui est un hom\'{e}omorphisme  au-dessus de
la partie lisse de $Y(\mathbb{C})$ et qui contracte une sph\`{e}re \'{e}vanescente $S^n\subset X_t(\mathbb{C})$
sur chaque point double de $Y$.
Par ailleurs on dispose de l'application de d\'{e}singularisation $g : Z(\mathbb{C}) \rightarrow Y (\mathbb{C})$, qui est aussi
un hom\'{e}omorphisme au-dessus de la partie lisse de $Y (\mathbb{C})$ et qui contracte une quadrique $Q_i$ de
dimension $n-1$  sur chaque point double de $Y$.
Pour $p\not=0$, les faisceaux $R^pf_*\mathbb{Z}$ et $R^pg_*\mathbb{Z}$ sont support\'{e}s sur les points doubles
de $Y$ et on a donc pour $p\not=0$, et $q>0$,
 $$H^q(Y(\mathbb{C}),R^pg_*\mathbb{Z})=0=H^q(Y(\mathbb{C}),R^pf_*\mathbb{Z}).$$
 Au niveau  $E_2$ de la suite spectrale de Leray de $f$, seuls les termes
 $$ H^0(Y(\mathbb{C}), R^3f_*\mathbb{Z})\,\,{\rm et}\,\,H^3(Y(\mathbb{C}), R^0f_*\mathbb{Z})=H^3(Y(\mathbb{C}),\mathbb{Z})$$
 contribuent donc \`{a} $H^3(X_t(\mathbb{C}),\mathbb{Z})$, et de m\^{e}me pour $g$ et $Z$.
  De plus,  on a $R^pf_*\mathbb{Z}=0$ pour $p\not=0,n$, donc en particulier $R^3f_*\mathbb{Z}=0$
pour $n\not=3$.
Enfin, pour $n\not=2,3$, le faisceau
$R^2g_*\mathbb{Z}$ est constitu\'{e} d'un exemplaire de $\mathbb{Z}=H^2(Q_i,\mathbb{Z})$ support\'{e} sur chaque point double de $Y$ et de plus  le g\'{e}n\'{e}rateur $\alpha_i$ de $H^2(Q_i,\mathbb{Z})$ est engendr\'{e} par la classe $c_1(\mathcal{O}_{Z}(Q_i))_{\mid Q_i}$. On en d\'{e}duit que pour $n\not=2,3$, la fl\`{e}che
$$H^2(Z(\mathbb{C}),\mathbb{Z})\rightarrow H^0(Y(\mathbb{C}),R^2g_*\mathbb{Z})$$ est surjective, de sorte que la diff\'{e}rentielle
$$d_3:H^0(Y(\mathbb{C}),R^2g_*\mathbb{Z})\rightarrow H^3(Y(\mathbb{C}),R^0g_*\mathbb{Z})=H^3(Y(\mathbb{C}),\mathbb{Z})$$
est nulle. Cet \'{e}nonc\'{e} est bien s\^{u}r aussi valable pour $f$ puisqu'on a $R^2f_*\mathbb{Z}=0$ pour $n\not=2$.
Pour $n\not=2,3$, les suites spectrales de Leray de $f$ et $g$ fournissent donc des isomorphismes
 $$H^3(Z(\mathbb{C}),\mathbb{Z})\stackrel{g^*}{\loi}H^3(Y(\mathbb{C}),\mathbb{Z})
\stackrel{f^*}{\oi}H^3(X_t(\mathbb{C}),\mathbb{Z}).$$
Dans le cas $n = 2$, la r\'{e}solution simultan\'{e}e montre directement que $Z(\mathbb{C})$ et $X_t(\mathbb{C})$
sont hom\'{e}omorphes, d'o\`{u} la conclusion.
\cqfd

\begin{rema}{\rm Dans le cas $n=3$, les r\'{e}sultats de
\cite{endrass} montrent qu'on peut effectivement avoir une famille de vari\'{e}t\'{e}s de Fano
de dimension $3$ param\'{e}tr\'{e}e par une courbe lisse sur
$\mathbb{C}$, \`{a} fibres au pire nodales, dont les fibres lisses n'ont pas de torsion
dans leur cohomologie de Betti de degr\'{e} $3$, tandis que les fibres singuli\`{e}res
ont une d\'{e}singularisation poss\'{e}dant de la torsion dans leur cohomologie de Betti de degr\'{e} $3$.}
\end{rema}

\begin{prop}\label{singularitesordinaires}
Soit $\pi : X \to \Gamma$ un morphisme dominant de vari\'{e}t\'{e}s projectives, lisses, connexes sur $\C$, avec $\Gamma$ une courbe.
Supposons la fibre g\'{e}n\'{e}rique $V=X_{\eta}$ g\'{e}om\'{e}triquement int\`{e}gre,
 sur le corps $F=\C(\Gamma)$.   Soit $n$ un entier positif.

Supposons que $H^2(V,\O_{V})=0$ et que la cohomologie de Betti
enti\`{e}re de degr\'{e} 3 d'une fibre lisse de $\pi$ est sans torsion.

(i) Soit  $\Gamma^{0} \subset \Gamma$ l'ouvert maximal au-dessus duquel $\pi$
 est lisse et $X^{0}=\pi^{-1}(\Gamma^{0} )$.
 La restriction naturelle $H^3_{nr}(X^0,\mu_{n}^{\otimes 2}) \to H^3_{nr}(V,\mu_{n}^{\otimes 2})$
est un isomorphisme.

Supposons de plus que les fibres singuli\`{e}res de $\pi$  n'ont que des singularit\'{e}s quadratiques ordinaires. Alors :

(ii)   L'inclusion naturelle $H^3_{nr}(X,\mu_{n}^{\otimes 2}) \to  H^3_{nr}(X^0,\mu_{n}^{\otimes 2})=H^3_{nr}(V,\mu_{n}^{\otimes 2})$
 a un conoyau fini.

(iii) Si   la dimension de $V=X_{\eta}$ est diff\'{e}rente de 3,
alors la restriction \`{a} la fibre g\'{e}n\'{e}rique donne un isomorphisme  $H^3_{nr}(X,\mu_{n}^{\otimes 2}) \oi
 H^3_{nr}(V,\mu_{n}^{\otimes 2})$.
 \end{prop}

{\bf D\'{e}monstration.} Les hypoth\`{e}ses d'annulation faites sur
 $H^2(V,\mathcal{O}_{V})$ et sur la torsion
 dans la cohomologie enti\`{e}re de degr\'{e} $3$ d'une fibre lisse restent valables pour {\it toute} fibre $Y$
   de $\pi$ au-dessus de $\Gamma^0$. Soit $\alpha \in H^3_{nr}(V,\mu_{n}^{\otimes 2})$.
Le r\'{e}sidu de $\alpha$ en une fibre lisse $Y$ de $\pi$
appartient \`{a} $H^2_{nr}(Y,\mu_n)$, qui est nul gr\^{a}ce \`{a} ces hypoth\`{e}ses. Ceci montre (i).

 En \'{e}clatant sur $X$ les points singuliers des fibres sp\'{e}ciales, on obtient un
mod\`{e}le $X'$ muni d'un morphisme $\pi' : X' \to \Gamma$ dont chaque fibre
non lisse est de la forme $Z+2\sum H_{i}$, avec $Z$ lisse irr\'{e}ductible (obtenu par \'{e}clatement
de la fibre correspondante $Y$ de $\pi$ en ses points doubles) et
chaque diviseur exceptionnel $H_{i}$ isomorphe \`{a}
$\mathbb{P}^n$, les intersections $H_{i} \cap H_{j}$ \'{e}tant vides pour $i\neq j$,
et   $Q_{i}:=Z\cap H_{i}$ \'{e}tant isomorphe \`{a} une hypersurface quadratique lisse de  $H_{i}=\mathbb{P}^n$.
On a la suite exacte de localisation
$$0 \to H^2_{nr}(H_{i},\mu_n) \to H^2_{nr}(H_{i}\setminus  Q_{i},\mu_n) \to  H^1(Q_{i},\mathbb{Z}/n)$$
avec $H^2_{nr}(H_{i},\mu_n)=H^1(Q_{i},\mathbb{Z}/n)=0$. On a donc $H^2_{nr}(H_{i}\setminus Q_{i},\mu_n)=0$.

En une fibre singuli\`{e}re de $\pi':X' \to \Gamma$, le r\'{e}sidu de $\alpha$ au point
g\'{e}n\'{e}rique de $H_{i}$
appartient \`{a} $H^2_{nr}(H_{i}\setminus Q_{i},\mu_n)$, donc est nul.
On en d\'{e}duit que le r\'{e}sidu de $\alpha$ au point g\'{e}n\'{e}rique de $Z$ est sans r\'{e}sidu le long de
$Q_i$ et donc
appartient \`{a} $H^2_{nr}(Z,\mu_n)$.

  Le (i) du  lemme \ref{9mai}  implique que, sous les hypoth\`{e}ses de la proposition
\ref{singularitesordinaires}, le groupe
$H^2_{nr}(Z,\mu_n)$ est fini. Avec les arguments pr\'{e}c\'{e}dents, on conclut que $H^3_{nr}(X',\mu_{n}^{\otimes 2})\rightarrow H^3_{nr}(X^0,\mu_{n}^{\otimes 2})$ a un conoyau de torsion, ce qui donne (ii).

Enfin le (ii) du lemme implique que, sous les hypoth\`{e}ses de la proposition
\ref{singularitesordinaires},
 si $n\not=3$, on a $H^2_{nr}(Z,\mu_n)=0$, et donc le r\'{e}sidu de $\alpha$ au point g\'{e}n\'{e}rique de $Z$ est
nul, d'o\`{u} le point (iii).
\cqfd

\bigskip

 \begin{theo}\label{Hitrivialbis}
 Soit $\pi : X \to \Gamma$ un morphisme dominant de vari\'{e}t\'{e}s projectives, lisses, connexes sur $\C$,
de fibre g\'{e}n\'{e}rique $V=X_{\eta}$ g\'{e}om\'{e}triquement int\`{e}gre, de dimension 2,
 sur le corps $F=\C(\Gamma)$.

Supposons $H^{i}(V,\O_{V})=0$ pour $i=1,2$ et
supposons que la cohomologie de Betti
enti\`{e}re de degr\'{e} 3 d'une fibre lisse de $\pi$ est sans torsion.
 Supposons de plus que les fibres singuli\`{e}res de $\pi$  n'ont que des singularit\'{e}s quadratiques ordinaires.

Si $H^3_{nr}(X,\Q/\Z(2))=0$, alors $I(X_{\eta})=1$.
 \end{theo}

{ \bf D\'{e}monstration.} Ceci r\'{e}sulte de la combinaison de la proposition \ref{singularitesordinaires} et
de la proposition \ref{indicesurface}. \cqfd

\begin{rema}
{\rm On notera que cet  \'{e}nonc\'{e} est plus faible que le r\'{e}sultat obtenu par  application du th\'{e}or\`{e}me
\ref{propindice1}.}
\end{rema}

\begin{theo}\label{corconclus8}
Soit $\pi : X \to \Gamma$ un morphisme surjectif de vari\'{e}t\'{e}s  connexes, projectives et lisses sur $\C$,
 o\`u $\Gamma$ est une courbe et la fibre g\'{e}n\'{e}rique $V=X_{\eta}$ est une surface
  g\'{e}om\'{e}triquement connexe sur $F=\C(\Gamma)$.
  Supposons $H^{i}(V,\O_{V})=0$ pour $i=1,2$, que le groupe de
 N\'{e}ron-Severi de la fibre
  g\'{e}n\'{e}rique g\'{e}om\'{e}trique n'a pas de torsion, et que la fl\`{e}che
  $\deg : CH_{0}(\ovV) \to \Z$ est un isomorphisme.

  (i) On a alors des inclusions $H^3_{nr}(X,\Q/\Z(2)) \subset H^3_{nr}(V,\Q/\Z(2)) \simeq \Z/I(V).$

  (ii) Si les fibres singuli\`{e}res de $\pi$ n'ont que des singularit\'{e}s ordinaires, on a des isomorphismes
  $H^3_{nr}(X,\Q/\Z(2)) = H^3_{nr}(V,\Q/\Z(2)) \simeq \Z/I(V).$
 \end{theo}

 {\bf D\'{e}monstration.}  L'\'{e}nonc\'{e} (i) est une cons\'{e}quence imm\'{e}diate de la proposition
 \ref{indicesurfacesviaKtheorie}. L'\'{e}galit\'{e}  $H^3_{nr}(X,\Q/\Z(2)) = H^3_{nr}(V,\Q/\Z(2))$ dans (ii)
a \'{e}t\'{e} \'{e}tablie \`{a} la proposition  \ref{singularitesordinaires}. \cqfd

\begin{rema} {\rm Comparons ce r\'{e}sultat  avec  le th\'{e}or\`{e}me \ref{propindice2}.
\`A la diff\'{e}rence de ce th\'{e}or\`{e}me,
on a suppos\'{e} ici $CH_{0}(\ovV) \oi \Z$, ce qui selon une conjecture de Bloch
devrait r\'{e}sulter de l'hypoth\`{e}se $H^{i}(V,\O_{V})=0$ pour $i=1,2$.
Cette conjecture est connue pour les surfaces qui ne sont pas de type g\'{e}n\'{e}ral,
et aussi pour un certain nombre de surfaces de type g\'{e}n\'{e}ral.
L'hypoth\`{e}se $CH_{0}(\ovV) \oi \Z$ implique que  $H^3_{nr}(V,\Q/\Z(2))$
est d'exposant fini, et  il en est donc de m\^{e}me de   $H^3_{nr}(X,\Q/\Z(2))$,
qui co\"incide donc
avec $Z^4(X)$.  L'\'{e}nonc\'{e} (ii), o\`u l'on suppose les singularit\'{e}s des fibres le plus simples possibles,
n'est donc pas plus fort que le  th\'{e}or\`{e}me   \ref{propindice2}. Mais   l'\'{e}nonc\'{e} (i) donne en particulier que
  $I(X_{\eta})=1$ implique $H^3_{nr}(X,\Q/\Z(2)=0$ et donc $Z^4(X)=0$, et ceci est \'{e}tabli sans hypoth\`{e}se sur
  les fibres singuli\`{e}res, alors qu'au th\'{e}or\`{e}me  \ref{propindice2} on n'autorise que
  des   singularit\'{e}s quadratiques dans les fibres.}
\end{rema}

\section{Appendice: Action des correspondances sur les groupes $H^i(\mathcal{H}^j)$\label{secappendice}}

Dans ce qui suit $X$ et $Y$ sont lisses et projectives sur $\C$.
On consid\`{e}re les cycles  $Z\subset X\times Y$
de codimension $r+\dim\,X$, c'est-\`{a}-dire tels que  $\dim\,Y-\dim\,Z=r$.
\begin{prop} \label{theoappendice} Pour tout groupe ab\'{e}lien $A$, le  cycle $Z\in CH^{r+\dim\,X}(X\times Y)/alg$ induit un homomorphisme
$$Z_*:H^p(X,\mathcal{H}_X^q(A))\rightarrow H^{p+r}(Y,\mathcal{H}^{q+r}_Y(A(r))).$$
Ces actions sont compatibles \`{a} la composition des
correspondances.
\end{prop}
{\bf D\'{e}monstration.}
Tout d'abord, par la formule de Bloch--Ogus (cf. Corollaire \ref{formuledebo}),  $Z$ a une classe $$[Z]_{BO}\in
H^{d+r}(X\times Y,\mathcal{H}_{X\times Y}^{d+r}(\mathbb{Z}(d+r))),$$ o\`{u} $d=\dim\,X$. Par ailleurs
on a clairement des fl\`{e}ches
$$pr_1^*:
H^p(X,\mathcal{H}^q(A))\rightarrow H^p(X\times Y,\mathcal{H}_{X\times Y}^q(A)),$$ et enfin on dispose
de produits (induits par les cup-produits $\mathcal{H}_W^l(A)\otimes \mathcal{H}^s_W(\mathbb{Z}(l))\rightarrow \mathcal{H}_W^{l+s}(A(l))$):
$$H^p(W,\mathcal{H}_W^l(A))\otimes H^q(W,\mathcal{H}_W^s(\mathbb{Z}(l)))\rightarrow H^{p+q}(W,\mathcal{H}_W^{l+s}(A(l)))$$
 pour toute vari\'{e}t\'{e} $W$ et tout entier $l$. On en d\'{e}duit que
l'on a une fl\`{e}che
$$[Z]_{BO}\circ pr_1^*:H^p(X,\mathcal{H}_X^q(A))\rightarrow H^{p+d+r}(X\times Y,\mathcal{H}_{X\times Y}^{q+d+r}(A(d+r))).$$
Pour conclure, il suffit donc de construire
$$pr_{2*}:H^i(X\times Y,\mathcal{H}_{X\times Y}^j(A(d+r)))\rightarrow H^{i-d}( Y,\mathcal{H}^{j-d}_Y(A(r)))$$
pour $X$ projective lisse de dimension $d$. (On peut \'{e}videmment remplacer ici la projection par n'importe quelle application propre.)

Or ceci r\'{e}sulte de la r\'{e}solution de Bloch--Ogus (th\'{e}or\`{e}me \ref{boresolution})
qui donne la formule explicite :
\begin{eqnarray}\nonumber\label{groupe}H^i(X\times Y,\mathcal{H}_{X\times Y}^j(A(d+r)))=\hskip10cm\\
\frac{\Ker\,[\oplus_{x\in (X\times Y)^{(i)}}H^{j-i}_{B}(\mathbb{C}(x),A(d+r-i))\stackrel{\partial}{\rightarrow}
\oplus_{x\in (X\times Y)^{(i+1)}}H^{j-i-1}_{B}(\mathbb{C}(x),A(d+r-i-1))]}
{\Im\,[\oplus_{x\in (X\times Y)^{(i-1)}}H^{j-i+1}_{B}(\mathbb{C}(x),A(d+r-i+1))
\stackrel{\partial}{\rightarrow}\oplus_{x\in (X\times Y)^{(i)}}H^{j-i}_{B}(\mathbb{C}(x),A(d+r-i))]}.
\end{eqnarray}
Un point $x$ de codimension $i$ a une adh\'{e}rence $Z_x\subset X\times Y$
qui s'envoie sur une sous-vari\'{e}t\'{e} $Z'_x\in Y$ de codimension $\geq i-d$, avec \'{e}galit\'{e} lorsque
$pr_2:Z_x\rightarrow Z'_x$ est g\'{e}n\'{e}riquement finie. On notera $x'$ le point de $Y$ correspondant \`{a}
$Z'_x$.
Pour un $x$ ne satisfaisant pas la  condition que $x'\in Y^{(i-d)}$, on dira que
$pr_{2*}:H^{j-i}_{B}(\mathbb{C}(x),A(d+r-i))\rightarrow H^{j-i}_{B}(\mathbb{C}(x'),A(d+r-i))$ est nulle, et
dans le cas contraire,
$pr_{2*}:H^{j-i}_{B}(\mathbb{C}(x),A(d+r-i))\rightarrow H^{j-i}_{B}(\mathbb{C}(x'),A(d+r-i))$ sera induite par
l'application  propre $pr_{2\mid Z_x}: pr_{2\mid Z_x}^{-1}(V)\rightarrow V$ pour des ouverts suffisamment petits  $V$  de $Z'_x$.
Nous avons alors besoin du lemme suivant:
\begin{lemm} \label{lemmeappendice1} Les fl\`{e}ches $pr_{2*} $ ci-dessus commutent aux fl\`{e}ches de r\'{e}sidus, c'est-\`{a}-dire les diff\'{e}rentielles $\partial$ apparaissant dans (\ref{groupe}).
\end{lemm}
Ce lemme entra\^{\i}ne que la somme des fl\`{e}ches  $pr_{2*}$ passe \`{a} la cohomologie pour donner
une fl\`{e}che $pr_{2*}$ du groupe
(\ref{groupe}) dans
$$\frac{\Ker\,[\oplus_{x\in  Y^{(i-d)}}H^{j-i}_{B}(\mathbb{C}(x),A(d+r-i))\stackrel{\partial}{\rightarrow}\oplus_{x\in Y^{(i+1-d)}}H^{j-i-1}_{B}(\mathbb{C}(x),A(d+r-i-1))]}{\Im\,[\oplus_{x\in Y^{(i-d-1)}}H^{j-i+1}_{B}(\mathbb{C}(x),A(d+r-i+1))\stackrel{\partial}{\rightarrow}\oplus_{x\in  Y^{(i-d)}}H^{j-i}_{B}(\mathbb{C}(x),A(d+r-i))]}
$$
$$=H^{i-d}(Y,\mathcal{H}_Y^{j-d}(A(r)))
,$$
o\`{u} la derni\`{e}re \'{e}galit\'{e} est encore donn\'{e}e par le th\'{e}or\`{e}me de r\'{e}solution \ref{boresolution}, appliqu\'{e} \`{a} $Y$ cette fois.

Pour conclure la preuve de la proposition \ref{theoappendice}, il reste \`{a} \'{e}tablir le lemme \ref{lemmeappendice1} ainsi qu'\`{a} montrer le lemme suivant:
\begin{lemm}\label{lemmeappendice2} L'action ainsi construite $Z\mapsto Z_*$, de $CH^{\dim\,X+r}(X\times Y)/alg$
dans $$Hom\,(H^p(X,\mathcal{H}_X^q(A)), H^{p+r}(Y,\mathcal{H}^{q+r}_Y(A(r))) ),$$ est compatible \`{a} la composition des correspondances.
\end{lemm}

\cqfd

{\bf D\'{e}monstration du lemme \ref{lemmeappendice1}.} Soit $Z\subset X\times Y$ une sous-vari\'{e}t\'{e} de codimension
$i$. Soit $D'\subset Z'=pr_2(Z)$ une sous-vari\'{e}t\'{e} de codimension $i-d-1$ dans $Y$. Supposons pour simplifier que $Z$ et $Z'$ sont normales.
On veut montrer que pour tout entier $l$, le compos\'{e}
$$Res_{Z',D'}\circ pr_{2*}: H^k_{B}(\mathbb{C}(Z),A(l))\rightarrow H^k_{B}(\mathbb{C}(Z'),A(l))\rightarrow H^{k-1}_{B}(\mathbb{C}(D'),A(l-1))$$
est \'{e}gal \`{a}
$pr_{2*}\circ \sum_{\overset{D\subset Z}{pr_2(D)=D'}} pr_{2*}\circ  Res_{Z,D}$: $$H^k_{B}(\mathbb{C}(Z),A(l))\rightarrow
\bigoplus_{\overset{D\subset Z}{pr_2(D)=D'}}H^{k-1}_{B}(\mathbb{C}(D),A(l-1))\rightarrow H^{k-1}_{B}(\mathbb{C}(D'),A(l-1)).$$
Supposons tout d'abord que $\dim\,Z'=\dim\,Z$. Alors, rappelant que les fl\`{e}ches de r\'{e}sidus sont induites par la suite exacte de cohomologie \`{a}  supports pour les paires constitu\'{e}es
d'une vari\'{e}t\'{e} et d'un diviseur, le r\'{e}sultat r\'{e}sulte simplement du fait qu'on a un morphisme
propre $pr_2$ entre les  paires $({Z}_{0},\cup_{ pr_2(D)=D'}D)$ et $(Z'_{0},D')$, o\`{u} $Z'_{0}$ est un ouvert  lisse de $Z'$ contenant un ouvert de Zariski dense de $D'$, et ${Z}_{0}$ est l'image inverse de $Z'_{0}$ dans $Z$. De plus ce morphisme reste propre lorsqu'on passe aux compl\'{e}mentaires des diviseurs consid\'{e}r\'{e}s.
Sans hypoth\`{e}se de normalit\'{e}, il suffit de normaliser les vari\'{e}t\'{e}s consid\'{e}r\'{e}es.

Il reste \`{a} \'{e}tudier ce qui se passe lorsque
$\dim\,Z'<\dim\,Z$, ce qui n'est possible que si  $\dim\,Z'=\dim\,Z-1$, et $Z'=D'$. Dans ce cas, le premier compos\'{e} est nul par d\'{e}finition car $pr_{2*}: H^k_{B}(\mathbb{C}(Z),A(l))\rightarrow H^k_{B}(\mathbb{C}(Z'),A(l))$ est nul. Il suffit donc de montrer que pour toute classe $\alpha\in H^k_{B}(\mathbb{C}(Z),A(l))$, on a
$pr_{2*}(Res\,\alpha)=0$, o\`u $Res\,\alpha$ est une somme finie de classes de degr\'{e} $k-1$ sur des
diviseurs $D\subset Z$ tels que $pr_2:D\rightarrow D'$ soit g\'{e}n\'{e}riquement fini.
Ceci est \'{e}l\'{e}mentaire car l'application ${pr_2}:Z\rightarrow Z'$
induit une application $pr_{2}^0: Z_0\rightarrow Z_0'$ qui est propre et lisse  de dimension relative
$1$ au-dessus d'un ouvert de Zariski dense $Z'_0$ de $D'=Z'$. On dispose donc d'une application
$pr_{2*}^0: H^{k+1}(Z_0,A(l))\rightarrow H^{k-1}(Z'_0,A(l-1))$ et on a alors la relation suivante, pour
$\beta\in H^{k-1}(D,A(l-1))$, o\`u $D\subset Z^0$ est un diviseur lisse propre au-dessus de $Z'_0$:
$$pr_{2*}(\beta)=pr_{2*}^0\circ j_*(\beta),$$
o\`u le premier $pr_{2*}$ est relatif \`{a} l'application propre induite par $pr_2$ de $D$ sur $Z'_0$.
Cela conclut la preuve car on a
$$j_*\circ Res=0: H^k(Z^0\setminus D,A(l))\rightarrow H^{k-1}(D,A(l-1))\rightarrow H^{k+1}(Z^0,A(l)).$$

\cqfd

{\bf D\'{e}monstration du lemme \ref{lemmeappendice2}.} En examinant la preuve
de la compatibilit\'{e} de l'action des correspondances sur les groupes de Chow ou la cohomologie avec la composition des correspondances (cf. \cite[Prop. 21.17]{voisinlivre}), et la construction de l'action
$Z_*$ sur les $H^p(\mathcal{H}^q)$, on constate que la compatibilit\'{e} de l'action $ Z_*$ construite ci-dessus avec  la composition des correspondances est une cons\'{e}quence formelle des trois faits suivants, faciles \`{a} v\'{e}rifier dans notre cas:

(1) La compatibilit\'{e} de l'application classe de cycles
de Bloch--Ogus $Z\mapsto [Z]_{BO}$ avec le produit d'intersection.

(2)  La formule de
projection pour les applications $pr_{2*}:H^k_{B}(\mathbb{C}(Z),A(l))\rightarrow H^{k}_{B}(\mathbb{C}(Z'),A(l))$ introduites ci-dessus.
 Pour $\alpha\in H^p_{B}(\mathbb{C}(Z),A),\,\beta\in H^p_{B}(\mathbb{C}(Z'),\mathbb{Z}(l))$
$$pr_{2*}(\alpha\cup pr_2^*\beta)=pr_{2*}(\alpha)\cup \beta.$$

(3)  La formule de changement de base suivante : prenons un produit $X\times Y\times Z$ de vari\'{e}t\'{e}s projectives avec $\dim\,X=d$, et notons $p$ la projection
de $X\times Y$ sur $Y$, $P$ la projection
de $X\times Y\times Z$ sur $Y\times Z$, $q$ la projection de
$Y\times Z$ sur $Y$ et $Q$ celle  de $X\times Y\times Z$ sur $X\times Y$.
Alors on a pour $\alpha\in H^p(X\times Y,\mathcal{H}^q_{X\times Y}(A(l)))$:
$$q^*(p_*\alpha)=P_*(Q^*\alpha)\,\,{\rm dans}\,\,
H^{p-d}(X\times Y,\mathcal{H}^{q-d}_{Y\times Z}(A(l-d))).$$
\cqfd

\bigskip

{\it Remerciements.} Ce travail a \'{e}t\'{e} commenc\'{e} lors de l'atelier \og Rational curves and $\A^1$-homotopy theory \fg \, organis\'{e}
par Aravind Asok et Jason Starr \`{a}   l'American Institute of Mathematics (Palo Alto, Californie)  du 5 au 9 octobre 2009.}

\vfill\eject

Jean-Louis  Colliot-Th\'{e}l\`{e}ne

C.N.R.S., UMR  8628,

Math\'{e}matiques, B\^atiment 425,

Universit\'{e}   Paris-Sud,

F-91405 Orsay, France

 \smallskip

courriel :  jlct@math.u-psud.fr

 \medskip

 Claire Voisin,

 C.N.R.S.,
Institut de math\'{e}matiques de Jussieu,

Case 247,

4 Place Jussieu,

F-75005 Paris, France

\smallskip

courriel : voisin@math.jussieu.fr

\end{document}